\documentclass{amsart}

\usepackage{latexsym}

\usepackage{amsfonts,amssymb,euscript,epsfig}

\usepackage[all]{xy}

\newtheorem{theo}{Theorem}[section]
\newtheorem{defi}[theo]{Definition}

\newtheorem{lem}[theo]{Lemma}
\newtheorem{prop}[theo]{Proposition}
\newtheorem{rem}[theo]{Remark}
\newtheorem{coro}[theo]{Corollary}

\newcommand{\agot}{\ensuremath{\mathfrak{a}}}
\newcommand{\ggot}{\ensuremath{\mathfrak{g}}}
\newcommand{\hgot}{\ensuremath{\mathfrak{h}}}
\newcommand{\kgot}{\ensuremath{\mathfrak{k}}}

\newcommand{\pgot}{\ensuremath{\mathfrak{p}}}
\newcommand{\tgot}{\ensuremath{\mathfrak{t}}}
\newcommand{\Rgot}{\ensuremath{\mathfrak{R}}}

\newcommand{\Acal}{\ensuremath{\mathcal{A}}}
\newcommand{\Bcal}{\ensuremath{\mathcal{B}}}
\newcommand{\Dcal}{\ensuremath{\mathcal{D}}}
\newcommand{\Ecal}{\ensuremath{\mathcal{E}}}

\newcommand{\Hcal}{\ensuremath{\mathcal{H}}}
\newcommand{\Kcal}{\ensuremath{\mathcal{K}}}
\newcommand{\Lcal}{\ensuremath{\mathcal{L}}}
\newcommand{\Mcal}{\ensuremath{\mathcal{M}}}
\newcommand{\Ncal}{\ensuremath{\mathcal{N}}}
\newcommand{\Ocal}{\ensuremath{\mathcal{O}}}
\newcommand{\Qcal}{\ensuremath{\mathcal{Q}}}
\newcommand{\Scal}{\ensuremath{\mathcal{S}}}
\newcommand{\Xcal}{\ensuremath{\mathcal{X}}}
\newcommand{\Ycal}{\ensuremath{\mathcal{Y}}}
\newcommand{\Zcal}{\ensuremath{\mathcal{Z}}}
\newcommand{\Ucal}{\ensuremath{\mathcal{U}}}
\newcommand{\Vcal}{\ensuremath{\mathcal{V}}}

\newcommand{\Z}{\ensuremath{\mathbb{Z}}}
\newcommand{\C}{\ensuremath{\mathbb{C}}}
\newcommand{\Q}{\ensuremath{\mathbb{Q}}}
\newcommand{\R}{\ensuremath{\mathbb{R}}}
\newcommand{\N}{\ensuremath{\mathbb{N}}}
\newcommand{\tore}{\ensuremath{\mathbb{T}}}
\newcommand{\Lfibre}{\ensuremath{\mathbb{L}}}

\newcommand{\esp}{\ensuremath{\varepsilon}}
\newcommand{\f}{\ensuremath{\mathcal{C}^{\infty}}}
\newcommand{\fgene}{\ensuremath{\mathcal{C}^{-\infty}}}
\newcommand{\croc}{\ensuremath{\hookrightarrow}}

\newcommand{\indice}{\ensuremath{\hbox{\rm Index}}}

\newcommand{\Cr}{\ensuremath{\hbox{\rm Cr}}}
\newcommand{\Vol}{\ensuremath{\hbox{\rm vol}}}
\newcommand{\T}{\ensuremath{\hbox{\bf T}}}

\newcommand{\Char}{\ensuremath{\hbox{\rm Char}}}
\newcommand{\End}{\ensuremath{\hbox{\rm End}}}
\newcommand{\Hom}{\ensuremath{\hbox{\rm Hom}}}
\newcommand{\Thom}{\ensuremath{\hbox{\rm Thom}}}
\newcommand{\sthom}{\ensuremath{\hbox{\rm S-Thom}}}

\def \indH {{\rm Ind}^{^G}_{_H}}
\def \indB {{\rm Ind}^{^G}_{_{G_{\beta}}}}
\def \indT {{\rm Ind}^{^G}_{_{G_{\tau}}}}
\def \inds {{\rm Ind}^{^G}_{_{G_{\sigma}}}}

\def \carre {{\rm L}^{2}}
\def \HolH {{\rm Hol}^{^G}_{_H}}
\def \HolB {{\rm Hol}^{^G}_{_{G_{\beta}}}}
\def \HolHB {{\rm Hol}^{^{G_{\beta}}}_{_H}}
\def \Hols {{\rm Hol}^{^G}_{_{G_{\sigma}}}}
\def \HolHs {{\rm Hol}^{^{G_{\sigma}}}_{_H}}
\def \HolS {{\rm Hol}^{^G}_{_{S^1}}}

\def \spin {{\rm Spin}}
\def \spinc {{\rm Spin}^{\rm c}}
\def \so {{\rm SO}}
\def \u {{\rm U}}
\def \Det {{\rm Det}}
\def \Pspin {{\rm P}_{{\rm Spin}^{\rm c}}}
\def \Pso {{\rm P}_{{\rm SO}}}
\def \Pu {{\rm P}_{{\rm U}}}
\def \clif {{\bf c}}
\def \Clif {{\rm Cl}}

\begin{document}
    
\title{Localization of the Riemann-Roch character}

\author{Paul-Emile  PARADAN}

\date{\today}

\maketitle

{\center\small
Universit\'e Grenoble I, Institut Fourier\\
B.P. 74, 38402, Saint-Martin-d'H\`eres\\
e-mail: Paul-Emile.Paradan@ujf-grenoble.fr\\
}

\begin{abstract}
We present a K-theoretic approach to the Guillemin-Sternberg 
conjecture \cite{Guillemin-Sternberg82}, about the commutativity 
of geometric 
quantization and symplectic reduction, which was proved  
by Meinrenken \cite{Meinrenken96,Meinrenken98} and Tian-Zhang
\cite{Tian-Zhang98}. Besides providing a new proof of this 
conjecture for the full non-abelian group action case, our methods 
lead to a generalization for compact Lie group actions on manifolds 
that are not symplectic; these manifolds carry an invariant 
almost complex structure and an abstract moment map.
\end{abstract}


{\def\thefootnote{\relax}
\addtocounter{footnote}{-1}
}

{\small
\tableofcontents
}

\section{Introduction}

This article is devoted to the study of the `quantization commutes 
with reduction' principle of Guillemin-Sternberg 
\cite{Guillemin-Sternberg82}. The object of this paper is twofold. 
The first goal is to give a K-theoretic approach to this problem 
which provides a new 
proof of results obtained by Meinrenken \cite{Meinrenken98},
Meinrenken-Sjamaar \cite{Meinrenken-Sjamaar} and Tian-Zhang
\cite{Tian-Zhang98}. The second goal is to define an extension 
to the {\em non-symplectic} case. 

In the Kostant-Souriau framework 
one considers a prequantum line bundle $L$ over a compact 
symplectic manifold $(M,\omega)$ :  $L$ carries a 
Hermitian connection $\nabla^{L}$ 
with curvature form equal to $-\imath\omega$. Suppose now that a 
compact Lie group $G$, with Lie algebra $\ggot$, 
acts on $L\to M$, living the data 
$(\omega,\nabla^{L})$ invariant. Then the $G$-action on $(M,\omega)$ 
is Hamiltonian with moment map $f_{_{G}}:M\to\ggot^*$ given by 
the Kostant formula : 
$\Lcal^L(X)-\nabla^L_{X_M}=\imath \langle f_{_G},X\rangle, \ X\in\ggot$. 
Here $\Lcal^L(X)$ is the infinitesimal action of $X$ on the section 
of $L\to M$ and  $X_M$ is the vector field on $M$ generated by 
$X\in \ggot$. 

Choose now an invariant almost complex structure $J$ on $M$ 
that is compatible with $\omega$, in the sense that $\omega(-,J-)$ 
defines a Riemannian metric. It defines a quantization map
$$
RR^{^{G,J}}(M,-): K_G(M)\to R(G)\ ,
$$
from the equivariant $K$-theory of complex vector bundles over $M$ to 
the character ring of $G$. The 
`quantization commutes with reduction' Theorem tells us how the 
multiplicities of $RR^{^{G,J}}(M,L)$ behave (see Theorem {\bf C}).

\medskip 

Here our main goal is to compute the multiplicity of the trivial 
representation in $RR^{^{G,J}}(M,L)$, when the data $(L,J)$ are 
not associated to a symplectic form. 

\medskip

We consider a compact manifold $M$ on which a compact Lie group $G$ 
acts, and which carries a $G$-invariant almost complex structure $J$.
Let $L\to M$ be a $G$-equivariant Hermitian line bundle over $M$, 
equipped with a Hermitian connection $\nabla^L$ on $L$.  
This defines a map $f_{_L}:M\to \ggot^*$ by the equation 
\begin{equation}\label{eq.moment.L}
\Lcal^L(X)-\nabla^L_{X_M}=\imath \langle f_{_L},X\rangle,\quad 
X\in\ggot\ .
\end{equation}
(see \cite{B-G-V}[section 7.1]).  The map $f_{_L}$ is an {\em abstract 
moment map} in the sense of Karshon \cite{Karshon.98}, since $f_{_L}$ 
is equivariant, and for any $X\in\ggot$, the function $\langle 
f_{_L},X\rangle$ is locally constant on the submanifold $M^X:=\{m\in 
M,\ X_{M}(m)=0\}$.

If $0$ is a regular value of $f_{_L}$, $\Zcal:= f_{_L}^{-1}(0)$ is a 
smooth submanifold of $M$ which carries a locally free action of $G$.  
We consider the orbifold reduced space $\Mcal_{red}=\Zcal/G$ and we 
denote $\pi :\Zcal\to\Mcal_{red}$ the projection.  In 
Lemma \ref{lem.spinc.induit} we show that the almost complex structure 
$J$ induces an orientation $o_{red}$ on $\Mcal_{red}$ and a Spin$^{\rm 
c}$ structure on $(\Mcal_{red},o_{red})$.  Let $\Qcal(\Mcal_{red},-): 
K_{orb}(\Mcal_{red})\to \Z$ be the quantization map defined by the 
Spin$^{\rm c}$ structure and let $L_{red}\to\Mcal_{red}$ be the 
orbifold line bundle induced by $L$.


We obtain the following `quantization commutes with reduction' 
theorem.

\medskip

{\bf Theorem A}\hspace{5mm}{\it Let $L\to M$ be a $G$-equivariant 
Hermitian line bundle over $M$, equipped with a Hermitian connection 
$\nabla^L$ on $L$. Let $f_{_L}:M\to \ggot^*$ be the corresponding 
abstract moment map. If $0$ is a regular value of 
$f_{_{L}}$, we have 
\begin{equation}\label{eq.theorem.A}
 \left[ RR^{^{G,J}}(M,L^{\stackrel{k}{\otimes}})\right]^G=
 \Qcal\Big(\Mcal_{red},L_{red}^{\stackrel{k}{\otimes}}\Big), \quad 
 k\in\N-\{0\},
\end{equation}
if any of the following hold:

(i) $G=T$ is a torus; or

(ii) $k\in \N$ is large enough , so that the ball 
$\{ \xi\in\ggot^{*},\, \parallel\xi\parallel\leq\frac{1}{k}
\parallel\theta\parallel \}$
is contained in the set of regular values of $f_{_L}$. 
Here $\theta=\sum_{\alpha>0}\alpha$ is 
the sum of the positive roots of $G$, and $\parallel\cdot\parallel$ 
is a $G$-invariant Euclidean norm on $\ggot^*$.

Here, for $V\in R(G)$, we denote $[V]^G\in \Z$ the multiplicity of the 
trivial representation.
}

\medskip

A similar result was proved by Jeffrey-Kirwan \cite{Jeffrey-Kirwan97} 
in the Hamiltonian setting when one relaxes the condition of
positivity of $J$ with respect to the symplectic form. See also  
\cite{CdS-K-T} for a similar result in the Spin$^{\rm c}$ setting,
when $G=S^1$.

\medskip  

 As an example, let us apply Theorem {\bf A} to the counterexemple 
due to Vergne which shows that quantization does not always commute 
with reduction. Let $G=SU(2)$ and let $M$ be the $SU(2)$-coadjoint 
orbit passing through the
unique positive root $\theta$. Thus $M$ is the projective line bundle
$\C\mathbb{P}^{1}$ with $\omega$ equal to twice the standard
K\"ahler form. The prequantum line bundle is $L=\mathcal{O}(2)$ and
$RR^{^G}(M,L^{-1})=[RR^{^G}(M,L^{-1})]^G=-1$. Since
$\Mcal_{red}=\emptyset$ we have $[RR^{^G}(M,L^{-1})]^G\neq$ \break 
$\Qcal(\Mcal_{red},(L^{-1})_{red})$ : the condition $ii)$ of
Theorem {\bf A} does not hold since $\theta$ is not a regular value of
the moment map $M\croc\ggot^*$. But if we take
$(L^{-1})^{\stackrel{k}{\otimes}}$ with $k>1$ the condition $ii)$ is
satisfied, and thus $[RR^{^G}(M,(L^{-1})^{\stackrel{k}{\otimes}})]^G
=0$ for $k>1$. In fact a direct computation shows that
$-\,RR^{^G}(M,(L^{-1})^{\stackrel{k}{\otimes}})$ is the character of the
irreducible $SU(2)$-representation with highest weight $(k-1)\theta$ for
all $k\geq 1$.

\medskip

The result of  Theorem {\bf A} can be rewritten when $J$ defines an almost 
complex structure $J_{red}$ on $\Mcal_{red}$. It happens when  
the following decomposition holds 
\begin{equation}\label{eq.J.red}
\T M\vert_{\Zcal}=\T\Zcal\oplus J(\ggot_{\Zcal})\quad {\rm with}\quad 
\ggot_{\Zcal}:=\{X_{\Zcal},\, X\in\ggot\}\ .
\end{equation}
First we note that (\ref{eq.J.red}) always holds in the 
Hamiltonian case when $J$ is compatible with the symplectic form. 
Condition (\ref{eq.J.red}) already appears in the works of 
Jeffrey-Kirwan \cite{Jeffrey-Kirwan97}, and Cannas da Silva-Karshon-Tolman 
\cite{CdS-K-T}.

In all this paper we fix a $G$-invariant 
scalar product on $\ggot^*$ which induces an identification 
$\ggot\simeq\ggot^*$. Thus $f_{_G}$ can be considered as a map from
$M$ to $\ggot$, and we define the 
endomorphism $\Dcal$ of the bundle $\ggot\times\Zcal$  by : 
$\Dcal(X)=-d\,f_{_G}(J(X_{\Zcal}))$, for $X\in\ggot$. Condition 
(\ref{eq.J.red}) is then equivalent to : 
$\det\Dcal(z)\neq 0$ for all $z\in\Zcal$. The endomorphism $\Dcal$
defines a complex stucture $J_{\Dcal}$ on $\Zcal\times\ggot_{\C}$, 
so the vector bundle $\Zcal\times\ggot_{\C}$ inherits two irreducible 
complex spinor bundles  $\Zcal\times\wedge^{\bullet}_{\C}\ggot_{\C}$ 
and $\Zcal\times\wedge^{\bullet}_{J_{\Dcal}}\ggot_{\C}$ related by
$$
\wedge^{\bullet}_{J_{\Dcal}}\ggot_{\C}\times\Zcal = 
\wedge^{\bullet}_{\C}\ggot_{\C}\times\Zcal\otimes \pi^*L_{\Dcal}
$$
where $L_{\Dcal}\to\Mcal_{red}$ is a line bundle (see (\ref{eq.L.D})). 
In this case we prove in Proposition \ref{prop.RR.T.0.bis} 
that (\ref{eq.theorem.A}) has the following form
\begin{equation}\label{eq.theorem.A.bis}
 \left[ RR^{^{G,J}}(M,L^{\stackrel{k}{\otimes}})\right]^G
 =\pm RR^{^{J_{red}}}
 \Big(\Mcal_{red},L_{red}^{\stackrel{k}{\otimes}}\otimes 
 L_{\Dcal}\Big)\ , 
\end{equation}
where $\pm$ is the sign of $\det\Dcal$, and where 
$RR^{^{J_{red}}}(\Mcal_{red},-)$ is the Riemann-Roch character 
defined by $J_{red}$. 

\medskip

In this paper, we start from an abstract moment 
map $f_{_G}:M\to\ggot^*$ , and we extend the result of Theorem {\bf A} 
to the $f_{_G}$-moment bundles, and the 
$f_{_G}$-positive bundles. These notions were
introduced in the Hamiltonian setting by Meinrenken-Sjamaar 
 \cite{Meinrenken-Sjamaar} and Tian-Zhang \cite{Tian-Zhang98}.
Let us recall the definitions.

Let $H$ be a maximal torus of $G$ with Lie algebra $\hgot$.
 
\begin{defi}\label{moment.bundle}
A $G$-equivariant line bundle over $M$ is called a $f_{_G}$-moment
bundle if for all components $F$ of the fixed-point set $M^H$ 
the weight of the $H$-action on $L\vert_F$ is equal to $f_{_G}(F)$. 
\end{defi}

It is easy to see that the definition is independent of the choice of
the maximal torus. Note that $f_{_G}(F)\in \hgot^*=(\ggot^*)^H$,
since $f_{_G}$ is equivariant.  Any Hermitian line bundle $L$ is 
tautologically a moment bundle relative to the abstract moment 
map $f_{_L}$.

For any $\beta\in\ggot$, we denote by $\tore_{\beta}$ 
the torus of $G$ generated by $\exp_{G}(t.\beta),\, t\in \R$, and 
$M^{\beta}$ the submanifold of points fixed by $\tore_{\beta}$.

\begin{defi}\label{eq.mu.positif}
A complex $G$-vector bundle $E$ is called $f_{_G}$-{\em positive} 
if the following hold: for 
any $m\in M^{\beta}\cap f_{_G}^{-1}(\beta)$, we have
$$
  \langle \xi,\beta\rangle \geq 0
$$
for any weights $\xi$ of the $\tore_{\beta}$-action on $E_m$. 
A complex $G$-vector bundle $E$ is called $f_{_{G}}$-{\em strictly 
positive} when furthermore the last inequality 
is strict for any $\beta\neq 0$. 

For any $f_{_{G}}$-{\em strictly positive} complex vector bundle $E$, 
and any $\beta\in\ggot$ such that 
$M^{\beta}\cap f_{_G}^{-1}(\beta)\neq \emptyset$, we define 
$\eta_{_{E,\beta}}=\inf_{\xi}\langle \xi,\beta\rangle$, 
where $\xi$ runs over the set of weights for the 
$\tore_{\beta}$-action  on the fibers of each complex vector 
bundle $E\vert_{\Zcal}$, $\Zcal$ being a connected component
of $M^{\beta}$ that intersects $f_{_G}^{-1}(\beta)$.
\end{defi}

It is not difficult to see that a $f_{_G}$-moment bundle $L$ is
$f_{_{G}}$-strictly positive with  $\eta_{_{L,\beta}}=\parallel
\beta\parallel^2$, for any $\beta\in\ggot$ such that 
$M^{\beta}\cap f_{_G}^{-1}(\beta)\neq\emptyset$ (see Lemma \ref{L-a-positif}). 
The bundle $M\times\C\to M$ is the trivial 
example of $f_{_{G}}$-positive complex vector bundle over $M$.

\medskip

Let $\hgot_{+}$ be the choice of some positive Weyl chamber in $\hgot$. 
We prove in Lemma \ref{lem.C.f.G} that the set $\Bcal_{_{G}}:=
\{\beta\in\hgot_{+},\ M^{\beta}\cap f_{_{G}}^{-1}(\beta)\neq\emptyset\}$ 
is finite.

\medskip

{\bf Theorem B}\hspace{5mm}{\it
  Let $f_{_{G}}:M\to\ggot^{*}$ be an abstract moment map with $0$ as 
  regular value. Let $E$ be a $f_{_{G}}$-{\em strictly positive} $G$-complex 
  vector bundle  over $M$ (see Def. \ref{eq.mu.positif}). We have 
\begin{equation}\label{eq.theorem.B}
 \left[ RR^{^{G,J}}(M,E^{\stackrel{k}{\otimes}})\right]^G=
 \Qcal\Big(\Mcal_{red},E_{red}^{\stackrel{k}{\otimes}}\Big), 
 \quad  k\in\N-\{0\},
\end{equation}
if any of the following hold:

(i) $G=T$ is a torus; or

(ii) $k$ is large enough, so that 
$k.\eta_{_{E,\beta}}>\sum_{\alpha>0}\langle \alpha,\beta\rangle$, for 
any $\beta\in \Bcal_{_{G}}-\{0\}$; here the sum $\sum_{\alpha>0}$ 
is taken over the positive roots of $G$.

Moreover if (\ref{eq.J.red}) holds, (\ref{eq.theorem.B}) becomes 
$$
\left[ RR^{^{G,J}}(M,E^{\stackrel{k}{\otimes}})\right]^G=\pm 
 RR^{^{J_{red}}}\Big(\Mcal_{red},E_{red}^{\stackrel{k}{\otimes}}
 \otimes L_{\Dcal}\Big)\ .
$$
}

\bigskip

Let us explain why Theorem {\bf B} applied to a $G$-hermitian line 
bundle $L$ with the abstract moment map $f_{_{G}}=f_{_{L}}$ implies 
Theorem {\bf A}.  It is sufficient to prove that condition $ii)$ of 
Theorem {\bf A} implies condition $ii)$ of Theorem {\bf B}.  The 
curvature of $(L,\nabla^L)$ is $(\nabla^L)^2=-\imath\, \omega^L$, 
where $\omega^L$ is a differential $2$-form on $M$.  From the 
equivariant Bianchi formula (see Proposition 7.4 in \cite{B-G-V}) we 
get $\langle d f_{_{L}},X\rangle=-\omega^{L}(X_{M},-)$ for any 
$X\in\ggot$.  So, for any $\beta\in \Bcal_{_{G}}-\{0\}$, and $m\in 
M^{\beta}\cap f_{_{L}}^{-1}(\beta)$, the last equality gives $\langle 
d f_{_{L}}\vert_{m},\beta\rangle= 0$, hence $\beta$ is a critical 
value of $f_{_{L}}$.  Suppose now that $k\in\N$ is large enough so 
that the ball $\{ \xi\in\ggot^{*},\, \parallel\xi\parallel\leq
\frac{1}{k} \parallel\theta\parallel \}$ is included in the set of regular 
values of $f_{_{L}}$.  This gives first $\parallel\beta\parallel> 
\frac{1}{k} \parallel\theta\parallel$ and then
$\eta_{_{L,\beta}}=\parallel\beta\parallel^2> 
\frac{1}{k}\langle \theta,\beta\rangle$, 
for any $\beta\in \Bcal_{_{G}}-\{0\}$. $\Box$

\medskip

In the last section of this paper, we restrict ourselves to the 
Hamiltonian case.  In this situation, the
abstract moment map $f_{_G}$ and the almost complex structure  
$J$ are related by means of a $G$-invariant symplectic $2$-form $\omega$ : 

\begin{itemize}
  \item $f_{_{G}}$ is the moment map associated to a Hamiltonian action 
of $G$ over $(M,\omega)$ : $ d\langle f_{_G},X\rangle=-\omega(X_M,-)$ ,
 for $X\in\ggot$, and
 
  \item the data $(\omega, J)$ are {\em compatible} : 
$(v,w)\to \omega(v,Jw)$ is a Riemannian metric on $M$.
\end{itemize}

When $0$ is a regular value of $f_{_G}$, the compatible
data $(\omega, J)$ induce compatible data 
$(\omega_{red}, J_{red})$ on $\Mcal_{red}$. We have then a map  
$RR^{J_{red}}(\Mcal_{red},-)$. If $0$ is not a regular
value of $f_{_G}$, we consider elements $a$ in the principal
face $\tau$ of the Weyl chamber (see subsection \ref{non-regulier}). For
generic elements $a\in\tau\cap f_{_{G}}(M)$, the set
$\Mcal_{a}:=f_{_G}^{-1}(G\cdot a)/G$ is a symplectic orbifold and one
can consider the quantization map $RR^{J_{a}}(\Mcal_{a},-)$ 
relative to the choice of compatible almost complex structure $J_{a}$.

In this situation, we recover the results of 
\cite{Meinrenken98,Meinrenken-Sjamaar,Tian-Zhang98}.

\medskip

{\bf Theorem C}\ (Meinrenken, Meinrenken-Sjamaar, Tian-Zhang).\  
{\it Let $f_{_{G}}$ be the moment map 
associated to a Hamiltonian action  of $G$ over $(M,\omega)$, and let
$J$ be a $\omega$-compatible almost complex structure. Let 
$E\to M$ be a $G$-vector bundle. 

\noindent If $0\notin f_{_{G}}(M)$
and $E$ is $f_{_{G}}$-strictly positive, we have $[ RR^{^{G,J}}(M,E)]^G= 
0$. 

\noindent If $0\in f_{_{G}}(M)$  then : 
\begin{enumerate}
  \item[i)] If $0$ is a regular value, we have
$[ RR^{^{G,J}}(M,E)]^G=RR^{J_{red}}(\Mcal_{red},E_{red})$, if 
$E$ is $f_{_G}$-positive.
  \item[ii)]  If $0$ is not a regular value of $f_{_G}$ and $E=L$ is
a $f_{_G}$-moment bundle, we have 
$[ RR^{^{G,J}}(M,L)]^G=RR^{J_{a}}(\Mcal_{a},L_{a})$,
for every generic value $a$ of  $\tau\cap f_{_{G}}(M)$ sufficiently 
close to $0$. Here $L_a$ is the orbifold line bundle 
$L\vert_{f_{_G}^{-1}(G\cdot a)}/G$.
\end{enumerate}
}

\bigskip

We now turn to an introduction of our method. We associate to the 
abstract moment map $f_{_{G}}:M\to\ggot$ the vector field 
$$
\Hcal^{^G}_{m}=[f_{_{G}}(m)]_{M}.m,\quad m\in M\ ,
$$
and we denote by $C^{f_{_{G}}}$ the set where $\Hcal^{^G}$ vanishes.  
There are two important cases.  First, when the map $f_{_{G}}$ is 
constant, equal to an element $\gamma$ in the center of $\ggot$, the 
set $C^{f_{_{G}}}$ corresponds to the submanifold $M^{\gamma}$.  
Second, when $f_{_{G}}$ is the moment map associated with a 
Hamiltonian action of $G$ over $M$.  In this situation, Witten 
\cite{Witten} introduces the vector field $\Hcal^{^G}$ to propose, in 
the context of equivariant cohomology, a localization on the 
set of critical points of the function $\parallel f_{_{G}}\parallel^2$ 
: here $\Hcal^{^G}$ is the Hamiltonian vector field of 
$\frac{-1}{2}\parallel f_{_{G}}\parallel^2$, hence 
$\Hcal^{^G}_m=0\Longleftrightarrow d(\parallel 
f_{_{G}}\parallel^2)_m=0$.  This idea has been developed by the author 
in \cite{pep1,pep2}.

\medskip

Using a deformation argument in the context of transversally elliptic 
operator introduced by Atiyah \cite{Atiyah.74} and Vergne 
\cite{Vergne96}, we prove in 
section \ref{sec.general.procedure} that the 
map\footnote{We fix one for once a 
G-invariant almost complex structure $J$ and denote by $RR^{^{G}}$ the 
quantization map.} $RR^{^{G}}$ can be localized near
$C^{f_{_{G}}}$. More precisely, we have the finite 
decomposition $C^{f_{_{G}}}=\cup_{\beta\in\Bcal_{_{G}}}C_{\beta}^{^G}$ 
with $C_{\beta}^{^G}=G(M^{\beta}\cap f_{_{G}}^{-1}(\beta))$, and
\begin{equation}\label{RR.G.decompose}
RR^{^G}(M,E)= \sum_{\beta\in\Bcal_{_{G}}} RR^{^G}_{\beta}(M,E)\ .
\end{equation}
Each term $RR^{^G}_{\beta}(M,E)$ is a generalized character of $G$ 
that only depends on the behaviour of the data $M,E,J,f_{_{G}}$ near 
the subset $C_{\beta}^{^G}$.  In fact, $RR^{^G}_{\beta}(M,E)$ is the 
index of a transversally elliptic operator defined in an open 
neighbourhood of $C_{\beta}^{^G}$.

Our proof of Theorems {\bf B} and {\bf C} is in two steps. First we 
compute the term $RR^{^G}_{0}(M,E)$ which is the Riemann-Roch 
character localized near $f_{_{G}}^{-1}(0)$. After, we prove that
$[RR^{^G}_{\beta}(M,E)]^{G}=0$ for every $\beta\neq 0$.
For this purpose, the analysis of the localized Riemann-Roch 
characters $RR^{^G}_{\beta}(M,-):K_{G}(M)\to  
R^{-\infty}(G)$ is divided in three cases\footnote{$G_{\beta}$ is the 
stabilizer of $\beta$ in $G$.}  :
\medskip

\noindent{\bf Case 1 :} $\beta=0$

\noindent{\bf Case 2 :} $\beta\neq 0$ and $G_{\beta}=G$,

\noindent{\bf Case 3 :} $G_{\beta}\neq G$.

\medskip

We work out {\bf Case 1}  in subsection \ref{subsec.RR.G.O}.  We compute the 
generalized character $RR^{^G}_{0}(M,E)$ when $0$ is a regular value 
of $f_{_{G}}$.  We prove in particular that the multiplicity of the 
trivial representation in $RR^{^G}_{0}(M,E)$ is 
$\Qcal(\Mcal_{red},E_{red})$.  This last quantity is equal to $\pm 
RR^{J_{red}}(\Mcal_{red},E_{red}\otimes L_{\Dcal})$ when 
(\ref{eq.J.red}) holds.

\medskip

{\bf Case 2} is studied in section \ref{sec.loc.M.beta} for the 
particular situation where $f_{_{G}}$ is constant, equal to a 
G-invariant element $\beta\in\ggot$.  Then 
$C^{f_{_{G}}}=C_{\beta}^{^G}=M^{\beta}$, and (\ref{RR.G.decompose}) 
becomes $RR^{^G}(M,E)= RR^{^G}_{\beta}(M,E)$.  We prove then a 
localization formula (see (\ref{eq.localisation.1})) in the spirit of 
the Atiyah-Segal-Singer formula in equivariant K-theory 
\cite{Atiyah-Segal68,Segal68}.  Let us sketch out the result.

The normal bundle $\Ncal$ of $M^{\beta}$ in $M$ inherits a canonical 
complex structure $J_{\Ncal}$ on the fibers.  We denote by 
$\overline{\Ncal}\to M^{\beta}$ the complex vector bundle with the 
opposite complex structure.  The torus $\tore_{\beta}$ is included in 
the center of $G$, so the bundle $\overline{\Ncal}$ and the virtual 
bundle $\wedge_{\C}^{\bullet}\overline{\Ncal}:=\wedge_{\C}^{even} 
\overline{\Ncal}\stackrel{0}{\to}\wedge_{\C}^{odd}\overline{\Ncal}$ 
carry a $G\times \tore_{\beta}$-action: they can be considered as 
elements of $K_{G\times \tore_{\beta}}(M^{\beta})\, = 
\,K_{G}(M^{\beta})\otimes R(\tore_{\beta})$.  Let 
$K_{G}(M^{\beta})\,\widehat{\otimes} \,R(\tore_{\beta})$ be the vector 
space formed by the infinite formal sums $\sum_{a} E_{a}\, h^a$ taken 
over the set of weights of $\tore_{\beta}$, where $E_{a}\in 
K_{G}(M^{\beta})$ for every $a$.  The Riemann-Roch character $RR^{^G}$ 
can be extended to a map $RR^{^{G\times {\rm T}_{\beta}}}$ which 
satisfies the commutative diagram
\[
\xymatrix@C=2cm{
K_{G}(M^{\beta}) \ar[r]^{RR^{^G}} \ar[d] & R(G)\ar[d]^{k}\\
K_{G}(M^{\beta})\,\widehat{\otimes}\, R(\tore_{\beta})
\ar[r]^{RR^{^{G\times {\rm T}_{\beta}}}} & 
R(G)\,\widehat{\otimes}\, R(\tore_{\beta})\ .
   }
\] 
The arrow $k:R(G)\to R(G)\,\widehat{\otimes}\, R(\tore_{\beta})$ 
is the canonical map defined by $k(\phi)(g,h):=\phi(gh)$. We shall
notice that $[k(\phi)]^{G\times \tore_{\beta}}=[\phi]^{G}$.

In Section 5, we define an inverse, denoted by $\left[ 
\wedge_{\C}^{\bullet} \overline{\Ncal}\,\right]^{-1}_{\beta}$, of 
$\wedge_{\C}^{\bullet}\overline{\Ncal}$ in 
$K_{G}(M^{\beta})\,\widehat{\otimes}\, R(\tore_{\beta})$ which is 
polarized by $\beta$.  It means that $\left[ \wedge_{\C}^{\bullet} 
\overline{\Ncal}\,\right]^{-1}_{\beta}= \sum_{a} N_{a}\, h^a$ with 
$N_{a}\neq 0$ only if $\langle a, \beta\rangle\geq 0$.  We can state 
now our localization formula as the following equality in 
$R(G)\,\widehat{\otimes}\, R(\tore_{\beta})$ :
\begin{equation}\label{eq.localisation.1}
RR^{^G}(M,E)=RR^{^{G\times {\rm T}_{\beta}}}
\left(M^{\beta},E_{\vert M^{\beta}}\otimes 
\left[ \wedge_{\C}^{\bullet} \overline{\Ncal}\,\right]^{-1}_{\beta}\right)\ ,
\end{equation}
for every $E\in K_{G}(M)$.

In subsection \ref{subsec.RR.G.beta} we work out {\bf Case 2} for the 
general situation.  The map \break $RR^{^G}_{\beta}(M^{\beta},-)$ is 
the Riemann-Roch character on the $G$-manifold $M^{\beta}$, localized 
near $M^{\beta}\cap f_{_{G}}^{-1}(\beta)$, and we extend it to a map 
$RR^{^{G\times {\rm T}_{\beta}}}_{\beta}(M^{\beta},-):$ 
$K_{G}(M^{\beta})\,\widehat{\otimes}\, R(\tore_{\beta}) \to 
R^{-\infty}(G)\,\widehat{\otimes}\,R(\tore_{\beta})$. 
We prove then the following localization formula 
\begin{equation}\label{eq.localisation.2}
RR^{^G}_{\beta}(M,E)=RR^{^{G\times {\rm T}_{\beta}}}_{\beta}
\left(M^{\beta},E_{\vert M^{\beta}}\otimes 
\left[ \wedge_{\C}^{\bullet} \overline{\Ncal}\,\right]^{-1}_{\beta}\right)\ ,
\end{equation}
as an equality in 
$R^{-\infty}(G)\,\widehat{\otimes}\, R(\tore_{\beta})$. With 
(\ref{eq.localisation.2}) in hand, we see easily that \break 
$[RR^{^G}_{\beta}(M,E)]^{G}=0$ if the vector bundle $E$ is 
$f_{_{G}}$-strictly positive.

\medskip

Subsection \ref{subsec.induction.G.H} is devoted to {\bf Case 3}.  The 
abstract moment map $f_{_{G}}:M\to\ggot$ for the $G$-action on $M$ 
induces abstract moment maps $f_{_{G'}}:M\to\ggot'$ for every closed 
subgroup $G'$ of $G$.  For every $\beta\in\Bcal_{_{G}}$, we consider 
the Riemann-Roch characters $RR^{^G}_{\beta}(M,-)$, 
$RR^{^{G_{\beta}}}_{\beta}(M,-)$, and $RR^{^H}_{\beta}(M,-)$ localized 
respectively on $G(M^{\beta}\cap f_{_{G}}^{-1}(\beta))$, 
$M^{\beta}\cap f_{_{G}}^{-1}(\beta)$, and $M^{\beta}\cap 
f_{_{H}}^{-1}(\beta)$.  The major result of subsection 
\ref{subsec.induction.G.H} is the induction formulas proved in Theorem 
\ref{th.induction.G.H} and Corollary \ref{coro.induction.G.G.beta}, 
between these three characters.  I will explain briefly this result.

Let $W$ be the Weyl group associated to $(G,H)$. 
The choice of a Weyl chamber $\hgot^+$ in $\hgot$ determines 
a complex structure on the real vector space $\ggot/\hgot$.
Our induction formulas make a crucial use of the holomorphic induction
map $\HolH: R(H)\to R(G)$ (see (\ref{eq.holomorphe.G.H}) 
in Appendix B). Recall that $\HolH(h^{\lambda})$ is, for any 
weight $\lambda$, either equal to zero or to the 
character of an irreducible representation 
of $G$ (times $\pm 1$). In Theorem \ref{th.induction.G.H} we prove the 
following relation between $RR^{^{G}}_{\beta}(M,-)$ and 
$RR^{^{H}}_{\beta}(M,-)$
\begin{eqnarray}\label{eq.relation.induction.1}
  RR^{^{G}}_{\beta}(M,E)
  &=& \frac{1}{\vert W_{\beta}\vert}
  \HolH\left(\sum_{w\in W}w.RR^{^H}_{\beta}(M,E)\right)\\
  &=&\frac{1}{\vert W_{\beta}\vert}
  \HolH\left(RR^{^H}_{\beta}(M,E)\,\wedge_{\C}^{\bullet}
  \overline{\ggot/\hgot}\right)
  \nonumber \ ,
\end{eqnarray}  
where $W_{\beta}$ is the stabilizer of $\beta$ in $W$. 
In Corollary \ref{coro.induction.G.G.beta} we get the other relation:
\begin{equation}\label{eq.relation.induction.2}
 RR^{^{G}}_{\beta}(M,E)=\HolB \left(RR^{^{G_{\beta}}}_{\beta}(M,E)\,
 \wedge_{\C}^{\bullet}\overline{\ggot/\ggot_{\beta}}\right).
\end{equation}

Let us compare (\ref{eq.relation.induction.1}), with the 
Weyl integration formula\footnote{See Remark \ref{wedge.C-wedge.R}.}: 
for any $\phi\in R(G)$ we have
$\phi=\HolH\left(\phi_{\vert H}\right)=\HolH\left(\phi^{+}_{\vert H}
  \,  \wedge_{\C}^{\bullet}\overline{\ggot/\hgot}\right)$, 
where $\phi_{\vert H}$ is the restriction of $\phi$ to $H$, and 
$\phi^{+}_{\vert H}=\sum_{\lambda}m(\lambda)\, h^{\lambda}$ is the 
unique element in $R(H)\otimes\Q$ such that 
$\sum_{w\in W} w.\phi^{+}_{\vert H}= \phi_{\vert H}$
and  $m(\lambda)\neq 0$ only if $\lambda\in \hgot^+$. In 
(\ref{eq.relation.induction.1}), 
the $W$-invariant element 
$\frac{1}{\vert W_{\beta}\vert}\sum_{w\in W}w.RR^{^H}_{\beta}(M,E)$
plays the role of the restriction to $H$ of 
the character $\phi=RR^{^{G}}_{\beta}(M,E)$, and 
$\frac{1}{\vert W_{\beta}\vert}RR^{^H}_{\beta}(M,E)$ plays the role
of $\phi^{+}_{\vert H}$.

Since $\beta$ is a $G_{\beta}$-invariant element,  
(\ref{eq.relation.induction.2}) reduces 
the analysis of {\em Case 3} to the one of {\em Case 2}. 
From the result proved in {\em Case 2}, we have 
$[RR^{^{G_{\beta}}}_{\beta}(M^{\beta},E)]^{G_{\beta}}=0$ 
if the vector bundle $E$ is $f_{_{G}}$-strictly positive. 
But this does not implies in general that 
$[RR^{^G}_{\beta}(M,E)]^{G}=0$. We have
to take the tensor product of $E$ (so that 
$E^{\stackrel{k}{\otimes}}$ becomes more and more  
$f_{_{G_{\beta}}}$-strictly positive) to see that 
$[RR^{^G}_{\beta}(M,E^{\stackrel{k}{\otimes}})]^{G}=0$, when 
$\eta_{E^{\stackrel{k}{\otimes}},\beta}=
k.\eta_{_{E,\beta}}>\sum_{\alpha>0}\langle \alpha,\beta\rangle$.

\medskip

In the Hamiltonian setting considered in Section \ref{sec.Hamiltonien}, 
our strategy 
is the same, but at each step we obtain considerable refinements that 
are the principal ingredients of the proof of Theorem {\bf C}. 

{\bf Case 1 :}  When $0$ is a regular value of $f_{_{G}}$, we show 
that the Spin$^{\rm c}$ structure on $\Mcal_{red}$ is defined by 
$J_{red}$, hence $\Qcal(\Mcal_{red},-)=RR^{J_{red}}(\Mcal_{red},-)$. 
When $0$ is not a regular value of $f_{_{G}}$, 
we use the `shifting trick'  to compute the 
$G$-invariant part of  $RR^{^G}_{0}(M,E)$ 
(see subsection \ref{non-regulier}).

{\bf Case 2 :} For any $G$-invariant element $\beta\in\Bcal_{_G}$
with $\beta\neq 0$, we prove that the inverse 
$\left[ \wedge_{\C}^{\bullet} \overline{\Ncal}\,\right]^{-1}_{\beta}$ 
is of the form $\sum_{a} N_{a}\, h^a$ with $N_{a}\neq 0$ only if 
$\langle a, \beta\rangle > 0$ (in general we have only 
$\langle a, \beta\rangle \geq 0$). 

{\bf Case 3 :} For $\beta\in\Bcal_{_{G}}$ with $G_{\beta}\neq G$, we consider 
the open face $\sigma$ of the Weyl chamber which contains $\beta$, 
and the corresponding symplectic slice $\Ycal_{\sigma}$
which is a $G_{\beta}$-symplectic submanifold 
of $M$. The localized Riemann-Roch characters 
$RR^{^{G}}_{\beta}(M,E)$ and $RR^{^{G_{\beta}}}_{\beta}(\Ycal_{\sigma},-)$ 
are related by the following induction formula
$$
RR^{^{G}}_{\beta}(M,E)=
\HolB\left(RR^{^{G_{\beta}}}_{\beta}(\Ycal_{\sigma},
E\vert_{\Ycal_{\sigma}})\right)\ .
$$


\bigskip

{\bf Acknowledgments.} I am grateful to Mich\`ele Vergne for her 
interest in this work, especially for the useful discussions and 
insightful suggestions on a preliminary version of this paper.

\bigskip

\begin{center}
    {\bf Notation}
\end{center}

Throughout the paper $G$ will denote a compact, connected Lie group, 
and $\ggot$ its Lie algebra.  We let $H$ be a maximal torus in $G$, 
and $\hgot$ be its Lie algebra.  The integral lattice 
$\Lambda\subset\hgot$ is defined as the kernel of $\exp:\hgot\to H$, 
and the real weight lattice $\Lambda^* \subset\hgot^*$ is defined by : 
$\Lambda^*:=\hom(\Lambda,2\pi\Z)$.  Every $\lambda\in\Lambda^*$ 
defines a 1-dimensionnal $H$-representation, denoted $\C_{\lambda}$, 
where $h=\exp X$ acts by 
$h^{\lambda}:=e^{\imath\langle\lambda,X\rangle}$.  We let $W$ be the 
Weyl group of $(G,H)$, and we fix the positive Weyl chambers 
$\hgot_+\subset\hgot$ and $\hgot^*_+\subset\hgot^*$.  For any dominant 
weight $\lambda\in\Lambda^*_+:=\Lambda^*\cap\hgot^*_+$, we denote by  
$V_{\lambda}$ the $G$-irreducible representation with highest weight 
$\lambda$, and $\chi_{_{\lambda}}^{_G}$ its character.  We denote by 
$R(G)$ (resp.  $R(H)$) the ring of characters of finite-dimensional 
$G$-representations (resp.  $H$-representations).  We denote by 
$R^{-\infty}(G)$ (resp.  $R^{-\infty}(H)$) the set of generalized 
characters of $G$ (resp.  $H$).  An element $\chi\in R^{-\infty}(G)$ 
is of the form $\chi=\sum_{\lambda\in\Lambda^{*}_{+}}m_{\lambda}\, 
\chi_{_{\lambda}}^{_G}\,$, where $\lambda\mapsto m_{\lambda}, 
\Lambda^{*}_{+}\to\Z$ has at most polynomial growth.  In the same way, 
an element $\chi\in R^{-\infty}(H)$ is of the form 
$\chi=\sum_{\lambda\in\Lambda^{*}}m_{\lambda}\, h^{\lambda}$, where 
$\lambda\mapsto m_{\lambda}, \Lambda^{*}\to\Z$ has at most 
polynomial growth.

Some additional notation will be introduced later :     

\begin{itemize}
    
   \item[]$G_{\gamma}$ : stabilizer subgroup of $\gamma\in \ggot$
   
   \item[]$\tore_{\beta}$ : torus generated by $\beta\in \ggot$
    
   \item[]$M^{\gamma}$ : submanifold of points fixed by $\gamma\in\ggot$
   
   \item[]$\T M$ : tangent bundle of $M$
   
   \item[]$\T_{G} M$ : set of tangent vectors orthogonal to the 
   $G$-orbits in $M$
   
   \item[]$\fgene(G)^{G}$ : set of generalized functions on $G$, invariant 
   by conjugation 
   
   \item[]${{\rm Ind}^{^G}_{_{G_{\gamma}}}}:\fgene(G_{\gamma})^{G_{\gamma}}
   \to \fgene(G)^{G}$ : induction map
   
   \item[]${{\rm Hol}^{^G}_{_{G_{\gamma}}}}: R(G_{\gamma})
   \to R(G)$ : holomorphic induction map
   
   \item[]$RR^{^G}_{\beta}(M,-)$ : Riemann-Roch character localized on 
   $G.(M^{\beta}\cap f_{_{G}}^{-1}(\beta))$
   
   \item[]$\Char(\sigma)$ : characteristic set of the symbol $\sigma$
   
   \item[]$\Thom_{G}(M,J)$ : Thom symbol 
   
   \item[]$\Thom^{\gamma}_{G}(M)$ : Thom symbol 
   localized near $M^{\gamma}$ 
   
   \item[]$\Thom_{G,\beta}^f(M)$ : Thom symbol 
   localized near $G.(M^{\beta}\cap f_{_{G}}^{-1}(\beta))$. 
\end{itemize}


\bigskip

\section{Quantization of compact manifolds}\label{sec.quantization}

\medskip

Let $M$ be a compact manifold provided with an action of a compact 
connected Lie group $G$.  A $G$-invariant almost complex structure $J$ 
on $M$ defines a map $RR^{^{G,J}}(M,-) : K_{G}(M)\to R(G)$ from the 
equivariant $K$-theory of complex vector bundles over $M$ to the 
character ring of $G$.

Let us recall the definition of this map. The almost complex 
structure on $M$ gives the decomposition
$\wedge \T^{*} M \otimes \C =\oplus_{i,j}\wedge^{i,j}\T^* M$
of the bundle of differential forms. Using Hermitian structure in the tangent 
bundle $\T M$ of $M$, and in the fibers of $E$, we define a twisted 
Dirac operator
$$
\Dcal^{+}_{E}:\Acal^{0,even}(M,E)\to\Acal^{0,odd}(M,E)
$$
where $\Acal^{i,j}(M,E):=\Gamma(M,\wedge^{i,j}\T^{*}M\otimes_{\C}E)$ 
is the space of $E$-valued forms of type $(i,j)$. The Riemann-Roch
character $RR^{^{G,J}}(M,E)$ is defined as the index of the elliptic
operator $\Dcal^{+}_{E}$:
$$
RR^{^{G,J}}(M,E)= [Ker\Dcal^{+}_{E}] - [Coker\Dcal^{+}_{E}].
$$
In fact, the virtual character $RR^{^{G,J}}(M,E)$ is independent of 
the choice of the Hermitian metrics on the vector bundles $\T M$ 
and $E$.

If $M$ is a compact complex analytic manifold, and $E$ is an 
holomorphic complex vector bundle, we have 
$RR^{^{G,J}}(M,E)=\sum_{q=0}^{q=dimM}(-1)^q [\Hcal^q(M,\Ocal(E))]$,
where $\Hcal^q(M,\Ocal(E))$ is the $q$-th cohomology group of the 
sheaf $\Ocal(E)$ of the holomorphic sections of $E$ over $M$.

In this paper, we shall use an equivalent definition of the map 
$RR^{^{G,J}}$.
We associate to an invariant almost complex structure $J$ 
the symbol $\Thom_{G}(M,J)\in K_{G}(\T M)$ defined as follows. Consider
a Riemannian structure $q$ on $M$ such that the endomorphism $J$ is 
orthogonal relatively to $q$, and let $h$ be the following Hermitian
structure on $\T M$ : $ h(v,w)=q(v,w) -\imath q(Jv,w)$ for $v,w\in 
\T M$.
Let $p:\T M\to M$ be the canonical projection. The symbol 
$\Thom_{G}(M,J):p^{*}(\wedge_{\C}^{even} \T M)\to p^{*}
(\wedge_{\C}^{odd} \T M)$  is equal, at 
$(x,v)\in \T M$, to the Clifford map
\begin{equation}\label{eq.thom.complex}
 Cl_{x}(v)\ :\ p^{*}(\wedge_{\C}^{even} \T M)\vert_{(x,v)}
\longrightarrow p^{*}(\wedge_{\C}^{odd} \T M)\vert_{(x,v)},
\end{equation}
where $Cl_{x}(v).w= v\wedge w - c_{h}(v).w$ for $w\in 
\wedge_{\C}^{\bullet} \T_{x}M$. Here $c_{h}(v):\wedge_{\C}^{\bullet} 
\T_{x}M\to\wedge^{\bullet -1} \T_{x}M$ denotes the 
contraction map relatively to $h$ : for $w\in \T_{x}M$ we have
$c_{h}(v).w=h(w,v)$. Here $(\T M, J)$ is considered as a 
complex vector bundle over $M$.

The symbol $\Thom_{G}(M,J)$  determines the Bott-Thom 
isomorphism $\Thom_{J}:$\break  
$K_{G}(M)\longrightarrow K_{G}(\T M)$
by $\Thom_{J}(E):=\Thom_{G}(M,J)\otimes p^*(E), \ E\in K_{G}(M)$. To 
make the notation clearer, $\Thom_{J}(E)$ is the symbol 
$\sigma^{_{E}} : p^{*}(\wedge_{\C}^{even} \T M\otimes E)\to 
p^{*}(\wedge_{\C}^{odd} \T M\otimes E)$ with
\begin{equation}
\sigma^{_{E}}(x,v):= Cl_{x}(v)\otimes Id_{E_{x}}\ ,\quad  (x,v)\in \T M,
    \label{eq:thom-iso-b}
\end{equation}
where $E_{x}$ is the fiber of $E$ at $x\in M$. 

Consider the index map $\indice_{M}^{G} : K_{G}(\T^{*}M)\to R(G)$
where $\T^{*}M$ is the cotangent bundle of $M$.  Using a $G$-invariant
auxiliary metric on $\T M$, we can identify the vector bundle
$\T^{*}M$ and $\T M$, and produce an `index' map $\indice_{M}^{G} :
K_{G}(\T M)\to R(G)$.  We verify easily that this map is independent
of the choice of the metric on $\T M$.

\begin{lem}\label{lem.quantization}
    We have the following commutative diagram
\begin{equation}
\xymatrix@C=2cm{
 K_{G}(M)\ar[r]^{\Thom_{J}} \ar[dr]_{RR^{^{G,J}}} & 
 K_{G}(\T M) \ar[d]^{\indice_{M}^{G}}\\
     & R(G)\ .}
\end{equation}     
        
\end{lem}

{\em Proof} : If we use the natural identification 
$(\wedge^{0,1}\T^{*}M,\imath)\cong
(\T M,J)$ of complex vector bundles over $M$, we see that the
principal symbol of the operator $\Dcal^{+}_{E}$ is equal to
$\sigma^{_E}$ (see \cite{Duistermaat96}).

\medskip

We will conclude with the following Lemma. Let $J^{0}, J^{1}$ be two 
$G$-invariant almost complex structures on $M$, and let $RR^{^{G,J^{0}}},
RR^{^{G,J^{1}}}$ be the respective quantization maps.

\medskip

\begin{lem} \label{lem.inv.homotopy}
The maps $RR^{^{G,J^{0}}}$ and $RR^{^{G,J^{1}}}$ are identical in 
the following cases:

\noindent i) There exists a $G$-invariant section $A\in 
\Gamma(M,\End(\T M))$, {\em homotopic to the identity} in  
$\Gamma(M,\End(\T M))^{G}$ such that $A_{x}$ is invertible, 
and $A_{x}.J^{0}_{x}=J^{1}_{x}.A_{x}$ for every $x\in M$.

\noindent ii)  There exists an homotopy $J^t,\ t\in[0,1]$ of 
$G$-invariant almost complex structures between $J^{0}$ and $J^{1}$.
\end{lem}

\medskip

{\em Proof of i)} :  Take a Riemannian structure $q^{1}$ on $M$ such 
that $J^{1}\in O(q^{1})$ and define another Riemannian structure 
$q^{0}$ by $q^{0}(v,w)=q^{1}(Av,Aw)$ so that $J^{0}\in O(q^{0})$. 
The section $A$ defines a bundle unitary map $\underline{A}:
(\T M,J^{0},h^{0})\to(\T M,J^{1},h^{1}),\ (x,v)\to (x,A_{x}.v)$, 
where $h^l(.,.):=q^l(.,.)-\imath q^l(J^l.,.),\ 
l=0,1$. This gives an isomorphism  $A^{\wedge}_{x}:
\wedge_{J^{0}}\T_{x}M\to \wedge_{J^{1}}\T_{x}M$ such that the 
following diagram is commutative

\[
\xymatrix@C=2cm{
\wedge_{J^{0}}\T_{x}M\ar[r]^{Cl_{x}(v)} \ar[d]_{A^{\wedge}_{x}} & 
\wedge_{J^{0}}\T_{x}M\ar[d]^{A^{\wedge}_{x}}\\
\wedge_{J^{1}}\T_{x}M \ar[r]^{Cl_{x}(A_{x}.v)} & 
\wedge_{J^{1}}\T_{x}M\ .
   }
\] 

Then $A^{\wedge}$ induces an isomorphism between the symbols
$\Thom_{G}(M,J^{0})$ and \break $\underline{A}^{*}(\Thom_{G}(M,J^{1}))
\ :\ (x,v)\to \Thom_{G}(M,J^{1})(x,A_{x}.v)$.  Here
$\underline{A}^{*}:K_{G}(\T M)\to K_{G}(\T M)$ is the map induced by
the isomorphism $\underline{A}$.  Thus the complexes \break
$\Thom_{G}(M,J^{0})$ and $\underline{A}^{*}(\Thom_{G}(M,J^{1}))$
define the same class in $K_{G}(\T M)$.  Since $A$ is
homotopic to the identity, we have $\underline{A}^{*}={\rm Identity}$. 
We have proved that $\Thom_{G}(M,J^{0})= \Thom_{G}(M,J^{1})$ in
$K_{G}(\T M)$, and by Lemma \ref{lem.quantization} this shows that
$RR^{^{G,J^{0}}}= RR^{^{G,J^{1}}}$.

{\em Proof of ii)} : We construct $A$ as in $i)$.  Take first 
$A^{1,0}:= Id -J^{1}J^{0}$ and remark that $A^{1,0}.J^{0}=J^{1}.A^{1,0}$.
Here we consider the homotopy $A^{1,0}_{u}:= Id -uJ^{1}J^{0},\ u\in [0,1]$.
If $-J^{1}J^{0}$ is close to $Id$, for example 
$\vert Id + J^{1}J^{0}\vert\leq 1/2$, the bundle map $A^{1,0}_{u}$ will be 
invertible for every $u\in [0,1]$. Then we can conclude with Point
$i)$.  In general we use the homotopy $J^t,\ t\in[0,1]$.  First, 
we decompose the interval $[0,1]$ in 
$0=t_{0}<t_1<\cdots<t_{k-1}<t_{k}=1$ and we consider the maps 
$A^{t_{l+1},t_{l}}:= Id -J^{t_{l+1}}J^{t_{l}}$, with the corresponding
homotopy $A^{t_{l+1},t_{l}}_{u},\ u\in [0,1]$, for $l=0,\ldots,k-1$.  
Because $-J^{t_{l+1}}J^{t_{l}}\to Id$ when $t\to t'$, the bundle maps 
$A^{t_{l+1},t_{l}}_{u}$ are invertible for all $u\in [0,1]$ if 
$t_{l+1}-t_{l}$ is small enough.  Then we take the $G$-equivariant 
bundle map $A:= \Pi_{l=0}^{k-1}A^{t_{l+1},t_{l}}$ with the homotopy
$A_{u}:= \Pi_{l=0}^{k-1}A^{t_{l+1},t_{l}}_{u},\ u\in [0,1]$. We have 
$A.J^{0}=J^{1}.A$ and $A_{u}$ is invertible for every $u\in [0,1]$, 
hence we conclude with the point $i)$.
$\Box$

\medskip


\section{Transversally elliptic symbols}\label{sec.T.G.M}

\medskip

We give here a brief review of the material we need in the next 
sections.  The references are 
\cite{Atiyah.74,B-V.inventiones.96.1,B-V.inventiones.96.2,Vergne96}.

Let $M$ be a {\em compact} manifold provided with a $G$-action.  Like 
in the previous section, we identify the tangent bundle $\T M$ and the 
cotangent bundle $\T^{*}M$ via a $G$-invariant metric $(.,.)_{_{M}}$ 
on $\T M$.  For any $X\in \ggot$, we denote by $X_{M}$ the following 
vector field : for $m\in M$, $X_{M}(m):= \frac{d}{dt}\exp(-tX).m 
|_{t=0}$.

If $E^{0},E^{1}$ are $G$-equivariant vector bundles over $M$, a 
morphism \break
$\sigma \in \Gamma(\T M,\hom(p^{*}E^{0},p^{*}E^{1}))$ of $G$-equivariant
complex vector bundles will be called a symbol. 
The subset of all $(x,v)\in \T M$ where $\sigma(x,v): E^{0}_{x}\to 
E^{1}_{x}$ is not invertible will be called the characteristic set 
of $\sigma$, and denoted $\Char(\sigma)$. 

We denote by $\T_{G}M$ the following subset of $\T M$ :
$$
   \T_{G}M\ = \left\{(x,v)\in \T M,\ (v,X_{M}(m))_{_{M}}=0 \quad {\rm for\ all}\ 
   X\in\ggot \right\} .
$$

A symbol $\sigma$ will be called {\em elliptic} if $\sigma$ is 
invertible outside a compact subset of $\T M$ ($\Char(\sigma)$ is 
compact), and it will be called {\em transversally elliptic} if the 
restriction of $\sigma$ to $\T_{G}M$ is invertible outside a compact 
subset of $\T_{G}M$ ($\Char(\sigma)\cap \T_{G}M$ is compact).  An 
elliptic symbol $\sigma$ defines an element of $K_{G}(\T M)$, and the 
index of $\sigma$ is a virtual finite dimensional representation of 
$G$ 
\cite{Atiyah-Segal68,Atiyah-Singer-1,Atiyah-Singer-2,Atiyah-Singer-3}.  
A transversally elliptic symbol $\sigma$ defines an element of 
$K_{G}(\T_{G}M)$, and the index of $\sigma$ is defined (see 
\cite{Atiyah.74} for the analytic index and 
\cite{B-V.inventiones.96.1,B-V.inventiones.96.2} for the cohomological 
one) and is a trace class virtual representation of $G$.  Remark that 
any elliptic symbol of $\T M$ is transversally elliptic, hence we have 
a restriction map $K_{G}(\T M)\to K_{G}(\T_{G}M)$ which makes the 
following diagram 
\begin{equation}\label{indice.generalise}
\xymatrix{
K_{G}(\T M)\ar[r]\ar[d]_{\indice_{M}^G} & 
K_{G}(\T_{G}M)\ar[d]^{\indice_{M}^G}\\
R(G)\ar[r] &  R^{-\infty}(G)\ .
   }
\end{equation} 
commutative. 

\subsection{Index map on non-compact manifolds}
\label{subsec.indice.ouvert} 

Let $U$ be a non-compact $G$-manifold. Lemma 3.6 and Theorem 3.7
of \cite{Atiyah.74} tell us that for any open G-embedding
$j:U\croc M$ into a compact manifold we have a pushforward map
$j_{*}:K_{G}(\T_{G}U)\to K_{G}(\T_{G}M)$ such that the 
composition 
$$
K_{G}(\T_{G}U)\stackrel{j_{*}}{\longrightarrow} K_{G}(\T_{G}M)
\stackrel{\indice_{M}^G}{\longrightarrow} R^{-\infty}(G)
$$
is independent of the choice of $j:U\croc M$.

\begin{lem}
    Let $U$ be a G-invariant open subset of a
    $G$-manifold $\Xcal$. If $U$ is relatively compact, 
    there exists an open G-embedding
   $j:U\croc M$ into a compact $G$-manifold.
\end{lem}   

{\em Proof :} Here we follow the proof given by Boutet de Monvel in 
\cite{Boutet.70}. Let $\chi\in\f(\Xcal)^{G}$ be a function with 
compact support, 
such that $0\leq \chi\leq 1$ and $\chi=1$ on $U$. 
Let $q:\Xcal\times \R\to \R$ 
be the function defined by $q(m,t)=\chi(m)-t^2$. The 
interval $(-\infty, 1]$ is the image of $q$, and the fibers
$q^{-1}(\esp)$ are compact for every $\esp > 0$. According to Sard's 
Theorem there exists a regular value $0<\esp_{0}<1$ of $q$.
The set $q^{-1}(\esp_{0})$ is then a compact $G$-invariant submanifold 
of $\Xcal\times \R$, and $j:U\to q^{-1}(\esp_{0})$, 
$m\mapsto (m,\sqrt{1-\esp_{0}})$ is an open embedding. $\Box$

\begin{coro}\label{hyp.indice}
    The index map $\indice_{U}^{G}:K_{G}(\T_{G}U)\to R^{-\infty}(G)$
    is defined when $U$ is a $G$-invariant relatively compact
    open subset of a $G$-manifold.
\end{coro}

\subsection{Excision lemma}\label{subsec.excision}

Let $j:U\hookrightarrow M$ be the inclusion map of a $G$-invariant 
open subset on a compact manifold, and let 
$j_{*}:K_{G}(\T_{G}U)\to K_{G}(\T_{G}M)$ be the pushforward map. We have 
two index maps $\indice_{M}^{G}$, and $\indice^{G}_{U}$ such that 
$\indice_{M}^{G}\circ j_{*}= \indice^{G}_{U}$. Suppose that 
$\sigma$ is a transversally elliptic symbol on $\T M$ with 
characteristic set contained in $\T M |_{U}$. Then, the restriction 
$\sigma |_{U}$ of $\sigma$ to $\T U$ is a transversally elliptic 
symbol on $\T U$, and
\begin{equation}
    j_{*}(\sigma|_{U})=\sigma \quad {\rm in}\quad K_{G}(\T_{G}M). 
     \label{eq:excision}
\end{equation}
In particular, it gives 
$\indice_{M}^{G}(\sigma)=\indice^{G}_{U}(\sigma|_{U})$.

\medskip

\subsection{Locally free action}\label{subsec.free.action}

Let $G$ and $H$ be compact Lie groups and let $M$ be a {\em compact} 
$G\times H$ manifold

In a first place, we suppose that $G$ acts freely on $M$,  and we 
denote by $\pi : M\to M/G$ the principal fibration. Note that the map $\pi$ is 
$H$-equivariant. In this situation we have $\T_{G\times H}M\widetilde{=}
\pi^{*}(\T_{H}(M/G))$, and thus an isomorphism
\begin{equation}
  \pi^{*}\ :\ K_{H}(\T_{H}(M/G))\longrightarrow 
  K_{G\times H}(\T_{G\times H}M)\ .
    \label{eq:free.action}
\end{equation}

We rephrase now  Theorem 3.1 of Atiyah in \cite{Atiyah.74}. 
Let $\{W_{a}, a\in \hat{G} \}$ be a completed set of 
inequivalent irreducible representations of $G$. 

For each irreducible $G$-representation $V_{\mu}$, 
we associate the complex vector bundle 
$\underline{V}_{\mu}:=M\times_{H}V_{\mu}$ on $M/G$ and denote by  
$\underline{V}_{\mu}^{*}$
its dual. The group $H$ acts trivially on $V_{\mu}$, this makes 
$\underline{V}_{\mu}^{*}$ a $H$-vector bundle.

\begin{theo}[Atiyah]\label{thm.atiyah.1}
    If $\sigma\in K_{H}(\T_{H}(M/G))$, then we have the following 
    equality in $R^{-\infty}(G\times H)$
\begin{equation}\label{eq:atiyah.1}
\indice^{G\times H}_{M}(\pi^{*}\sigma)\ =\ \sum_{\mu\in \Lambda^{*}_{+}}
\indice^{H}_{M/G}(\sigma\otimes \underline{V}_{\mu}^{*}) . V_{\mu}\quad.
\end{equation}
In particular the $G$-invariant part of $\indice^{G\times 
H}_{M}(\pi^{*}\sigma)$ is $\indice^{H}_{M/G}(\sigma)$.
\end{theo}

A classical example is when $M=G$, $G=G_{r}$ acts by right 
multiplications on $G$, and $G=G_{l}$ acts by left multiplications 
on $G$. Then the zero map $\sigma_{0}:G\times \C\to G\times\{0\}$ 
defines a $G_{r}\times G_{l}$-transversally elliptic symbol 
associated to the zero differential operator $\f(G)\to 0$. 
This symbol is equal to the pullback of 
$\C\in K_{G_{r}}(\T_{G_{r}}\{{\rm point}\})\widetilde{=}R(G_{r})$. 
In this case $\indice^{G_{r}\times G_{l}}_{G}(\sigma_{0})$ is equal 
to $L^{2}(G)$, the $L^{2}$-index 
of the zero operator on $\f(G)$. The $G_{r}$-vector bundle 
$\underline{V}^{*}_{\mu}\to\{{\rm point}\}$ is just the vector space 
$V_{\mu}^{*}$ with the canonical action of $G_{r}$. Finally, 
(\ref{eq:atiyah.1}) is the Peter-Weyl decomposition of $L^{2}(G)$ in 
$R^{-\infty}(G_{r}\times G_{l})$:
$\ L^2(G)=\sum_{\mu\in \Lambda^{*}_{+}}V_{\mu}^{*}\otimes V_{\mu}$.

\medskip

We suppose now that  $G$ acts locally freely on  $M$. The quotient 
$\Xcal:=M/G$ is an orbifold, a space with finite-quotient 
singularities.  One considers on $\Xcal$ the $H$-equivariant {\em proper} 
orbifold vector bundles and the corresponding $R(H)$-module $K_{orb,H}(\Xcal)$
\cite{Kawasaki81}. In the same way we consider the $H$-equivariant 
proper elliptic 
symbols on the orbifold $\T \Xcal$ and the corresponding $R(H)$-module 
$K_{orb,H}(\T\Xcal)$. The principal fibration  $\pi: M\to \Xcal$ induces 
isomorphisms $K_{orb,H}(\Xcal)\simeq K_{G\times H}(M)$ and 
$K_{orb,H}(\T\Xcal)\simeq K_{G\times H}(\T_{H}M)$ that we 
both denote by $\pi^{*}$. The index map
\begin{equation}\label{index-orbifold}
    \indice^{H}_{\Xcal}: K_{orb,H}(\T\Xcal)\to R(H)
\end{equation}
is defined by the following equation: for any $\sigma\in K_{orb,H}(\T\Xcal)$,    
$\indice^{H}_{\Xcal}(\sigma):=[\indice^{G\times 
H}_{M}(\pi^{*}\sigma)]^{G}$. 

We are particularly interested in the case where the bundle 
$\T_{G}M\to M$ carries a $G\times H$-equivariant almost complex structure $J$. 
Taking the quotient by $G$, it defines  a $H$-equivariant almost complex structure 
$J_{\Xcal}$ on the orbifold tangent bundle $\T\Xcal\to\Xcal$. Like in the
smooth case, we have the Thom symbol 
$\Thom_{H}(\Xcal,J_{\Xcal})\in K_{orb,H}(\T\Xcal)$ and a Riemmann-Roch 
character $RR^{^H}:K_{orb,H}(\Xcal)\to R(H)$  related as in Lemma 
\ref{lem.quantization}.

\medskip
\subsection{Induction}\label{subsec.induction.def}

Let $i:H\croc G$ be a 
closed subgroup with Lie algebra $\hgot$, and $\Ycal$ be a  
$H$-manifold (as in Corollary 
\ref{hyp.indice}). We have two principal bundles $\pi_{1}: 
G\times\Ycal\to \Ycal$ for the $G$-action, and $\pi_{2}:
G\times\Ycal\to \Xcal:=G\times_{H}\Ycal$ for the diagonal $H$-action. 
The map $i_{*}: K_{H}(\T_{H}\Ycal)\to K_{G}(\T_{G}\Xcal)$ is well defined 
by the following commutative diagram
\begin{equation}\label{eq:G.H.induction}
\xymatrix@C=25mm{
K_{H}(\T_{H}\Ycal)\ar[r]^{\pi_{1}^{*}} \ar[dr]_{i_{*}} & 
 K_{G\times H}(\T_{G\times H}(G\times\Ycal)) \ar[d]^{(\pi_{2}^{*})^{-1}}\\
     & K_{G}(\T_{G}\Xcal)\ ,}
\end{equation}    
 since $\pi_{1}^{*}$ and $\pi_{2}^{*}$ are isomorphisms.
 
Let us show how to compute $i_{*}(\sigma)$,  for an $H$-transversally elliptic 
symbol $\sigma \in \Gamma( \T Y, \hom(E^{0},E^{1}))$, 
where $E^{0}, E^{1}$ are $H$-equivariant vector bundles 
over $\T\Ycal$. First we notice\footnote{\label{eq.espace.tangent}
These identities come from the following $G\times H$-equivariant 
isomorphism of vector bundles over $G\times\Ycal$:
$\T_{H}(G\times\Ycal)\to G\times(\ggot/\hgot \oplus \T \Ycal),
(g,m;\frac{d}{dt}_{\vert t=0}(g.e^{tX})+ v_{m})\mapsto
(g,m; pr_{\ggot/\hgot}(X)+v_{m})$. Here $pr_{\ggot/\hgot}:
\ggot\to \ggot/\hgot$ is the orthogonal projection.}
that $\T(G\times_{H}\Ycal)\widetilde{=}
G\times_{H}(\ggot/\hgot \oplus \T \Ycal)$, and  
$\T_{G}(G\times_{H}\Ycal)\widetilde{=}G\times_{H}( \T_{H} \Ycal)$.
So we extend trivially 
$\sigma$ to $\ggot/\hgot\oplus\T\Ycal$, and we define 
$i_{*}(\sigma)\in \Gamma(G\times_{H}(\ggot/\hgot \oplus \T\Ycal), 
\hom(G\times_{H}E^{0},G\times_{H}E^{1}))$ by 
$i_{*}(\sigma)([g;\xi,x,v]):=\sigma(x,v)$ for 
$g\in G$, $\xi\in\ggot/\hgot$ and $(x,v)\in \T\Ycal$.

To express the $G$-index of $i_{*}(\sigma)$  in terms of the $H$-index 
of $\sigma$, we need the induction map
\begin{equation}
    \indH : \fgene(H)^{H}\longrightarrow \fgene(G)^{G}\ ,
    \label{eq:fonction.induction}
\end{equation}
where $\fgene(H)$ is the set of generalized functions on $H$, and 
the $H$ and $G$ invariants are taken with the conjugation action.
The map $\indH$ is defined as follows : for $\phi\in\fgene(H)^{H}$, 
we have $\int_{G}\indH(\phi)(g)f(g)dg= {\rm cst}
\int_{H}\phi(h)f|_{H}(h)dh$, 
for every $f\in\f(G)^{G}$, where 
${\rm cst}=\Vol(G,dg)/\Vol(H,dh)$. 

We can now recall Theorem 4.1 of Atiyah in \cite{Atiyah.74}.

\begin{theo}\label{thm.atiyah.2}
Let $i:H\to G$ be the inclusion of a closed subgroup, let $\Ycal$ be a 
$H$-manifold satisfying the hypothesis of Corollary \ref{hyp.indice}, and set 
$\Xcal=G\times_{H}\Ycal$. Then we have the commutative diagram
\[
\xymatrix{
K_{H}(\T_{H}\Ycal)\ar[r]^{i_{*}}\ar[d]_{\indice_{\Ycal}^H} & 
K_{G}(\T_{G}\Xcal)\ar[d]^{\indice_{\Xcal}^G}\\
\fgene(H)^{H}\ar[r]_{\indH} &  \fgene(G)^{G}\ .
   }
\] 
\end{theo}

\subsection{Reduction}\label{subsec.reduction}

Let us recall a multiplicative property of the index for the product 
of manifold. Let  a compact Lie group $G$ acts on two 
manifolds $\Xcal$ and $\Ycal$, and assume that another compact 
Lie group $H$ 
acts on $\Ycal$ commuting with the action of $G$.
The external product of complexes on $\T\Xcal$ and $\T \Ycal$ induces
a multiplication (see \cite{Atiyah.74} and \cite{Vergne96}, section 2):

\begin{eqnarray}
K_{G}(\T \Xcal)\times K_{G\times H}(\T \Ycal)&\longrightarrow &
K_{G\times H}(\T (\Xcal\times \Ycal)) \\
(\sigma_{1},\sigma_{2})&\longmapsto & \sigma_{1}\odot  
\sigma_{2} \ .\nonumber
    \label{eq:produit.usuel}
\end{eqnarray}

Let us recall the definition of this external product. Let 
$E^{\pm},F^{\pm}$ be $G\times H$-equivariant Hermitian vector bundles over 
$\Xcal$ and $\Ycal$ respectively, and let $\sigma_{1}:E^{+}\to E^{-}$,
$\sigma_{2}:F^{+}\to F^{-}$ be $G\times H$-equivariant symbols. We consider 
the $G\times H$-equivariant symbol 
$$
\sigma_{1}\odot \sigma_{2}:E^{+}\otimes F^{+}\oplus E^{-}\otimes F^{-}
\longrightarrow E^{-}\otimes F^{+} \oplus E^{+}\otimes F^{-}
$$ 
defined by

\begin{equation}\label{eq:produit.externe}
\sigma_{1}\odot \sigma_{2}=
\left(
\begin{array}{cc}
  \sigma_{1}\otimes I & -I\otimes \sigma_{2}^{*}\\
I \otimes \sigma_{2} & \sigma_{1}^{*}\otimes I
\end{array}
\right)\ .
\end{equation}

We see that the set $\Char(\sigma_{1}\odot \sigma_{2})\subset
\T\Xcal\times\T\Ycal$ is equal to 
$\Char(\sigma_{1})\times \Char(\sigma_{2})$. 
This exterior product defines the $R(G)$-module structure on 
$K_G(\T\Xcal)$, by taking  $\Ycal=point$ and $H=\{e\}$. If we take 
$\Xcal=\Ycal$ and $H=\{e\}$, the product on $K_G(\T\Xcal)$ is defined 
by 
\begin{equation}\label{eq.produit.anneau}
  \sigma_{1}\tilde{\odot}\, \sigma_2:=s_{_{\Xcal}}^*(\sigma_{1}\odot 
  \sigma_{2}) \ ,
\end{equation}  
where $s_{_{\Xcal}}:\T \Xcal\to\T \Xcal\times\T \Xcal$ is the diagonal 
map.

In the transversally elliptic case we need to be careful in the  
definition of the exterior product, since
$\T_{G\times H}(\Xcal\times \Ycal)\neq\T_{G} \Xcal\times\T_{H} \Ycal$.

\begin{defi}
Let $\sigma$ be a $H$-transversally elliptic symbol on $\T \Ycal$. 
This symbol is
called $H$-{\em transversally-good} if the characteristic set 
of $\sigma$ intersects $\T_{H}\Ycal$ in a compact subset of $\Ycal$.
\end{defi}

Recall Lemma 3.4 and Theorem 3.5 of Atiyah in \cite{Atiyah.74}. 
Let $\sigma_{1}$ be a $G$-transversally 
elliptic symbol on $\T \Xcal$, and $\sigma_{2}$ be a $H$-transversally 
elliptic symbol on $\T \Ycal$ that is $G$-equivariant. Suppose 
furthermore that $\sigma_{2}$ is $H$-{\em transversally-good}, 
then the product $\sigma_{1}\odot \sigma_{2}$ is 
$G\times H$-transversally 
elliptic. Since every class of $K_{G\times H}(\T_{H}\Ycal)$ can be 
represented by an $H$-{\em transversally-good} elliptic symbol, we 
have a multiplication
\begin{eqnarray}\label{eq:produit.transversal}
K_{G}(\T_{G}\Xcal)\times K_{G\times H}(\T_{H}\Ycal)&\longrightarrow &
K_{G\times H}(\T_{G\times H}(\Xcal\times \Ycal) )\\
(\sigma_{1},\sigma_{2})&\longmapsto &  \sigma_{1}\odot  
\sigma_{2} \ .\nonumber
\end{eqnarray}

Suppose now that the manifolds $\Xcal$ and $\Ycal$ satisfy 
the condition of Corollary \ref{hyp.indice}. So, the index maps 
$\indice_{\Xcal}^{G}$, $\indice_{\Ycal}^{G\times H}$, and  
$\indice_{\Xcal\times \Ycal}^{G\times H}$ are well defined.
According to Theorem 3.5 of \cite{Atiyah.74}, we know that

\begin{equation} \label{eq:formule.G.H.produit}
\indice_{\Xcal\times \Ycal}^{G\times H} (\sigma_{1}\odot \sigma_{2})=
\indice_{\Xcal}^{G}(\sigma_{1})\cdot\indice_{\Ycal}^{G\times H}
(\sigma_{2})
\quad{\rm in}\quad \ R^{-\infty}(G\times H)\ ,
\end{equation}

for any $\sigma_{1}\in K_{G}(\T_{G}\Xcal)$ and 
$\sigma_{2}\in K_{G\times H}(\T_{H}(\Xcal\times H))$.

\bigskip

{\em In the rest of this subsection we suppose that the subgroup
$H\subset G$ is the stabilizer of an element $\gamma\in\ggot$.
The manifold $G/H$ carries a $G$-invariant complex structure
$J_{\gamma}$ defined by the element $\gamma$: at 
$e\in G/H$, the map $J_{\gamma}(e)$ equals 
$ad(\gamma).(\sqrt{- ad(\gamma)^{2}})^{-1}$ 
on $\T_{e}(G/H)=\ggot/\hgot$.
} 

\bigskip

We recall now the definition of the map 
$r^{\gamma}_{_{G,H}}:K_{G}(\T_{G}\Xcal)\to K_{H}(\T_{H}\Xcal)$ 
introduced by Atiyah in \cite{Atiyah.74}. We consider the manifold 
$\Xcal\times G$ with two actions of $G\times H$: 
for $(g,h)\in G\times H$ and $(x,a)\in \Xcal\times G$, 
we have $(g,h).(x,a):=(g.x,gah^{-1})$ on 
$\Xcal\stackrel{1}{\times} G$, and we have $(g,h).(x,a):=(h.x,gah^{-1})$ on 
$\Xcal\stackrel{2}{\times} G$.

The map $\Theta :\Xcal\stackrel{2}{\times} G\to 
\Xcal\stackrel{1}{\times} G,\ (x,a)\mapsto (a.x,a)$ 
is $G\times H$-equivariant, and induces $\Theta^{*}:
K_{G\times H}(\T_{G\times H}(\Xcal\stackrel{1}{\times}G))\to
K_{G\times H}(\T_{G\times H}(\Xcal\stackrel{2}{\times} G))$. 
The $G$-action is free on $\Xcal\stackrel{2}{\times} G$, so the quotient map
$\pi :\Xcal\stackrel{2}{\times} G\to\Xcal$ induces
an isomorphism $\pi^{*}:K_{H}(\T_{H}\Xcal)
\to K_{G\times H}(\T_{G\times H}(\Xcal\stackrel{2}{\times}G)) $. 
We denote by $\sigma^{_\gamma}_{\ggot/\hgot}\in K_{G\times H}(\T_{H}G)$ 
the pullback of the Thom class $\Thom_{G}(G/H,J_{\gamma})\in 
K_{G}(\T(G/H))$, via the quotient map $G\to G/H$.

Consider the manifold  $\Ycal=G$ with the action of $G\times H$ 
defined by $(g,h).a=gah^{-1}$  for $a\in G$, and $(g,h)\in G\times H$. 
Since the symbol $\sigma^{_\gamma}_{\ggot/\hgot}$ is $H$-transversally 
good on $\T G$, the product by $\sigma^{_\gamma}_{\ggot/\hgot}$ 
induces, by (\ref{eq:produit.transversal}), the map 
\begin{eqnarray*}
K_{G}(\T_{G}\Xcal)&\longrightarrow &
K_{G\times H}(T_{G\times H}(\Xcal\stackrel{1}{\times}G) )\\
\sigma&\longmapsto &  \sigma\odot  
\sigma^{_\gamma}_{\ggot/\hgot} \ .
\end{eqnarray*} 

\begin{defi}[Atiyah]\label{def.G.H.restriction} 
    Let $H$ the stabilizer of 
    $\gamma\in \ggot$ in $G$. The map $r^{\gamma}_{_{G,H}}:
    K_{G}(\T_{G}\Xcal)\to K_{H}(\T_{H}\Xcal)$ is 
    defined for every $\sigma\in K_{G}(\T_{G}\Xcal)$ by
$$
r^{\gamma}_{_{G,H}}(\sigma):=(\pi^{*})^{-1}
\circ\Theta^{*}(\sigma\odot \sigma^{_\gamma}_{\ggot/\hgot})\ .
$$
\end{defi}

Theorem 4.2 in \cite{Atiyah.74} tells us that the following diagram is 
commutative

\begin{equation}\label{indice.r.G.H}
\xymatrix@C=2cm@R=10mm{
K_{G}(\T_{G}\Xcal)\ar[r]^{r^{_\gamma}_{_{G,H}}}
\ar[d]_{\indice_{\Xcal}^G} & 
K_{H}(\T_{H}\Xcal)\ar[d]^{\indice_{\Xcal}^H}\\
\fgene(G)^{G} &  \fgene(H)^{H}\ar[l]_{\indH}.
   }
\end{equation} 

\medskip

We show now a more explicit description of the map 
$r^{\gamma}_{_{G,H}}$. Consider the moment map 
$$
\mu_{_G}:\T^{*}\Xcal\to \ggot^{*}
$$
for the (canonical) Hamiltonian action of $G$ on the 
symplectic manifold $\T^{*}\Xcal$. If we identify $\T \Xcal$ with 
$\T^{*}\Xcal$ via a $G$-invariant metric, and $\ggot$ with 
$\ggot^{*}$ via a $G$-invariant scalar product, the `moment map' is a 
map $\mu_{_G}:\T \Xcal\to \ggot$ defined as follows. If 
$E^{1},\cdots,E^{l}$ is an orthonormal basis of $\ggot$, we have
$\mu_{_G}(x,v)=\sum_{i}(E^{i}_{M}(x),v)_{_{M}}E^i$ for 
$(x,v)\in\T \Xcal$. The moment map admits the decomposition 
$\mu_{_G}=\mu_{_H}+\mu_{_{G/H}}$, relative to the $H$-invariant 
orthogonal 
decomposition of the Lie algebra $\ggot=\hgot\oplus\hgot^{\perp}$. 
It is important to note that $\T_{G}\Xcal=\mu_{_G}^{-1}(0)$, 
$\T_{H}\Xcal=\mu_{_H}^{-1}(0)$, and $\T_{G}\Xcal=
\T_{H}\Xcal\cap\mu_{_{G/H}}^{-1}(0)$.

The real vector space  $\ggot/\hgot$ is endowed with the complex 
structure defined by $\gamma$. Consider over $\T \Xcal$ the 
$H$-equivariant symbol 
\begin{eqnarray*}
    \sigma^{\Xcal}_{_{G,H}} : \T \Xcal\times
    \wedge_{\C}^{even}\ggot/\hgot&
    \longrightarrow&
    \T \Xcal\times\wedge_{\C}^{odd}\ggot/\hgot\\
    (x,v;w)&\longrightarrow& (x,v;w')\ ,
\end{eqnarray*}
with $w'=Cl(\mu_{_{G/H}}(x,v)).w$. Here  $\hgot^{\perp}\simeq
\ggot/\hgot$, and $Cl(X):\wedge_{\C}\ggot/\hgot
\to\wedge_{\C}\ggot/\hgot$, $X\in\ggot/\hgot$, denotes the Clifford 
action. This symbol has $\mu_{_{G/H}}^{-1}(0)$ for characteristic set.
For any symbol $\sigma$ over $\T\Xcal$, with characteristic set 
${\rm Char}(\sigma)$, the product $\sigma\,\tilde{\odot}\,
\sigma^{\Xcal}_{_{G,H}}$,
defined at (\ref{eq.produit.anneau}),
is a symbol over $\T \Xcal$ with  characteristic set 
${\rm Char}(\sigma\,\tilde{\odot}\,\sigma^{\Xcal}_{_{G,H}})=
{\rm Char}(\sigma)\cap\mu_{_{G/H}}^{-1}(0)$. Then, if $\sigma$
is a $G$-transversally elliptic symbol over $\T \Xcal$, the product
$\sigma\,\tilde{\odot}\,\sigma^{\Xcal}_{_{G,H}}$ is a 
$H$-transversally elliptic symbol.

\medskip

\begin{prop}\label{prop.restriction.bis}
    The map $r^{\gamma}_{_{G,H}}:K_{G}(\T_{G}\Xcal)\to 
    K_{H}(\T_{H}\Xcal)$ has the 
    following equivalent definition: for every $\sigma\in 
    K_{G}(\T_{G}\Xcal)$
    $$
    r^{\gamma}_{_{G,H}}(\sigma)=\sigma\,\tilde{\odot}\,\sigma^{\Xcal}_{_{G,H}}
    \quad {\rm in} \quad K_{H}(\T_{H}\Xcal).
    $$
\end{prop}

\medskip

{\em Proof} : We have to show that for every $\sigma\in K_{G}(\T_{G}\Xcal)$,
$\sigma\tilde{\odot}\sigma^{\Xcal}_{_{G,H}}=$ \break 
$(\pi^{*})^{-1}\circ\Theta^{*}(\sigma\odot\sigma^{_\gamma}_{\ggot/\hgot})$
in $K_{H}(\T_{H}\Xcal)$. Let $p_{_G}:\T G\to G$  
 and $p_{\Xcal}:\T \Xcal\to \Xcal$ be the canonical projections. 
The symbol $\sigma^{_\gamma}_{\ggot/\hgot}:p^{*}_{_G}(G\times
\wedge_{\C}^{even}\ggot/\hgot)\to p^{*}_{_G}(G\times\wedge_{\C}^{odd}
\ggot/\hgot)$ is defined by 
$\sigma^{_\gamma}_{\ggot/\hgot}(a, Z)=Cl(Z_{\ggot/\hgot})$ for $(a,Z)\in\T 
G\simeq G\times\ggot$, where $Z_{\ggot/\hgot}$ is the 
$\ggot/\hgot$-component of $Z\in \ggot$.

Consider $\sigma:p_{\Xcal}^{*}E_{0}\to p_{\Xcal}^{*}E_{1}$, a 
$G$-transversally elliptic symbol on $\T \Xcal$, where $E_{0}, E_{1}$ are 
$G$-complex vector bundles over $\Xcal$. The product $\sigma\odot
\sigma^{_\gamma}_{\ggot/\hgot}$ acts on the bundles
$p_{\Xcal}^{*}E_{\bullet}\otimes p^{*}_{G}(G\times\wedge_{\C}^{\bullet}
\ggot/\hgot)$ at $(x,v;a,Z)\in \T(\Xcal\times G)$ by
$$
\sigma(x,v)\odot Cl(Z_{\ggot/\hgot}) .
$$
The pullback $\sigma_{o}:=\Theta^{*}(\sigma\odot \sigma_{\ggot/\hgot})$ acts
on the bundle $G\times (p_{\Xcal}^{*}E_{\bullet}\otimes \wedge_{\C}^{\bullet}
\ggot/\hgot)$ (here we identify $\T(\Xcal\times G)$ with 
$G\times(\ggot\oplus\T \Xcal)$).
At $(x,v;a,Z) \in \T(\Xcal\times G)$ we have 
$$
\sigma_{o}(x,v;a,Z)=\sigma\odot 
\sigma^{_\gamma}_{\ggot/\hgot}(a.x,v';a,Z'),\quad {\rm with}
$$
$(v',Z')=\left([\T_{(x,a)}\Theta]^{*}\right)^{-1}(v,Z)$. Here
$\T_{(x,a)}\Theta:\T_{(x,a)}(\Xcal\times G)\to\T_{(a.x,a)}(\Xcal\times G)$ is the 
tangent map of $\Theta$ at $(x,a)$, and $[\T_{(x,a)}\Theta]^{*}:\T_{(a.x,a)}
(\Xcal\times G)\to\T_{(x,a)}(\Xcal\times G)$ its transpose. A small 
computation shows that $Z'=Z +\mu_{_{G}}(v)$ and $v'=a.v$. Finally, we get
$$
\sigma_{o}(x,v;a,Z)=\sigma(a.x,a.v)\odot Cl(Z_{\ggot/\hgot}
+\mu_{_{G/H}}(v)) .
$$
Hence, the symbol $(\pi^{*})^{-1}(\sigma_{o})$ acts on the bundle 
$p_{\Xcal}^{*}E_{\bullet}\otimes \wedge_{\C}^{\bullet}\ggot/\hgot$ by
$$
(\pi^{*})^{-1}(\sigma_{o})(x,v)=\sigma(x,v)\odot Cl(\mu_{_{G/H}}(v)).
$$
$\Box$

\medskip

For any $G$-invariant function
$\phi\in\f(G)^{G}$, the Weyl integration formula can be 
written\footnote{See Remark \ref{wedge.C-wedge.R}.} 
\begin{equation}\label{eq.Weyl.integration}
\phi=\indH\left(\phi_{\vert H}\wedge^{\bullet}_{\C}\ggot/\hgot\right)\ 
{\rm in}\ \fgene(G)^{G}\ .
\end{equation}
where $\phi_{\vert H}\in\f(H)^{H}$ is the restriction to $H=G_{\gamma}$.  
Equality (\ref{eq.Weyl.integration}) remains true for any
$\phi\in\fgene(G)^{G}$ that admits a restriction to $H$. 

\medskip

\begin{lem}\label{restriction.function.G.H}
    Let $\sigma$ be a $G$-transversally elliptic symbol. Suppose 
    furthermore that $\sigma$ is $H$-transversally elliptic. This
    symbol defines two classes $\sigma\in K_{G}(\T_{G}\Xcal)$ and
    $\sigma_{\vert H}\in K_{H}(\T_{H}\Xcal)$ with 
    the relation\footnote{Here we note  
    $\sigma_{\vert H}\otimes\wedge^{\bullet}_{\C}\ggot/\hgot$ for the 
    difference $\sigma_{\vert H}\otimes\wedge^{even}_{\C}\ggot/\hgot\,-\,
    \sigma_{\vert H}\otimes\wedge^{odd}_{\C}\ggot/\hgot$.}
    $r^{\gamma}_{_{G,H}}(\sigma)=
    \sigma_{\vert H}\otimes\wedge^{\bullet}_{\C}\ggot/\hgot$.
    Hence for the generalized character
    $\indice^{G}_{\Xcal}(\sigma)\in R^{-\infty}(G)$ we have a  
    `Weyl integration' formula
    \begin{equation}\label{eq.restriction.function.G.H}
    \indice^{G}_{\Xcal}(\sigma)=
    \indH\left(\indice^{H}_{\Xcal}(\sigma_{\vert H})
    \wedge^{\bullet}_{\C}\ggot/\hgot\right)\ .
    \end{equation}
\end{lem}    

{\em Proof :} If $\sigma$ is $H$-transversally elliptic, the
symbol $(x,v)\to \sigma(x,v)\odot Cl(\mu_{_{G/H}}(v))$ is homotopic to
$(x,v)\to \sigma(x,v)\odot Cl(0)$ in $K_{H}(\T_{H}\Xcal)$. Hence 
$\sigma_{\vert H}\odot\sigma^{\Xcal}_{_{G,H}}=
\sigma_{\vert H}\otimes\wedge_{\C}^{\bullet}\ggot/\hgot$ in 
$K_{H}(\T_{H}\Xcal)$. (\ref{eq.restriction.function.G.H}) 
follows from the diagram (\ref{indice.r.G.H}).
$\Box$

\begin{coro}\label{coro.restriction.G.H}
Let $\sigma$ be a $G$-transversally elliptic symbol which 
furthermore is $H$-transversally elliptic, and let 
$\phi\in\fgene(G)^{G}$ which admits a restriction to $H$. We have
$$
\phi =\indice^{G}_{\Xcal}(\sigma)\Longleftrightarrow 
\phi_{\vert H} =\indice^{H}_{\Xcal}(\sigma_{\vert H})\ .
$$
\end{coro}

In fact, if we come back to the definition of the analytic index given by 
Atiyah \cite{Atiyah.74}, one can show the following stronger result.  
If $\sigma$ be a $G$-transversally elliptic symbol which is also 
$H$-transversally elliptic, then 
$\indice^{G}_{\Xcal}(\sigma)\in\fgene(G)^G$ admits a {\em restriction} to
$H$ equal to $\indice^{H}_{\Xcal}(\sigma_{\vert H})\in\fgene(H)^H$.

\medskip

\section{Localization - The general procedure}\label{sec.general.procedure}

\medskip

We recall briefly the notations. Let $(M,J,G)$ be a compact $G$-manifold
provided with a $G$-invariant almost complex structure. We denote by 
$RR^{^{G,J}} : K_{G}(M)\to R(G)$ (or simply $RR^{^{G}}$), the corresponding 
quantization map. We choose a $G$-invariant Riemannian
metric $(.,.)_{_{M}}$ on $M$. We define in this section a general procedure 
to localize the quantization map  through the use of a 
$G$-equivariant vector field $\lambda$. This idea of localization 
goes back, when $G$ is a circle group, to Atiyah \cite{Atiyah.74}
(see Lecture 6) and Vergne \cite{Vergne96} (see part II).

We denote by $\Phi_{\lambda}:M\to\ggot^{*}$ the map defined by 
$\langle\Phi_{\lambda}(m),X\rangle:= (\lambda_{m}, X_{M}\vert_{m})_{_{M}}$ 
for $X\in\ggot$. We denote by $\sigma^{_E}(m,v),\ (m,v)\in \T M$ 
the elliptic symbol associated to 
$\Thom_{G}(M)\otimes p^*(E)$ for $E\in K_{G}(M)$ (see section 
\ref{sec.quantization}). 

Let $\sigma^{_E}_{1}$ be the following $G$-invariant elliptic symbol
\begin{equation}
    \sigma^{_E}_{1}(m,v):=\sigma^{_E}(m,v-\lambda_{m}),\quad (m,v)\in \T M.
    \label{eq:sigma.1}
\end{equation}

The symbol $\sigma^{_E}_{1}$ is obviously homotopic to $\sigma^{_E}$, 
so they define the same class in $K_{G}(\T M)$. The characteristic 
set $\Char(\sigma^{_E})$ is $M\subset \T M$, but we see easily
that $\Char(\sigma^{_E}_{1})$  is equal to the graph of the vector field
$\lambda$, and 
$$
\Char(\sigma^{_E}_{1})\cap\T_{G}M=\left\{(m,\lambda_{m})\in\T M, \ \ m\in 
\{\Phi_{\lambda}=0\} \right\}.
$$

We will now decompose the elliptic symbol $\sigma^{_E}_{1}$ in $K_{G}(\T_{G}M)$
near 
$$
C_{\lambda}:=\{\Phi_{\lambda}=0\}\ .
$$

If a $G$-invariant subset $C$ is a union of  
{\em connected components} of $C_{\lambda}$
there exists a $G$-invariant open neighbourhood $\Ucal^{c} \subset M$ 
of $C$  such that $\Ucal^c \cap C_{\lambda}=C$ and 
$\partial\Ucal^c \cap C_{\lambda}=\emptyset$. We associate to the 
subset $C$ the symbol $\sigma^{_E}_{^{C}}:=\sigma^{_E}_{1}\vert_{\Ucal^c}\in 
K_{G}(\T_{G}\Ucal^c)$ 
which is the  restriction of $\sigma^{_E}_{1}$ to $\T\Ucal^c$. It is 
well defined since $\Char(\sigma^{_E}_{1}\vert_{\Ucal^c})\cap\T_{G}
\Ucal^{c} =\{(m,\lambda_{m})\in\T M,\ m\in C\}$ is compact.

\medskip

\begin{prop}\label{prop.localisation}
   Let $C^a,a\in A$, be a finite collection of disjoint $G$-invariant subsets of
   $C_{\lambda}$, each of them being  a union of connected components 
   of $C_{\lambda}$, and let $\sigma^{_E}_{^{C^a}}\in K_{G}(\T_{G}\Ucal^a)$
   be the localized symbols. If $C_{\lambda}=\cup_{a}C^a$, we have
   $$
   \sigma^{_E}=\sum_{a\in A} i^{a}_{*}(\sigma^{_E}_{^{C^a}})\quad 
   {\rm in}\quad K_{G}(\T_{G}M),
   $$
   where $i^{a}:\Ucal^a \hookrightarrow M$ is the inclusion and 
   $i^a_{*}:K_{G}(\T_{G}\Ucal^a)\to K_{G}(\T_{G}M)$ is the 
   corresponding direct image.
\end{prop}

\medskip

{\em Proof :} This is a consequence of the property of excision 
(see subsection \ref{subsec.excision}). We 
consider disjoint neighbourhoods $\Ucal^a$ of $C^a$, and take 
$i:\Ucal=\cup_{a}\Ucal^a \hookrightarrow M$. Let $\chi_{a}\in\f(M)^{G}$ be a 
test function (i.e.  $0\leq \chi_{a}\leq 1$) with compact support on $\Ucal^a$ 
such that $\chi_{a}(m)\neq 0$ if $m\in C^a$. Then the function 
$\chi:=\sum_{a}\chi_{a}$ is a $G$-invariant test function with support in  
$\Ucal$ such that $\chi$ never vanishes on $C_{\lambda}$. 

Using the $G$-equivariant symbol 
$\sigma^{_E}_{^{\chi}}(m,v)
:=\sigma^{_E}(m,\chi(m)v-\lambda_{m}),\ (m,v)\in \T M$,
we prove the following :

\noindent i) the symbol $\sigma^{_E}_{^{\chi}}$ is $G$-transversally 
elliptic and $\Char(\sigma^{_E}_{^{\chi}})\subset\T M\vert_{\Ucal}$,

\noindent ii) the symbols $\sigma^{_E}_{^{\chi}}$ and  $\sigma^{_E}_{1}$
are equal in $K_{G}(\T_{G}M)$, and

\noindent iii) the restrictions $\sigma^{_E}_{^{\chi}}\vert_{\Ucal}$ and 
$\sigma^{_E}_{1}\vert_{\Ucal}$ are equal in $K_{G}(\T_{G}\Ucal)$.

With Point i) we can apply the excision property to 
$\sigma^{_E}_{^{\chi}}$, hence $\sigma^{_E}_{^{\chi}}=
i_{*}(\sigma^{_E}_{^{\chi}}\vert_{\Ucal})$.  By ii) and iii), the last 
equality gives $\sigma^{_E}_{1}=i_{*}(\sigma^{_E}_{1}\vert_{\Ucal})=
\sum_{a} i^{a}_{*}(\sigma^{_E}_{^{C^a}})$.   

{\em Proof of i)}. The point $(m,v)$ belongs to $\Char(\sigma^{_E}_{^{\chi}})$
if and only if $\chi(m)v=\lambda_{m}(*)$. If $m$ is not included in $\Ucal$, 
we have $\chi(m)=0$ and the equality $(*)$ becomes   
$\lambda_{m}=0$. But $\{\lambda =0\}\subset C_{\lambda}\subset 
\Ucal$, thus $\Char(\sigma^{_E}_{^{\chi}})\subset \T M\vert_{\Ucal}$.
The point $(m,v)$ belongs to $\Char(\sigma^{_E}_{^{\chi}})\cap \T_{G}M$
 if and only if $\chi(m)v=\lambda_{m}$ and $v$ is orthogonal to the 
 $G$-orbit in $m$. This imposes $m\in C_{\lambda}$, and finally we 
 see that $\Char(\sigma^{_E}_{^{\chi}})\cap \T_{G} M\simeq C_{\lambda}$
is compact because the function $\chi$ never vanishes on $C_{\lambda}$.

{\em Proof of ii)}. We consider the symbols $\sigma^{_E}_{t},\, t\in[0,1]$ 
defined by
$$
\sigma^{_E}_{t}(m,v)=\sigma^{_E}(m,(t+(1-t)\chi(m))v-\lambda_{m}).
$$
We see as above that 
$\sigma^{_E}_{t}$ is an homotopy of 
$G$-transversally elliptic symbols on $\T M$.
 
{\em Proof of iii)}. Here we use the homotopy 
$\sigma^{_E}_{t}\vert_{\Ucal},\ t\in[0,1]$. 

$\Box$

Because $RR^{^G}(M,E)=\indice_{M}^{G}(\sigma^{_E})\in R(G)$, we obtain from
 Proposition \ref{prop.localisation} the following decomposition 
\begin{equation}\label{eq:decomposition.RR}
    RR^{^G}(M,E)= \sum_{a\in A}\indice_{\Ucal^a}^{G}(
    \sigma^{_E}_{^{C^a}})\quad {\rm in}\quad R^{-\infty}(G).    
\end{equation}

The rest of this article is devoted to the description, in some 
particular cases, of the Riemann-Roch character localized near $C^a$:
\begin{eqnarray}\label{eq:RR.localise}
 RR_{C^a}^{^G}(M,-) \ :\ K_{G}(M)&\longrightarrow& R^{-\infty}(G)\\
    E\ \  &\longmapsto & 
    \indice_{\Ucal^a}^{G}(\sigma^{_E}_{^{C^a}}).\nonumber   
\end{eqnarray}

\medskip

\section{Localization on \protect $M^{\beta}$}\label{sec.loc.M.beta}

\medskip

Let $(M,J,G)$ be a compact $G$-manifold provided with a $G$-invariant 
almost complex structure. Let $\beta$ be an element in the {\em center} 
of the Lie algebra of $G$, and consider the $G$-invariant vector field 
$\lambda:=\beta_{M}$ generated by the infinitesimal action of $\beta$.
In this case we have obviously
$$
\{\Phi_{\beta_{M}}=0\}=\{\beta_{M}=0\}=M^{\beta}\ .
$$
In this section, we compute the localization of the quantization map
on the submanifold $M^{\beta}$ following the technique explained in  
section \ref{sec.general.procedure}. We first need to understand 
the case of a vector space. 

The principal results of this section, i.e. Proposition 
\ref{prop.indice.thom.beta.V} and Theorem \ref{th.localisation.pt.fixe} 
were obtained by Vergne \cite{Vergne96}[Part II], in the Spin case for 
an action of the circle group.

\subsection{Action on a vector space}
 
Let $(V,q,J)$ be a real vector space equipped with a complex 
structure $J$ and an euclidean metric $q$ such that $J\in O(q)$.
Suppose that a compact Lie group $G$ acts on $(V,q,J)$ in a unitary 
way, and that there exists $\beta$ in the center of $\ggot$ 
such that
$$
V^{\beta}=\{0\}.
$$
We denote by $\tore_{\beta}$ the torus generated by $\exp(t.\beta), 
t\in\R$, and $\tgot_{\beta}$ its Lie algebra.

The complex $\Thom_{G}(V,J)$ does not define an element in 
$K_{G}(\T V)$ because its characteristic set is $V$.
\begin{defi}
    Let $\Thom^{\beta}_{G}(V)\in K_{G}(\T_{G}V)$ be the 
    $G$-transversally\footnote{One can verify that $\Char(\Thom^{\beta}_{G}(V))\cap 
    \T_{G}V=\{(0,0)\}$.} elliptic complex defined by
    $$
    \Thom^{\beta}_{G}(V)(x,v):=\Thom_{G}(V)(x,v- \beta_{V}(x)) \quad
    {\rm for}\quad (x,v)\in \T V.
    $$
    \label{def.thom.beta.V}
\end{defi}

Before computing the index of $\Thom^{\beta}_{G}(V)$ explicitely, we 
compare it with the pushforward $j_{!}(\C)\in K_{G}(\T V)$ where
$j: \{0\}\croc V$ is the inclusion and $\C\to\{0\}$ is the trivial 
line bundle. Recall that $\indice_{V}^G(j_{!}(\C))=1$.

We denote by $\overline{V}$ the real vector space $V$ endowed with 
the complex structure $-J$, and 
$\wedge_{\C}^{\bullet} \overline{V}:=\wedge_{\C}^{even} 
\overline{V}-\wedge_{\C}^{odd} \overline{V} $ the 
corresponding element in $R(G)$.  

\begin{lem}\label{lem.indice.wedge.V.inverse}
We have $\wedge_{\C}^{\bullet}\overline{V}\, .\, \Thom^{\beta}_{G}(V)= 
j_{!}(\C)$ in $K_{G}(\T_{G}V)$, hence 
$$
\wedge_{\C}^{\bullet} \overline{V}.\,\indice_{V}^G(\Thom^{\beta}_{G}(V))=1\quad
{\rm in} \quad R^{-\infty}(G).
$$
\end{lem}

{\em Proof} : The class $j_{!}(\C)$ is represented by the symbol
$\sigma_{o}:\T V\times \wedge^{even}_{\C}(V\otimes\C)\to
\T V\times \wedge^{odd}_{\C}(V\otimes\C),\ (x,v,w)\mapsto
(x,v,Cl(x+\imath v).w)$. If we use the following
isomorphism of complex $G$-vector spaces
\begin{eqnarray*}
 V\otimes\C &\longrightarrow &V\oplus \overline{V}\\
 x +\imath v &\longmapsto &(v-J(x), v+J(x))\ ,
\end{eqnarray*}
we can write $\sigma_{o}=\sigma_{-}\odot\sigma_{+}$, where
the symbols\footnote{$V_{+}=V$ and $V_{-}=\overline{V}$.}
$\sigma_{\pm}$ act on $\T V\times \wedge^{\bullet}_{\C}V_{\pm}$
through the Clifford maps $\sigma_{\pm}(x,v)=Cl(v \mp J(x))$.
Finally we see that the following $G$-transversally elliptic 
symbols on $\T V$ are homotopic
\begin{eqnarray*}
 Cl(v + J(x))&\odot& Cl(v - J(x))\\
 Cl(v + J(x))&\odot& Cl(v - \beta_{V}(x))\\
 Cl(0)&\odot& Cl(v - \beta_{V}(x)) \ .
\end{eqnarray*}
The Lemma is proved since $(x,v)\to Cl(0)\odot Cl(v - \beta_{V}(x))$
represents the class $\wedge_{\C}^{\bullet} \overline{V}\, .\, 
\Thom^{\beta}_{G}(V)$ in $K_{G}(\T_{G}V)$. $\Box$

\medskip

We compute now the index of 
$\Thom^{\beta}_{G}(V)$. For $\alpha\in \tgot_{\beta}^{*}$, we define 
the $G$-invariant subspaces\footnote{We denote by $z\cdot v:=x.v+y.J(v),\ 
z=x+\imath y\in \C$, the action of $\C$ on the complex vector space 
$(V,J)$, and $zw=v\otimes zz',\ w=v\otimes z'\in V\otimes\C$ 
the canonical action of $\C$ on $V\otimes\C$.}
$V(\alpha):=\{v\in V,\ \rho(\exp X)(v)=
e^{\imath\langle\alpha,X\rangle}\cdot v,\ \forall X\in 
\tgot_{\beta}\}$, and 
$(V\otimes\C)(\alpha):=\{v\in V\otimes\C,\ \rho(\exp X)(v)=
e^{\imath\langle\alpha,X\rangle}v,\ \forall X\in 
\tgot_{\beta}\}$.

An element $\alpha \in \tgot_{\beta}^{*}$, is called a weight for
the action of $\tore_{\beta}$ on $(V,J)$ (resp. on $V\otimes\C$) 
if $V(\alpha)\neq 0$ (resp.  $(V\otimes\C)(\alpha)\neq 0$). We denote by 
$\Delta(\tore_{\beta},V)$ (resp. $\Delta(\tore_{\beta},V\otimes\C)$) 
the set of weights for the action of $\tore_{\beta}$ on $V$ 
(resp. $V\otimes\C$). We shall note that 
$\Delta(\tore_{\beta},V\otimes\C)=\Delta(\tore_{\beta},V)
\cup -\Delta(\tore_{\beta},V)$.

\begin{defi}
We denote by $V^{+,\beta}$ the following $G$-stable subspace of $V$
$$
V^{+,\beta}:=\sum_{\alpha\in \Delta_{+}(\tore_{\beta},V)}V(\alpha)\ ,
$$
where $\Delta_{+}(\tore_{\beta},V)=\{\alpha\in\Delta(\tore_{\beta},V),
\ \langle\alpha,\beta\rangle > 0\}$. In the same way, we denote by 
$(V\otimes\C)^{+,\beta}$ the following $G$-stable subspace of 
$V\otimes\C$:
$(V\otimes\C)^{+,\beta}
:=\sum_{\alpha\in \Delta_{+}(\tore_{\beta},V\otimes\C)}
(V\otimes\C)(\alpha)$, where $\Delta_{+}(\tore_{\beta},V\otimes\C)
=\{\alpha\in\Delta(\tore_{\beta},V\otimes\C),\ 
\langle\alpha,\beta\rangle > 0\}$.
\end{defi}

For any representation $W$ of $G$, we denote by $\det W$ the 
representation $\wedge_{\C}^{max}W$. In the same way, if $W\to M$ is 
a $G$ complex vector bundle we denote by $\det W$ the corresponding line 
bundle. 

\begin{prop}\label{prop.indice.thom.beta.V}
    We have the following equality in $R^{-\infty}(G)$ :
    $$
    \indice^{G}_{V}(\Thom^{\beta}_{G}(V))=(-1)^{\dim_{\C}V^{+,\beta}}\
    \det V^{+,\beta}\otimes\sum_{k\in \N} 
    S^k((V\otimes\C)^{+,\beta})\ ,
    $$
    where $S^k((V\otimes\C)^{+,\beta})$ is the $k$-th symmetric 
    product over $\C$ of $(V\otimes\C)^{+,\beta}$.
\end{prop}

Proposition \ref{prop.indice.thom.beta.V} and Lemma 
\ref{lem.indice.wedge.V.inverse} give the two important properties
of the generalized function 
$\chi:=\indice_{G}^{V}(\Thom^{\beta}_{G}(V))$.
First $\chi$ is an inverse, in $R^{-\infty}(G)$, of the function 
$g\in G\to \det_{V}^{\C}(1-g^{-1})$ which is  the trace of 
the (virtual) representation 
$\wedge_{\C}^{\bullet}\overline{V}$. Second, the decomposition 
of $\chi$ into irreducible characters of $G$ is of the form 
$\chi=\sum_{\lambda} m_{\lambda}\chi_{_{\lambda}}^{_{G}}$ with 
$m_{\lambda}\neq 0
\Longrightarrow \langle \lambda,\beta\rangle \geq 0$.

\begin{defi}\label{wedge.V.inverse}
    For any $R(G)$-module $A$, we denote by 
    $A\,\widehat{\otimes}\,R(\tore_{\beta})$, the $R(G)\otimes 
    R(\tore_{\beta})$-module formed by the infinite formal sums
    $\sum_{\alpha} E_{\alpha}\, h^{\alpha}$ taken over the set 
    of weights of $\tore_{\beta}$,  where $E_{\alpha}\in A$ 
    for every $\alpha$.

    We denote by $\left[\wedge_{\C}^{\bullet}\overline{V}\,
    \right]^{-1}_{\beta}$ 
    the infinite sum $(-1)^{r}\,\det V^{+,\beta}\otimes\sum_{k\in \N} 
    S^k((V\otimes\C)^{+,\beta})$, with $r=\dim_{\C}V^{+,\beta}$.
    It can be considered either as an element of $R^{-\infty}(G)$, 
    $R(G)\,\widehat{\otimes}\, R(\tore_{\beta})$, or
    $R^{-\infty}(\tore_{\beta})$.
    
    Let $\Vcal\to \Xcal$ be a $G$-complex vector bundle such that 
    $\Vcal^{\beta}=\Xcal$. The torus $\tore_{\beta}$ acts on the 
    fibers  of $\Vcal\to \Xcal$, so we can polarize the 
    $\tore_{\beta}$-weights and  define the vector bundles 
    $\Vcal^{+,\beta}$ and $(\Vcal\otimes\C)^{+,\beta}$. 
    In this case, the infinite sum 
    $\left[\wedge_{\C}^{\bullet}\overline{\Vcal}\,
    \right]^{-1}_{\beta}:=(-1)^{\dim_{\C}\Vcal^{+,\beta}}\,
    \det \Vcal^{+,\beta}\otimes\sum_{k\in \N} 
    S^k((\Vcal\otimes\C)^{+,\beta})$ 
    is an inverse of $\wedge_{\C}^{\bullet}\overline{\Vcal}$ 
    in \break
    $K_G(\Xcal)\,\widehat{\otimes}\, R(\tore_{\beta})$.
\end{defi}

The rest of this subsection is devoted to the proof of Proposition 
\ref{prop.indice.thom.beta.V}. The case $V^{+,\beta}=V$ or 
$V^{+,\beta}=\{0\}$ is considered by  Atiyah 
\cite{Atiyah.74} (see Lecture 6) and Vergne 
\cite{Vergne96} (see Lemma 6, Part II).

Let $H$ be a maximal torus of $G$ containing $\tore_{\beta}$. 
The symbol $\Thom^{\beta}_{G}(V)$ is also $H$-transversally 
elliptic and let $\Thom^{\beta}_{H}(V)$ be the corresponding 
class in $K_{H}(\T_{H}V)$. 
Following Corollary \ref{coro.restriction.G.H}, we can reduce the 
proof of Proposition \ref{prop.indice.thom.beta.V} to the case where 
the group $G$ is equal to the torus $H$.

\medskip

\underline{Proof of Th. \ref{prop.indice.thom.beta.V} for a
torus action.}

\medskip

We first recall the index theorem proved by Atiyah in Lecture 6 of 
\cite{Atiyah.74}. Let $\tore_{m}$ the circle
group act  on $\C$ with the representation $t^m,\ m>0$. We have two classes 
$\Thom^{\pm}_{\tore_{m}}(\C)\in K_{\tore_{m}}(\T_{\tore_{m}}(\C))$ 
that correspond
respectively to $\beta=\pm \imath\in Lie(S^{1})$. Atiyah denotes these
elements $\overline{\partial}^{\pm}$.

\begin{lem}[Atiyah]\label{lem.atiyah.1}
    We have, for $m>0$, the following equalities in $R^{-\infty}(\tore_{m})$:
    \begin{eqnarray*}
        \indice^{\tore_{m}}_{\C}(\Thom^{+}_{\tore_{m}}(\C))
        &= \left[\frac{1}{1-t^{-m}}\right]^{+}=&-t^m.\sum_{k\in\N}(t^m)^k\\
        \indice^{\tore_{m}}_{\C}(\Thom^{-}_{\tore_{m}}(\C))
        &=\left[\frac{1}{1-t^{-m}}\right]^{-}=&\sum_{k\in\N}(t^{-m})^k\ .
    \end{eqnarray*}     
\end{lem}
Here we follow the notation of Atiyah: $[\frac{1}{1-t^{-m}}]^{+}$  and 
$[\frac{1}{1-t^{-m}}]^{-}$ are the 
Laurent expansions of the meromorphic function 
$t\in\C\to\frac{1}{1-t^{-m}}$ around $t=0$ and $t=\infty$ respectively.

 From this Lemma we can compute the index of 
 $\Thom^{\pm}_{\tore_{m}}(\C)$ when $m<0$. Suppose $m<0$
 and consider the morphism $\kappa : \tore_{m}\to\tore_{\vert m\vert}, 
 t\to t^{-1}$. Using the induced morphism $\kappa^{*}: K_{\tore_{\vert m\vert}}
 (\T_{\tore_{\vert m\vert}}(\C))\to 
 K_{\tore_{m}}(\T_{\tore_{m}}(\C))$, we see that 
 $\kappa^{*}(\Thom^{\pm}_{\tore_{\vert m\vert}}(\C))= 
 \Thom^{\mp}_{\tore_{m}}(\C)$. This gives 
$\indice^{\tore_{m}}_{\C}(\Thom^{+}_{\tore_{m}}(\C))$ $=
\kappa^*(\sum_{k\in\N}(t^{-\vert m\vert})^k)=
 \sum_{k\in\N}(t^{-m})^k$ and 
 $\indice^{\tore_{m}}_{\C}(\Thom^{-}_{\tore_{m}}(\C))=
 \kappa^*(-t^{\vert m\vert}.\sum_{k\in\N}(t^{\vert m\vert})^k)
=-t^{m}\sum_{k\in\N}(t^{m})^k$.

We can summarize these different cases as follows.
\begin{lem}\label{lem.atiyah.}
Let $\tore_{\alpha}$ the circle group act on $\C$ with the 
representation $t\to t^{\alpha}$ for $\alpha\in\Z\setminus\{0\}$.
Let $\beta \in Lie(\tore_{\alpha})\simeq \R$ a non-zero element. 
We have the following equalities in $R^{-\infty}(\tore_{\alpha})$:
$$
\indice^{\tore_{\alpha}}_{\C}
\left(\Thom^{\beta}_{\tore_{\alpha}}(\C)\right)(t)
= \left[\frac{1}{1-u^{-1}}\right]^{\esp}_{u=t^{\alpha}}\ ,
$$
where $\esp$ is the sign of $\langle\alpha,\beta\rangle$.
\end{lem}

We decompose now the vector space $V$ into an orthogonal sum
$V=\oplus_{i\in I}\C_{\alpha_{i}}$, where $\C_{\alpha_{i}}$ is
a $H$-stable subspace of dimension 1 over $\C$ equipped with the 
representation $t\in H\to t^{\alpha_{i}}\in \C$. Here the set $I$ parametrizes 
the weights for the action of $H$ on $V$, counted with their 
multiplicities. 
Consider the circle group $\tore_{i}$ with the trivial action on 
$\oplus_{k\neq i}\C_{\alpha_{k}}$ and with the canonical action on 
$\C_{\alpha_{i}}$.
We consider $V$ equipped with the action of 
$H\times\Pi_{k}\tore_{k}$. The symbol $\Thom_{H}^{\beta}(V)$ is 
$H\times\Pi_{k}\tore_{k}$-equivariant and is either 
$H$-transversally elliptic, $H\times\Pi_{k}\tore_{k}$-transversally 
elliptic (we denote by $\sigma_{B}$ the corresponding class), or 
$\Pi_{k}\tore_{k}$-transversally elliptic (we denote by $\sigma_{A}$ 
the corresponding class). We have the following canonical morphisms :

\begin{eqnarray}\label{morphismes.1.2}
K_{H}(\T_{H}V)\longleftarrow 
&K_{H\times\Pi_{k}\tore_{k}}(\T_{H}V)&\longrightarrow 
K_{H\times\Pi_{k}\tore_{k}}(\T_{H\times\Pi_{k}\tore_{k}}V)\\
\Thom^{\beta}_{H}(V)
\longleftarrow &\sigma_{B_{1}}&
\longrightarrow \sigma_{B}\ ,\nonumber
\end{eqnarray}    

\begin{eqnarray*}
K_{H\times\Pi_{k}\tore_{k}}(\T_{H\times\Pi_{k}\tore_{k}}V)
\leftarrow& K_{H\times\Pi_{k}\tore_{k}}(\T_{\Pi_{k}\tore_{k}}V)&
\rightarrow K_{\Pi_{k}\tore_{k}}(\T_{\Pi_{k}\tore_{k}}V)\\
\sigma_{B}\leftarrow &\sigma_{B_{2}}&\rightarrow \sigma_{A} \ .
\end{eqnarray*}

We consider the following characters:

\noindent - $\phi(t)\in R^{-\infty}(H)$ the $H$-index of 
$\Thom_{H}^{\beta}(V)$,

\noindent - $\phi_{B}(t,t_{1},\cdots,t_{l})\in 
R^{-\infty}(H\times\Pi_{k}\tore_{k})$ the $H\times\Pi_{k}\tore_{k}$-index 
of $\sigma_{B}$ (the same for $\sigma_{B_{1}}$ and $\sigma_{B_{2}}$).

\noindent - $\phi_{A}(t_{1},\cdots,t_{l})\in R^{-\infty}(\Pi_{k}\tore_{k})$ 
the $\Pi_{k}\tore_{k}$-index of $\sigma_{A}$.

They satisfy the relations 

\noindent i) $\phi(t)=\phi_{B}(t,1,\cdots,1)$ and 
$\phi_{B}(1,t_{1},\cdots,t_{l})=\phi_{A}(t_{1},\cdots,t_{l}).$

\noindent ii) $\phi_{B}(tu,t_{1}u^{-\alpha_{1}},\cdots,t_{l}u^{-\alpha_{1}})=
\phi_{B}(t,t_{1},\cdots,t_{l})$, for all $u\in H$.

Point i) is a consequence of the morphisms (\ref{morphismes.1.2}). 
Point ii) follows from the fact 
that the elements $(u,u^{-\alpha_{1}},\cdots,u^{-\alpha_{l}}),\ u\in 
H$ act trivially on $V$.

The symbol $\sigma_{A}$ can be expressed through the map
\begin{eqnarray*}
K_{\tore_{1}}(\T_{\tore_{1}}\C_{\alpha_{1}})
\times K_{\tore_{2}}(\T_{\tore_{2}}\C_{\alpha_{2}})\times\cdots
\times K_{\tore_{l}}(\T_{\tore_{l}}\C_{\alpha_{l}})&\longrightarrow&
K_{\Pi_{k}\tore_{k}}(\T_{\Pi_{k}\tore_{k}}V)\\
(\sigma_{1},\sigma_{2},\cdots,\sigma_{l})\longmapsto
\sigma_{1}\odot\sigma_{2}\odot\cdots\odot\sigma_{l}\ .
\end{eqnarray*}

Here we have $\sigma_A=\odot_{k=1}^l\Thom_{\tore_{k}}^{\esp_{k}}
(\C_{\alpha_{k}})$ in $K_{\Pi_{k}\tore_{k}}(\T_{\Pi_{k}\tore_{k}}V)$, 
where $\esp_{k}$ is the sign of  $\langle\alpha_{k},\beta\rangle$. 
Finally, we get 
\begin{eqnarray*}
\phi(u)&=&\phi_{B}(u,1,\cdots,1) =
\phi_{B}(1,u^{\alpha_{1}},\cdots,u^{\alpha_{1}})\\
&=& \phi_{A}(u^{\alpha_{1}},\cdots,u^{\alpha_{1}})=
\Pi_{k}\left[\frac{1}{1-t^{-1}}\right]^{\esp_{k}}_
{t=u^{\alpha_{k}}}.
\end{eqnarray*}

To finish the proof, it suffices to note that the following 
identification of $H$-vector spaces holds : $V^{+,\beta}\simeq
\oplus_{\esp_{k}> 0}\C_{\alpha_{k}}$ and
$(V\otimes\C)^{+,\beta}\simeq \oplus_{k}\C_{\esp_{k}\alpha_{k}}$. $\Box$

\medskip

\subsection{Localization of the quantization map on \protect $M^{\beta}$}
\label{sec.loc.application.moment}

Let $\beta\neq 0$ be a $G$-invariant element of $\ggot$.
The localization formula that we prove for the 
Riemann-Roch character $RR^{^{G}}(M,-)$ will hold 
in\footnote{An element of $\widehat{R}(G)$ is simply a 
formal sum $\sum_{\lambda}m_{\lambda}\chi_{\lambda}^{_{G}}$ with 
$m_{\lambda} \in \Z$ for all $\lambda$.} 
$\widehat{R}(G):=\hom_{\Z}(R(G),\Z)$.

Let $\Ncal$ be the normal bundle of $M^{\beta}$ in $M$.
For $m\in M^{\beta}$, we have the decomposition $\T_{m}M=
\T_{m}M^{\beta}\oplus\Ncal\vert_{m}$. The linear action of $\beta$
on $T_{m}M$ precises this decomposition. The map $\Lcal^{M}(\beta):
\T_{m}M\to\T_{m}M$ commutes with the map $J$ and
satisfies $\T_{m}M^{\beta}=\ker(\Lcal^{M}(\beta))$. Here we take 
$\Ncal\vert_{m}:={\rm Image}(\Lcal^{M}(\beta))$. Then the almost
complex structure $J$ induces a $G$-invariant almost complex structure 
$J_{\beta}$ on $M^{\beta}$,  and a complex structure $J_{\Ncal}$ on 
the fibers of $ \Ncal\to M^{\beta}$. We have then a quantization map
$RR^{^G}(M^{\beta},-):K_{G}(M^{\beta})\to R(G)$. The torus $\tore_{\beta}$ 
acts linearly on the fibers of the complex vector bundle $\Ncal$. 
Thus we associate the polarized complex $G$-vector bundles
$\Ncal^{+,\beta}$ and $(\Ncal\otimes\C)^{+,\beta}$ (see Definition  
\ref{wedge.V.inverse}).

\begin{theo}\label{th.localisation.pt.fixe}
    For every $E\in K_{G}(M)$, we have the following equality in \break
    $\widehat{R}(G)$ :
    $$
    RR^{^{G}}(M,E)=(-1)^{r_{\Ncal}}\sum_{k\in\N}
    RR^{^{G}}(M^{\beta},E\vert_{M^{\beta}}
    \otimes\det\Ncal^{+,\beta}\otimes S^k((\Ncal\otimes\C)^{+,\beta}) \ ,
    $$
    where $r_{\Ncal}$ is the locally constant function on $M^{\beta}$
     equal to the complex rank of $\Ncal^{+,\beta}$.
\end{theo}    
 
\medskip

Before proving this result let us rewrite this localization formula 
in a more synthetic way. The $G\times\tore_{\beta}$-Riemann-Roch character 
$RR^{^{G \times {\rm T}_{\beta}}}(M^{\beta},-)$  is extended canonically 
to a map from $K_{G}(M^{\beta})\,\widehat{\otimes}\, R(\tore_{\beta})$ to
$R(G)\,\widehat{\otimes}\, R(\tore_{\beta})$ (see Definition 
\ref{wedge.V.inverse}). 
Note that the surjective morphism $G\times \tore_{\beta}\to G, 
(g,t)\mapsto g.t$ induces maps $R(G)\to  R(G)\otimes R(\tore_{\beta})$, 
$K_{G}(M)\to K_{G\times\tore_{\beta}}(M)$, both denoted $k$, with
the tautological relation $k(RR^{^{G}}(M,E))=
RR^{^{G \times {\rm T}_{\beta}}}(M,k(E))$. To simplify, we 
will omit the morphism $k$ in our notations.

Let $\overline{\Ncal}$ be the normal bundle $\Ncal$ with the opposite 
complex structure. With the convention of Definition \ref{wedge.V.inverse}
the element $\wedge_{\C}^{\bullet}\overline{\Ncal}\in 
K_{G\times \tore_{\beta}}(M^{\beta})\simeq
K_{G}(M^{\beta})\otimes R(\tore_{\beta})$ admits a polarized inverse
$\left[\wedge_{\C}^{\bullet}\overline{\Ncal}\,\right]^{-1}_{\beta}\, \in\, 
K_{G}(M^{\beta})\,\widehat{\otimes}\, R(\tore_{\beta})$.
Finally the result of Theorem \ref{th.localisation.pt.fixe} can be 
written as the following equality in $R(G)\,\widehat{\otimes}\, 
R(\tore_{\beta})$ :
\begin{equation}\label{eq.loc.M.beta.simplifie}
RR^{^{G}}(M,E)=RR^{^{G \times {\rm T}_{\beta}}}
\left(M^{\beta},E\vert_{M^{\beta}}\otimes
\left[\wedge_{\C}^{\bullet}\overline{\Ncal}\,\right]^{-1}_{\beta}\right)\ .
\end{equation}

\medskip

Note that Theorem \ref{th.localisation.pt.fixe} gives a proof 
of some rigidity properties 
\cite{Atiyah-Hirzebruch70,Meinrenken-Sjamaar}. Let $H$ be a maximal 
torus of $G$. Following  Meinrenken 
and Sjamaar, a $G$-equivariant complex vector bundle 
$E\to M$ is called {\em rigid} if the action of $H$ on 
$E\vert_{M^H}$ is trivial. Take $\beta\in\hgot$ such that $M^{\beta}=
M^{H}$, and apply Theorem \ref{th.localisation.pt.fixe}, with
$\beta$ and $-\beta$, to $RR^{^H}(M,E)$, with $E$ rigid. 

If we take $+\beta$, Theorem \ref{th.localisation.pt.fixe} shows that 
$h\in H\to RR^{^H}(M,E)(h)$ is of the form $h\in H\to 
\sum_{a\in\hat{H}}n_{a}h^a$ with 
$n_{a}\neq 0\Longrightarrow \langle a,\beta\rangle\geq 0$.
(see Lemma \ref{lem.multiplicites.tore}). If we take $-\beta$, 
we find $RR^{^H}(M,E)(h)= \sum_{a\in\hat{H}}n_{a}h^a$, with  
$n_{a}\neq 0 \Longrightarrow -\langle a,\beta\rangle\geq 0$. Comparing 
the two results, and using the genericity of $\beta$, we see that 
$RR^{^H}(M,E)$ is a {\em constant} function on $H$, hence $RR^{^G}(M,E)$ is 
then a constant function on $G$. We can now
rewrite the equation of Theorem \ref{th.localisation.pt.fixe}, 
where we keep on the right hand side
the {\em constant} terms:
\begin{equation}\label{prop.E.rigid}
RR^{^G}(M,E)=\sum_{F\subset M^{H,+}}RR(F,E\vert_{F})\ .
\end{equation}
Here the summation is taken over all connected components $F$ of $M^{H}$ 
such that $\Ncal_{F}^{+,\beta}=0$ (i.e. we have 
$\langle \xi,\beta\rangle < 0$ for all weights $\xi$ of the $H$-action
on the normal bundle $\Ncal_{F}$ of $F$).

\medskip

{\em Proof of Theorem \ref{th.localisation.pt.fixe} :}

\medskip
Let $\Ucal$ be a $G$-invariant tubular neighborhood\footnote{To simplify the 
notation, we keep the notation $M^{\beta}$ even if we work in fact
on a connected component of the submanifold $M^{\beta}$.} of 
$M^{\beta}$ in $M$. We know from section \ref{sec.general.procedure} that 
$RR^{^G}(M,E)=\indice^{G}_{\Ucal}(\Thom^{\beta}_{G}(M,J)\otimes 
E\vert_{\Ucal})$ where 
$$
\Thom^{\beta}_{G}(M,J)(m,w):=\Thom_{G}(\Vcal,J)
(m,w-\beta_{\Ncal}(m)), \quad {\rm }\quad (m,w)\in\T\Ucal.
$$ 
Let $\phi:\Vcal\to \Ucal$ be $G$-invariant diffeomorphism with a $G$-invariant 
neighbourhood $\Vcal$ of $M^{\beta}$ in  the normal bundle $\Ncal$. 
We denote by $\Thom^{\beta}_{G}(\Vcal,J)$ the symbol 
$\phi^{*}(\Thom^{\beta}_{G}(M,J))$. 
Here we still denote by $J$ the almost complex structure transported on 
$\Vcal$ via the diffeomorphism  $\Ucal\simeq\Vcal$.

Let $p:\Ncal\to M^{\beta}$ be the canonical projection. The choice of 
a $G$-invariant connection on $\Ncal$ induces an isomorphism of $G$-vector bundles 
over $\Ncal$:
\begin{eqnarray}\label{eq:trivialisation.T.N}
\T\Ncal&\tilde{\longrightarrow}& 
p^{*}\left(\T M^{\beta}\oplus \Ncal\right)\\
w&\longmapsto& \T p(w)\oplus (w)^{V}\nonumber
\end{eqnarray}
Here $w\to (w)^{V},\ \T\Ncal\to p^{*}\Ncal$ is the 
projection which associates to a tangent vector its {\em vertical}
part (see \cite{B-G-V}[section 7] or \cite{pep1}[section 4.1]).
The map $\widetilde{J}:=p^{*}(J_{\beta}\oplus J_{\Ncal})$
defines an almost complex structure on the manifold $\Ncal$ which 
is constant over the fibers of $p$. With this 
new almost complex structure $\widetilde{J}$ we construct
the $G$-transversally elliptic symbol over $\Ncal$
$$
\Thom^{\beta}_{G}(\Ncal)(n,w)=\Thom_{G}(\Ncal,\widetilde{J})
(n,w-\beta_{\Ncal}(n)), \quad {\rm }\quad (n,w)\in\T\Ncal.
$$
We denote by $i: \Vcal\to\Ncal$ the inclusion map, and
$i_{*}:K_{G}(\T_{G}\Vcal)\to K_{G}(\T_{G}\Ncal)$ the induced map.

\begin{lem}\label{lem.J.modifie}
    We have
    $$
    i_{*}(\Thom^{\beta}_{G}(\Vcal,J))=
    \Thom^{\beta}_{G}(\Ncal)
    \quad {\rm in} \quad K_{G}(\T\Ncal).
    $$
\end{lem}
 
{\em Proof} : We proceed as in Lemma 
\ref{lem.inv.homotopy}. The complex 
structure $J_{n},\ n\in\Vcal$ and $\widetilde{J}_{n},\ 
n\in\Ncal$ are equal on $M^{\beta}$, and are related by the 
homotopy $J^t_{(x,v)}:=J_{(x,t.v)},\ u\in[0,1]$ for  $n=(x,v)\in \Vcal$. 
Then, as in 
Lemma \ref{lem.inv.homotopy}, we can construct an invertible bundle 
map $A\in\Gamma(\Vcal,\End(\T\Vcal))^{G}$, which is homotopic 
to the identity and such that $A.J= \widetilde{J}.A$ on $\Vcal$.
We conclude as in Lemma \ref{lem.inv.homotopy} that the symbols
$\Thom^{\beta}_{G}(\Vcal,J)$ and 
$\Thom^{\beta}_{G}(\Ncal)\vert\Vcal$ are equal in $K_{G}(\T\Vcal)$. 
Then the Lemma follows from the excision property. $\Box$ 

\medskip

Since $E\simeq p^{*}(E\vert_{M^{\beta}})$, for any $G$-complex vector 
bundle $E$ over $\Ncal$, the former Lemma tells us that 
$RR^{^G}(M,E)=\indice^{G}_{\Ncal}(\Thom^{\beta}_{G}(\Ncal)\otimes 
p^{*}(E\vert_{M^{\beta}}))$.

We consider now the Hermitian vector bundle $\Ncal\to M^{\beta}$ with the 
action of $G\times \tore_{\beta}$. First we use the decomposition
$\Ncal=\oplus_{\alpha}\Ncal^{\alpha}$ relatively to the unitary 
action of $\tore_{\beta}$ on the fibers of $\Ncal$. Let
$N^{\alpha}$ be an Hermitian vector space of dimension equal to the 
rank of $\Ncal^{\alpha}$, equipped with
the representation $t\to t^{\alpha}$ of $\tore_{\beta}$.
Let $U$ be the group of $\tore_{\beta}$-{\em equivariant unitary 
maps} of the vector space $N:= \oplus_{\alpha}N^{\alpha}$, and let $R$ 
be the $\tore_{\beta}$-equivariant unitary frame of 
$(\Ncal,J_{\Ncal})$ framed on $N$. 
Note that $R$ is provided with a 
$U\times G$-action and a trivial action of 
$\tore_{\beta}$ :  for $x\in M^{\beta}$,  any element of $R\vert_{x}$ 
is a $\tore_{\beta}$-equivariant unitary map from $N$
to $\Ncal\vert_{x}$.
The manifold $\Ncal$ is isomorphic to 
$R\times_{U}N$, where $G$ acts on $R$ and 
$\tore_{\beta}$ acts on $N$.

We denote by $\Thom^{\beta}_{G\times\tore_{\beta}}(\Ncal)$ the 
$G\times\tore_{\beta}$ canonical extension of 
$\Thom^{\beta}_{G}(\Ncal)$. It can be considered
as a $G$, $G\times\tore_{\beta}$, or $\tore_{\beta}$-transversally
elliptic symbol. Here we consider $\Thom^{\beta}_{G\times\tore_{\beta}}(\Ncal)$
as an element of $K_{G\times\tore_{\beta}}
(\T_{\tore_{\beta}}(R\times_{U}N))$.
Recall that we have two isomorphisms
\begin{equation}
    \pi_{N}^{*}: 
    K_{G\times\tore_{\beta}}(\T_{\tore_{\beta}}(R\times_{U}N))
    \tilde{\longrightarrow}
K_{G\times\tore_{\beta}\times U}
(\T_{\tore_{\beta}\times U}(R\times N)),
    \label{eq:pi.N.iso}
\end{equation}
\begin{equation}
    \pi^{*}: 
    K_{G}(\T M^{\beta})
    \widetilde{\longrightarrow}
    K_{G\times U}(\T_{U}R),
    \label{eq:pi.iso}
\end{equation}
where $\pi_{N}:R\times N\to R\times_{U}N\simeq \Ncal$ and
$\pi:R\to R/U\simeq M^{\beta}$ are the quotient maps relative 
to the free $U$-action. 
Following (\ref{eq:produit.transversal}), we have a product
\begin{equation}
   K_{G\times U}(\T_{U}R)
   \times
   K_{\tore_{\beta}\times U}(\T_{\tore_{\beta}}N)
   \longrightarrow
   K_{G\times\tore_{\beta}\times U}
   (\T_{\tore_{\beta}\times U}(R\times N)) \ .
    \label{eq:G.T.beta.U}
\end{equation}

The following Thom classes 

\noindent - $\Thom^{\beta}_{G\times\tore_{\beta}}(\Ncal)\in 
K_{G\times\tore_{\beta}}(\T_{\tore_{\beta}}(R\times_{U}N))$,

\noindent - $\Thom^{\beta}_{\tore_{\beta}\times U}(N)\in 
K_{\tore_{\beta}\times U}(\T_{\tore_{\beta}} N)$, and

\noindent - $\Thom_{G}(M^{\beta})\in K_{G}(\T M^{\beta})$

are related by the following equality in $K_{G\times\tore_{\beta}\times U}
(\T_{\tore_{\beta}\times U}(R\times N))$ :
\begin{equation}
    \pi^{*}_{N}\Thom^{\beta}_{G\times\tore_{\beta}}(\Ncal)=
    (\pi^{*}\Thom_{G}(M^{\beta}))\odot
    \Thom^{\beta}_{\tore_{\beta}\times U}(N).
    \label{eq:Thom.egalite}
\end{equation}
We will justify (\ref{eq:Thom.egalite}) later. Every $E\in 
K_{G}(M)$, when restrict to $M^{\beta}$, admit the decomposition 
$E\vert_{M^{\beta}}=
\sum_{a\in\widehat{\tore_{\beta}}} E^a\otimes \C_{a}$ in 
$K_{G\times\tore_{\beta}}(M^{\beta})\simeq 
K_{G}(M^{\beta})\otimes R(\tore_{\beta})$. Multiplication of 
(\ref{eq:Thom.egalite}) by $E$ gives 
$$
 \pi^{*}_{N}(\Thom^{\beta}_{G\times\tore_{\beta}}(\Ncal)
 \otimes E\vert_{M^{\beta}}) =
 \sum_{a\in\widehat{\tore_{\beta}}}
    \pi^{*}(\Thom_{G}(M^{\beta})\otimes E^a)\odot
    (\Thom^{\beta}_{\tore_{\beta}\times U}(N)\otimes \C_{a}).
$$
    
Following (\ref{eq:formule.G.H.produit}) and Theorem (\ref{thm.atiyah.1}), 
the last equality gives, after taking the index and the $U$-invariant :
\begin{eqnarray}\label{eq:indice.egalite}
 \lefteqn{RR^{^{G\times{\rm T}_{\beta}}}(M,E)=}\nonumber\\
& & \sum_{a}\left[\sum_{i\in\widehat{U}} RR^{^G}(M^{\beta},E^a\otimes 
\underline{W}_{i}^{*})\cdot W_{i}\cdot\indice^{\tore_{\beta}\times U}
\left(\Thom^{\beta}_{\tore_{\beta}\times U}(N)\right)\cdot 
    \C_{a}\right]^U\ .
\end{eqnarray}

Here we used that $RR^{^{G\times{\rm T}_{\beta}}}(M,E)$ 
is equal to the $U$-invariant part of \break 
$\indice^{G\times\tore_{\beta}\times U}(\pi^{*}_{N}
(\Thom^{\beta}_{G\times\tore_{\beta}}(\Ncal)\otimes E\vert_{M^{\beta}}))$, 
and the index of 
$\pi^{*}(\Thom_{G}(M^{\beta})\otimes E^a)$ is
equal to $\sum_{i\in\widehat{U}}RR^{^G}(M^{\beta},E^a\otimes 
\underline{W}_{i}^{*}).W_{i}$.

Now we observe that for any $L\in R(U)$, the $U$-invariant part of 
\break 
$\sum_{i\in\widehat{U}}RR^{^G}(M^{\beta},E\vert_{M^{\beta}}\otimes 
\underline{W}_{i}^{*}).W_{i}\otimes L$
is equal to $RR^{^G}(M^{\beta},E\vert_{M^{\beta}}\otimes \underline{L})$ 
with $\underline{L}=R\times_{U}L$. With the computation of 
$\indice^{\tore_{\beta}\times U}(
\Thom^{\beta}_{\tore_{\beta}\times U}(N))$ given in
Proposition \ref{prop.indice.thom.beta.V} we obtain finally
$$
RR^{^{G\times{\rm T}_{\beta}}}(M,E)=
(-1)^{r_{\Ncal}}\sum_{k\in\N}
RR^{^{G\times{\rm T}_{\beta}}}
\Big(M^{\beta},E\vert_{M^{\beta}}\otimes\det \Ncal^{+,\beta}
\otimes S^k((\Ncal\otimes\C)^{+,\beta})\Big) \ 
$$
which implies the equality of Theorem \ref{th.localisation.pt.fixe}. 

\medskip

We give now an explanation for (\ref{eq:Thom.egalite}), 
which is a direct consequence of the fact that
the almost complex  structure $\widetilde{J}$ admits the decomposition
$\widetilde{J}=p^{*}(J_{\beta}\oplus J_{\Ncal})$.
Hence $\wedge_{\C}^{\bullet}\T_{n}\Ncal$ equipped with the map 
$Cl_{n}(v-\beta_{\Ncal}(n)),\ v\in \T_{n}\Ncal$ is isomorphic
to $\wedge_{\C}^{\bullet}\T_{x}M^{\beta}\otimes
\wedge_{\C}^{\bullet}\Ncal\vert_{x}$
equipped with $Cl_{x}(v_{1})\odot Cl_{x}(v_{2}-\beta_{\Ncal}(n))$
where $x=p_{a}(n)$, and the vector $v\in\T_{n}\Ncal$ 
is decomposed, following the isomorphism (\ref{eq:trivialisation.T.N}),
in $v=v_{1}+v_{2}$ with $v_{1}\in \T_{x}M^{\beta}$ and $v_{2}\in \Ncal\vert_{x}$.
Note that the vector $w=\beta_{\Ncal}(n)\in \T_{n}\Ncal$ is 
vertical, i.e. $w=(w)^{V}$. $\Box$

\medskip 


\section{Localization via an abstract moment map}\label{sec.Localisation.f}

\medskip

Let $(M,J,G)$ be a compact $G$-manifold provided with a $G$-invariant 
almost complex structure. We denote by $RR^{^{G}} : K_{G}(M)\to R(G)$ the 
quantization map. Here we suppose that the $G$-manifold is equipped
with an {\em abstract moment map}  \cite{G-G-K3,Karshon.98}. 

\begin{defi}\label{moment-map}
 A smooth map $f_{_{G}}:M\to\ggot^{*}$ is called an 
 {\rm abstract moment map} if
 
 i) the map $f_{_{G}}$ is equivariant for the action of the group $G$, and
 
 ii)\footnote{Condition ii) is equivalent to the following :  
for every $X\in\ggot$, the fonction $\langle f_{_{G}},X \rangle$ 
is locally constant on $M^X$.} for every Lie subgroup $K\subset G$ with 
Lie algebra $\kgot$, the induced map $f_{_{K}}:M\to \kgot^{*}$ is 
locally constant on the submanifold 
$M^{K}$ of fixed points for the $K$-action (the map  
$f_{_{K}}$ is the composition of $f_{_{G}}$ with the projection 
$\ggot^{*}\to\kgot^{*}$).
\end{defi}

The terminology ``moment map'' is usually used when we work in the 
case of a Hamiltonian action. More precisely, when the manifold is 
equipped with a symplectic $2$-form $\omega$ which is $G$-invariant, 
a {\em moment map} $\Phi:M\to\ggot^{*}$ relative to 
$\omega$ is a $G$-equivariant map satisfying 
$d\langle\Phi,X\rangle=-\omega(X_{M},-),\ X\in\ggot$.

\medskip
 
For the rest of this paper we make the choice of a $G$-invariant 
scalar product over $\ggot^{*}$. This defines an identification 
$\ggot^{*}\simeq\ggot$, and we work with a given abstract moment map
$f_{_{G}}: M\to \ggot$.

\begin{defi}\label{def.H.vector}
Let $\Hcal^{^{G}}$ be the $G$-invariant vector field over $M$ defined by
$$
\Hcal^{^{G}} _{m}:=(f_{_{G}}(m)_{M})_{m},\quad \forall \ m\in M.
$$
\end{defi}

The aim of this section is to compute the localization, as in section
\ref{sec.general.procedure}, with the $G$-invariant vector field 
$\Hcal^{^{G}}$. We know that the Riemann-Roch character is localized near 
the set $\{\Phi_{\Hcal^{^{G}}}=0\}$, but we see that
$\{\Phi_{\Hcal^{^{G}}}=0\}=\{\Hcal^{^{G}}=0\}$.
We will denote by $C^{f_{_{G}}}$ this set.
Let $H$ be a maximal torus of $G$, with Lie algebra $\hgot$, and let 
$\hgot_{+}$ be a Weyl chamber in $\hgot$.

\begin{lem}\label{lem.C.f.G}
There exists a finite subset $\Bcal_{_{G}}\subset\hgot_{+}$, such that
$$
C^{f_{_{G}}}=\bigcup_{\beta\in\Bcal_{_{G}}}C^{^{G}}_{\beta},\quad 
{\rm with}\quad 
C^{^{G}}_{\beta}=G.(M^{\beta}\cap f_{_{G}}^{-1}(\beta)).
$$    
\end{lem}

{\em Proof} : We first observe that $\Hcal^{^{G}} _{m}=0$ if and only
if $f_{_{G}}(m)=\beta'$ and $\beta'_{M}\vert_{m}=0$, that is
$m\in M^{\beta'}\cap f_{_{G}}^{-1}(\beta')$, for some 
$\beta'\in\ggot$. For every $\beta'\in\ggot$,  there exists 
$\beta\in\hgot_{+}$, with $\beta'=g.\beta$ for some $g\in G$. Hence
$M^{\beta'}\cap f_{_{G}}^{-1}(\beta')=g.(M^{\beta}\cap 
f_{_{G}}^{-1}(\beta))$. We have shown that 
$C^{f_{_{G}}}=\bigcup_{\beta\in\hgot_{+}}C^{^{G}}_{\beta}$, and we 
need to prove that the set $\Bcal_{_{G}}:=\{\beta\in\hgot_{+},\
M^{\beta}\cap f_{_{G}}^{-1}(\beta)\neq\emptyset\}$ is finite. Consider the set 
$\{H_{1},\cdots,H_{l}\}$ of stabilizers for the action of the torus $H$ 
on the compact manifold $M$. For each $\beta\in\hgot$ we denote by
$\tore_{\beta}$ the subtorus of $H$
generated by $\exp(t.\beta),\ t\in\R$, and we observe that
\begin{eqnarray*}
    M^{\beta}\cap f_{_{G}}^{-1}(\beta)\neq\emptyset
&\Longleftrightarrow&
\exists H_{i}\ {\rm such\ that}\ 
\tore_{\beta}\subset H_{i}\ \ {\rm and}\ \ M^{H_{i}}\cap 
f_{_{G}}^{-1}(\beta)\neq\emptyset \\
&\Longleftrightarrow&
\exists H_{i}\ {\rm such\ that}\ \beta\in f_{_{G}}(M^{H_{i}})\cap 
Lie(H_{i}).\\
\end{eqnarray*}
But $f_{_{G}}(M^{H_{i}})\cap Lie(H_{i})\subset 
f_{_{H_{i}}}(M^{H_{i}})$ is a finite set after Definition \ref{moment-map}.
The proof is now completed. $\Box$

\begin{defi}\label{def.thom.beta.f}
Let $\Thom_{G,[\beta]}^f(M)\in K_{G}(\T_{G}\Ucal^{^{G,\beta}})$ defined
by 
$$
\Thom_{G,[\beta]}^f(M)(x,v):=\Thom_{G}(M)(x,v-\Hcal^{^{G}}_{x}), \quad
{\rm for}\quad (x,v)\in \T\Ucal^{^{G,\beta}}\ .
$$
Here $i^{^{G,\beta}}:\Ucal^{^{G,\beta}}\croc M$ is any $G$-invariant 
neighbourhood of $C^{^{G}}_{\beta}$ such that 
$\overline{\Ucal^{^{G,\beta}}}\cap C^{f_{_{G}}}= C^{^{G}}_{\beta}$.
\end{defi}

\begin{defi}\label{def.RR.beta}
For every $\beta\in \Bcal_{_{G}}$, we denote by $RR_{\beta}^{^{G}}(M,-):
K_{G}(M)\to R^{-\infty}(G)$ the localized Riemann-Roch character 
near $C^{^{G}}_{\beta}$, defined as in (\ref{eq:RR.localise}), by 
$$
RR_{\beta}^{^{G}}(M,E)=\indice^{G}_{\Ucal^{^{G,\beta}}}
\left(\Thom_{G,[\beta]}^f(M)\otimes E_{\vert\Ucal^{^{G,\beta}}}\right)\ ,
$$ 
for $E\in K_{G}(M)$. Note that the map $RR_{\beta}^{^{G}}(M,-)$ is 
well defined on a {\em non-compact} manifold $M$ when the abstract 
moment map is proper, since we can take 
$\Ucal^{^{G,\beta}}$ relatively compact and the index map 
$\indice^{G}_{\Ucal^{^{G,\beta}}}$ is then defined (see Corollary 
\ref{hyp.indice}).
\end{defi}

According to Proposition \ref{prop.localisation}, we have the partition 
$RR^{^{G}}(M,-)=$ \break 
$\sum_{\beta\in\Bcal_{_{G}}}RR_{\beta}^{^{G}}(M,-)$, and 
the rest of this article is devoted to the analysis of the maps 
$RR_{\beta}^{^{G}}(M,-),\ \beta\in \Bcal_{_{G}}$. 

In subsections \ref{subsec.RR.G.beta} and \ref{subsec.induction.G.H} 
we prove that $[RR_{\beta}^{^{G}}(M,E)]^G=0$, when $E$ is 
$f_{_{G}}$-{\em strictly positive} with 
$\eta_{_{E,\beta}}>\langle \theta,\beta\rangle$ (see Def. 
\ref{eq.mu.positif} for the notion of 
$f_{_{G}}$-{\em positivity}). The next two subsections are 
devoted to the computation of
$RR^{^{G}}_{0}(M,-)$ when $0$ is a regular value of 
the abstract moment map $f_{_{G}}$.


\subsection{Induced $\spinc$ structures} \label{subsec.spinc}

In this subsection we first review the notion of $\spinc$-structures
(see \cite{Lawson-Michel,Duistermaat96,Schroder}). After we show that the
almost complex structure $J$ on $M$ induces a $\spinc$-structure on 
$\Mcal_{red}$.

The group $\spin_n$ is the connected double cover of the group 
$\so_n$. Let $\eta :\spin_n\to\so_n$ be the covering map, and let 
$\esp$ be the element who generates the kernel. The group $\spinc_n$ 
is the quotient $\spin_n\times_{\Z_2}\u_1$, where $\Z_2$ acts by 
$(\esp,-1)$. There are two canonical group homomorphisms
$$
\eta:\spinc_n\to\so_n\quad ,\quad \Det :\spinc_n\to \u_1\ 
$$
such that  $\eta^{\rm c}=(\eta,\Det):\spinc_n\to\so_n\times \u_1$ is a 
double covering map.

Let $p:E\to M$ be an oriented Euclidean vector bundle of rank $n$, and 
let $\Pso(E)$ be its bundle of oriented orthonormal frames. A 
$\spinc$-structure on $E$ is a $\spinc_n$-principal bundle 
$\Pspin(E)\to M$, together with a $\spinc$-equivariant map 
$\Pspin(E)\to\Pso(E)$. The line bundle 
$\Lfibre:=\Pspin(E)\times_{\Det}\C$ is called the determinant line bundle 
associated to $\Pspin(E)$.  Whe have then a double covering 
map\footnote{If $P$, $Q$ are  principal bundle over $M$ respectively 
for the groups $G$ and $H$, we denote simply by $P\times Q$ their 
fibering product over $M$ which is a $G\times H$ principal bundle 
over $M$.}
\begin{equation}\label{eq.spin.covering}
\eta^{\rm c}_E \, : \, \Pspin(E)\longrightarrow\Pso(E)\times\Pu(\Lfibre) \ ,
\end{equation}
where $\Pu(\Lfibre):=\Pspin(E)\times_{\Det}\u_1$ is the associated 
$\u_1$-principal bundle over $M$.

A $\spinc$-structure on an oriented Riemannian manifold is a 
$\spinc$-structure on its tangent bundle. If a group $K$ acts on the 
bundle $E$, preserving the orientation and the Euclidean structure, 
we defines a $K$-equivariant $\spinc$-structure by requiring $\Pspin(E)$ 
to be a $K$-equivariant principal bundle, and (\ref{eq.spin.covering}) 
to be $(K\times\spinc_n)$-equivariant.

\medskip

We assume now that $E$ is of even rank $n=2m$. 
Let $\Delta_{2m}$ be the irreducible complex Spin representation of 
$\spinc_{2m}$. Recall that $\Delta_{2m}=\Delta_{2m}^+\oplus\Delta_{2m}^-$ 
inherits a canonical Clifford action $\clif :\R^{2m}\to\End_{\C}(\Delta_{2m})$ 
which is $\spinc_{2m}$-equivariant, and which interchanges the 
graduation : $\clif(v):\Delta_{2m}^{\pm}\to\Delta_{2m}^{\mp}$, for 
every $v\in\R^{2m}$. Let 
\begin{equation}\label{eq.spinor.bundle}
  \Scal(E):=\Pspin(E)\times_{\spinc_{2m}}\Delta_{2m}
 \end{equation} 
be the irreducible complex spinor bundle over $E\to M$. The 
orientation on the fibers of $E$ defines a graduation  
$\Scal(E):=\Scal(E)^+\oplus\Scal(E)^-$. Let $\overline{E}$ be the bundle 
$E$ with opposite orientation. A $\spinc$ structure on $E$ induces a 
$\spinc$ on $\overline{E}$, with the same determinant line bundle, and such that 
$\Scal(\overline{E})^{\pm}=\Scal(E)^{\mp}$. 

More generaly, we associated 
to an Euclidean vector bundle $p: E\to M$ its Clifford bundle 
$\Clif(E)\to M$. A complex vector bundle $\Scal\to M$ is 
called a complex spinor bundle over $E\to M$ if it is a 
left-$\Clif(E)$-module; moreover $\Scal$ is called irreducible if
$\Clif(E)\otimes\C\simeq\End_{\C}(\Scal)$. In fact the notion 
of $\spinc$-structure (in terms of principal bundle) on a Euclidean 
bundle $E\to M$ is equivalent to the existence of an irreducible 
complex spinor bundle over $E\to M$  \cite{Schroder}.

Since $E=\Pspin(E)\times_{\spinc_{2m}}\R^{2m}$, the bundle $p^*\Scal(E)$
is isomorphic to \break 
$\Pspin(E)\times_{\spinc_{2m}}(\R^{2m}\oplus\Delta_{2m})$.
\begin{defi}
  Let $\sthom(E): p^*\Scal(E)^+\to p^*\Scal(E)^-$ be the  symbol 
  defined by 
  \begin{eqnarray*}
    \Pspin(E)\times_{\spinc_{2m}}(\R^{2m}\oplus\Delta_{2m}^+)
    &\longrightarrow&
   \Pspin(E)\times_{\spinc_{2m}}(\R^{2m}\oplus\Delta_{2m}^-)\\
  {} [p;v,w]&\longmapsto &[p,v,\clif(v)w]\ .
  \end{eqnarray*} 

When $E$ is the tangent bundle of a manifold $M$, the symbol 
$\sthom(E)$ is denoted by $\sthom(M)$. If a group $K$ acts equivariantly on 
the $\spinc$-stucture, we denote by $\sthom_{K}(E)$ the equivariant 
symbol.
\end{defi}

The characteristic set of $\sthom(E)$ is $M\simeq\{{\rm zero\ section\ 
of}\ E\}$, hence it defines a class in $K(E)$ if $M$ is compact. When
$E=\T M$, the symbol $\sthom(M)$ corresponds to the {\em principal symbol} 
of the $\spinc$ Dirac operator associated to the $\spinc$-structure 
\cite{Duistermaat96}. When $M$ is compact, we define a quantization map 
$\Qcal(M,-):K(M)\to\Z$ by the relation
$\Qcal(M,E):=\indice_{M}(\sthom(M)\otimes E)$ : $\Qcal(M,E)$ is the
index of the $\spinc$ Dirac operator on $M$ twisted by $E$.

These notions extend to the orbifold case. Let $M$ be a manifold with 
a locally free action of a compact Lie group $G$. The quotient 
$\Xcal:=M/G$ is an orbifold, a space with finite quotient 
singularities. A $\spinc$ structure on $\Xcal$ is by definition a
$G$-equivariant $\spinc$ structure on the bundle $\T_{G}M\to M$; where 
$\T_{G}M$ is identified with the pullback of $\T\Xcal$ via the 
quotient map $\pi:M\to\Xcal$. We define in the same way 
$\sthom(\Xcal)\in K_{orb}(\T\Xcal)$, such that $\pi^{*}\sthom(\Xcal)=
\sthom_{G}(\T_{G}M)$. The pullback by $\pi$ induces an isomophism $\pi^{*}: 
K_{orb}(\T\Xcal)\simeq K_{G}(\T_{G}M)$. The quantization map 
$\Qcal(\Xcal,-)$ is defined by : $\Qcal(\Xcal,\Ecal)=
\indice_{\Xcal}(\sthom(\Xcal)\otimes\Ecal)$.

\begin{lem}\label{lem.q.change}
  Let $E\to M$ be an oriented $G$-bundle . Let $g_0,g_1$ be two 
  $G$-invariant metric on the fibers of $E$, and suppose that
  $(E,g_0)$ admits an equivariant $\spinc$-stucture denoted  by 
  $\Pspin(E,g_0)$. The trivial homotopy
  $g_t=(1-t).g_0 +t.g_1$ between the metrics, induces an equivariant 
  homotopy between the principal bundles $\Pso(E,g_0)$, $\Pso(E,g_1)$ which
  can be lift to an equivariant homotopy between $\Pspin(E,g_0)$ and a 
  $\spinc$-bundle over $(E,g_1)$. When the base $M$ is compact, the
  corresponding symbols $\sthom_G(E,g_0)$ and $\sthom_G(E,g_1)$ define 
  the same class in $K_G(E)$. 
\end{lem}

{\em Proof :} Let $\Scal$ be the irreducible complex spinor bundle 
associated to $\Pspin(E,g_0)$. We denote by 
$\clif_0:\Clif(E,g_0)\to\End_{\C}(\Scal)$ the 
corresponding Clifford action. Let $A_t$ be the unique 
$g_0$-symmetric endomorphism of $E$ such that 
$g_t(v,w)=g_0(A_t(v),A_t(w))$. The composition $\clif_0\circ A_t$ is 
then a Clifford action of $(E,g_t)$ on $\Scal$. It defines a 
$\spinc$-structure on the bundle $(E,g_t)$ which is homotopic to 
$\Pspin(E,g_0)$. $\Box$

\medskip

Consider now the case of a {\em complex} vector bundle $E\to M$, of 
complex rank $m$. The orientation on the fibers of $E$ is given by the 
complex structure $J$. 
Let $\Pu(E)$ be the bundle of unitary frames on $E$. We have a 
morphism ${\rm j}:\u_m\to\spinc_{2m}$ which makes the   
diagram\footnote{Here ${\rm i}:\u_m\croc\so_{2m}$ is the 
canonical inclusion map.}  
\begin{equation}
\xymatrix@C=2cm{
 \u_m\ar[r]^{\rm j} \ar[dr]_{{\rm i}\times\det} & 
 \spinc_{2m}\ar[d]^{\eta^{\rm c}}\\
     & \so_{2m}\times\u_1\ .}
\end{equation}
commutative \cite{Lawson-Michel}. Then 
\begin{equation}\label{eq.J.spin}
  \Pspin(E):=\spinc_{2m}\times_{\rm j}\Pu(E)
\end{equation}
defines a $\spinc$-structure over $E$, with bundle of irreducible 
spinors $\Scal(E)=\wedge^{\bullet}_{\C}E$ and determinant line bundle 
equal to $\det_{\C}E$. 

\begin{rem} Let $M$ be a manifold equipped with an almost complex 
structure $J$. The symbol $\sthom(M)$  defined by the $\spinc$-structure 
(\ref{eq.J.spin}), and the Thom symbol $\Thom(M,J)$  defined in section 
\ref{sec.quantization}  coincide.
\end{rem}

\bigskip

Consider our case of interest, where $M$ is a compact $G$-manifold 
equipped with an equivariant almost complex structure $J$ and with an
abstract moment map $f_{_G}:M\to\ggot^*$. Here we assume that $0$ is a
regular value of $f_{_G}$ : $\Zcal:=f_{_G}^{-1}(0)$ is a smooth
submanifold of $M$ with a locally free action of $G$. Let 
$\Mcal_{red}:= \Zcal/G$ be the corresponding `reduced' space, and 
let $\pi:\Zcal\to \Mcal_{red}$ be the projection map. 
On $\Zcal$ we have  an exact sequence
$0 \longrightarrow \T\Zcal \longrightarrow 
\T M\vert_{\Zcal} \stackrel{df_G}{\longrightarrow}
\ggot^*\times\Zcal \longrightarrow 0$,  
and $\T\Zcal=\T_G\Zcal\oplus \ggot_{\Zcal}$ where 
$\ggot_{\Zcal}\simeq\ggot\times\Zcal$ 
denotes the trivial bundle corresponding to the subspace of $\T\Zcal$
formed by the vector field generated by the infinitesimal action of
$\ggot$. So $\T M\vert_{\Zcal}$ admits the decomposition
\begin{equation}\label{eq.TM.decompose}
\T M\vert_{\Zcal}=\T_G\Zcal\oplus \ggot_{\Zcal}\oplus
\ggot^*\times\Zcal \ .
\end{equation}
The bundle $\pi^{*}(\T\Mcal_{red})$ is identified with $\T_G\Zcal$. 
Thus the decomposition (\ref{eq.TM.decompose}) can be rewritten
\begin{equation}\label{eq.TM.decompose.bis}
\T M\vert_{\Zcal}=\pi^{*}(\T\Mcal_{red})\oplus \ggot_{\C}\times\Zcal \ .
\end{equation} 
with the convention $\ggot_{\Zcal}=(\ggot\otimes i\R)\times\Zcal$ and
$\ggot^*\times\Zcal=(\ggot\otimes \R)\times\Zcal$.

\begin{lem} \label{lem.spinc.induit}
  The data $(J,f_{_G})$ induce :
  \begin{itemize}
    \item an orientation $o_{red}$ on $\Mcal_{red}$,
    
    \item a $\spinc$-structure ${\rm Q}_{red}$ on $(\Mcal_{red},o_{red})$.

  \end{itemize}
Moreover, the irreducible complex spinor bundle 
$\wedge^{\bullet}_{J}\T M$, when 
restricted to $\Zcal$, defines a complex spinor bundle over 
$\pi^{*}(\T\Mcal_{red})\oplus \ggot_{\C}\times\Zcal$ which is 
homotopic to $\pi^{*}\Scal(\Mcal_{red})\otimes 
\wedge_{\C}^{\bullet}\ggot_{\C}\times\Zcal$. 
\end{lem}  

{\em Proof} : Since $\ggot_{\C}\times\Zcal$ is canonically oriented by 
the complex multiplication by $i$, the orientation $o(J)$ on $M$  
determines an orientation $o(\Mcal_{red})$ on $\T \Mcal_{red}$ such 
that $o(J)= o(\Mcal_{red})\,o(\imath)$.

Let $g_0$ be the Riemannian metric on $\T M\vert_{\Zcal}$ equal to the 
restriction to $\Zcal$ of the Riemannian metric on $M$ (which is taken 
compatible with $J$).  If ${\rm P}$ is the $\spinc$-structure on $M$ 
determined by $J$ (see \ref{eq.J.spin}), the restriction ${\rm 
P}\vert_{\Zcal}$ is then a $\spinc$-structure on $(\T 
M\vert_{\Zcal},g_0)$.Let $g_1$ be a $G$-invariant metric on the bundle 
$\T M\vert_{\Zcal}$ which makes (\ref{eq.TM.decompose.bis}) an 
orthogonal sum, and which is constant on the the trivial bundle 
$\ggot_{\C}\times\Zcal$.  We know from Lemma \ref{lem.q.change} that 
the $\spinc$-structure ${\rm P}\vert_{\Zcal}$ on $(\T 
M\vert_{\Zcal},g_0)$ is homotopic to $\spinc$-structure ${\rm P}_1$ on 
$(\T M\vert_{\Zcal},g_1)$ (both are $G$-equivariant).
 
The $\so_{2k}\times\u_{l}$-principal bundle 
$\Pso(\pi^{*}(\T \Mcal_{red}))\times\Pu(\ggot_{\C}\times\Zcal)$ 
is a reduction\footnote{Here $2n=\dim M$, $2k=\dim \Mcal_{red}$ and 
$l=\dim(\ggot)$, so $n=k+l$.} of the $\so_{2n}$ principal bundle 
$\Pso(\pi^{*}(\T\Mcal_{red})\oplus\ggot_{\C}\times\Zcal)$, thus  
we have the commutative diagram

\begin{equation}\label{diagram.Q.seconde}
\xymatrix@C=2cm{
 {\rm Q}\ar[r]\ar[d] & 
 \Pso(\pi^{*}(\T\Mcal_{red}))\times\Pu(\ggot_{\C}\times\Zcal)
 \times\Pu(\Lfibre\vert_{\Zcal})\ar[d]\\
  {\rm P}_1\ar[r]   & 
  \Pso(\pi^{*}(\T\Mcal_{red})\oplus\ggot_{\C}\times\Zcal)
  \times\Pu(\Lfibre\vert_{\Zcal})\ ,
  }
\end{equation}
where $\Lfibre=\det_{\C}(\T M,J)$. Here ${\rm Q}$ is a 
$(\eta^{\rm c})^{-1}(\so_{2k}\times\u_{l})
\simeq\spinc_{2k}\times\u_{l}$-principal bundle. 
Finally we see that 
${\rm Q}_{red}={\rm Q}/(\u_{l}\times G)$ is a $\spinc$ structure 
on $\Mcal_{red}$ with determinant line bundle 
$\Lfibre_{red}=\det_{\C}(\T M\vert_{\Zcal})/G$. 

The irreducible complex spinor bundle $\wedge^{\bullet}_{J}\T M$, 
when restricted to $\Zcal$, is homotopic to $\Scal'={\rm 
P}_1\times_{\spinc_{2n}}\Delta_{2m}$. Using (\ref{diagram.Q.seconde}) 
we get 
\begin{eqnarray*}
  \Scal' &=& {\rm Q}\times_{(\spinc_{2k}\times\u_l)}
  \Big(\Delta_{2k}\otimes\wedge^{\bullet}\C^{l}\Big)\\
  &=& \Big(({\rm Q}/\u_l)\times_{\spinc_{2k}}\Delta_{2k}\Big)\otimes
      \Big(({\rm Q}/\spinc_{2k})\times_{\u_l}\wedge^{\bullet}\C^{l}\Big)\\
  &=& \pi^{*}\Scal(\Mcal_{red})\otimes 
  (\wedge^{\bullet}\ggot_{\C})\times\Zcal \ .
\end{eqnarray*}  
Here we have used the identifications 
${\rm Q}/\spinc_{2k}=\Pu(\ggot_{\C}\times\Zcal)$ and 
$\Pu(\ggot_{\C}\times\Zcal)\times_{\u_l}\wedge^{\bullet}\C^{l}=
(\wedge^{\bullet}\ggot_{\C})\times\Zcal$. $\Box$

\bigskip

We shall consider the particular case where $J$ defines an almost 
complex structure on $\Mcal_{red}$. It happens when the 
following decomposition holds
\begin{equation}\label{eq.J.induit}
\T M\vert_{\Zcal}=\T\Zcal\oplus J(\ggot_{\Zcal})\ .
\end{equation}
With (\ref{eq.J.induit}), $\T M\vert_{\Zcal}$ decomposes in 
$\T M\vert_{\Zcal}=\pi^{*}(\T\Mcal_{red})\oplus
\ggot_{\Zcal}\oplus J(\ggot_{\Zcal})$ : let us denote by 
$pr:\T M\vert_{\Zcal}\to \pi^{*}(\T\Mcal_{red})$ the corresponding 
projection. Since 
$\ggot_{\Zcal}\oplus J(\ggot_{\Zcal})$ is invariant by $J$, the 
endomorphism $J_{red}:=pr\circ J$ is a $G$-invariant almost complex 
structure on $\pi^{*}(\T\Mcal_{red})$.

Using the identification $\ggot\simeq\ggot^*$, one considers the 
endomorphism $\Dcal$ of the trivial bundle $\ggot\times\Zcal$ defined 
by
\begin{equation}\label{eq.Dcal}
\Dcal(X)=-df_{_G}(J(X_{\Zcal}))\ ,\quad {\rm for}\quad X\,\in\,\ggot.
\end{equation}
Condition (\ref{eq.J.induit}) is then equivalent to : 
$\det\Dcal(z)\neq 0$ for all $z\in\Zcal$.  We shall use the normalized 
map $\Dcal(\Dcal^t\Dcal)^{-1/2}$ which is an orthogonal map for the 
fixed Euclidean structure on $\ggot$ (to simplify we keep the same 
notation $\Dcal$ for it). Let $J_{\Dcal}$ be the complex structure 
on the trivial bundle $\ggot_{\C}\times\Zcal$ defined by the 
following matrix
$$
J_{\Dcal}:=\left(
\begin{array}{cc}
0   & -\Dcal  \\ \Dcal^{-1} & 0
\end{array}
\right)\ .
$$

\begin{lem} \label{lem.J.induit} 
Suppose that the decomposition 
(\ref{eq.J.induit}) holds.  On\footnote{Here we use the decompostion 
(\ref{eq.TM.decompose.bis}) of $\T M\vert_{\Zcal}$.} $\T 
M\vert_{\Zcal}=$ \break 
$\pi^{*}(\T\Mcal_{red})\oplus\ggot_{\C}\times\Zcal$ the almost complex 
structure $J$ is homotopic to $J_{red}\times J_{\Dcal}$.  Hence the 
irreducible complex spinor bundle $\wedge^{\bullet}_{J}\T M$, 
when restricted to $\Zcal$, defines a complex spinor bundle over 
$\pi^{*}(\T\Mcal_{red})\oplus \ggot_{\C}\times\Zcal$ which is 
homotopic to \break $\pi^{*}(\wedge^{\bullet}_{J_{red}}\T\Mcal_{red})
\otimes \wedge^{\bullet}_{J_{\Dcal}}\ggot_{\C}\times\Zcal$.  
\end{lem}

{\em Proof :} Trough the decomposition 
$\T M\vert_{\Zcal}=\pi^{*}(\T\Mcal_{red})\oplus \ggot_{\Zcal}\oplus 
J(\ggot_{\Zcal})$, the map $J$ is described by the matrix 
$$
\left(
\begin{array}{cc}
J_{red}  & 0  \\ A & \imath
\end{array}
\right)\ ,
$$
hence $J$ is homotopic to 
$$
J'=\left(
\begin{array}{cc}
J_{red}  & 0  \\ 0 & \imath
\end{array}
\right)\ .
$$
In the decomposition (\ref{eq.TM.decompose.bis}), $J'$ has the following 
matrix 
$$
\left(
\begin{array}{cc}
J_{red}  & B  \\ 0 & C
\end{array}
\right)\ ,
$$
with $C\in\End(\ggot_{\C}\times\Zcal)$ of the form
$$
\left(
\begin{array}{cc}
-\Dcal b\Dcal^{-1}& -\Dcal  \\ b^2\Dcal^{-1} + \Dcal^{-1} & b
\end{array}
\right)\ .
$$
Hence $J'$ is tied to $J_{red}\times J_{\Dcal}$ through 
the homotopies $t\to t\, B$ and $t\to t\, b$, $0\leq t\leq 1$. 
$\Box$


\subsection{The map \protect $RR^{^{G}}_{0}$}\label{subsec.RR.G.O}

The map $RR^{^{G}}_{0}(M,-): K_{G}(M)\to R^{-\infty}(G)$ is 
the Riemann-Roch character localized near $C^{^G}_{0}=
f_{_{G}}^{-1}(0)$ (see Definition \ref{def.RR.beta}). In particular, 
$RR^{^{G}}_{0}(M,-)$ is the zero map if $0$ does not belong to $f_{_{G}}(M)$.
In this subsection, we assume that $0\in f_{_{G}}(M)$ is a regular value of 
$f_{_{G}}$. We have proved in the past subsection that $J$ induces an
orientation $o(\Mcal_{red})$ on the reduced space $\Mcal_{red}$
together with a $\spinc$-structure on $(\Mcal_{red},o(\Mcal_{red}))$. 
Let $\sthom(\Mcal_{red})$ be the elliptic symbol defined by this
$\spinc$-structure and let $\Qcal(\Mcal_{red},-)$ be the corresponding
quantization map.

\begin{prop}\label{prop.RR.T.0}
    For every $G$-equivariant vector bundle $E\to M$, we have  
    \begin{equation}\label{eq.RR.G.0}
    RR_{0}^{^{G}}(M,E)=\sum_{\mu\in\Lambda^{*}_{+}}
    \Qcal(\Mcal_{red},E_{red}\otimes\underline{V_{\mu}}^{*}) . V_{\mu} 
   \quad {\rm in}\quad R^{-\infty}(G) \ ,
    \end{equation}
    Here $E_{red}=E/G$ is the orbifold vector bundle on 
    $\Mcal_{red}$ induced by $E$, and $\underline{V_{\mu}}=
    \Zcal\times_{G}V_{\mu}$. In particular, the 
    $G$-invariant part of $RR_{0}^{^{G}}(M,E)$ is equal to \break 
    $\Qcal(\Mcal_{red},E_{red})\in \Z$.
\end{prop}

Equality (\ref{eq.RR.G.0}) is 
obtained by Vergne \cite{Vergne96}[Part II] in the case of a Hamiltonian 
action of the circle group on a compact symplectic manifold. 
    
\medskip

Suppose now that the decomposition (\ref{eq.J.induit}) holds. The 
trivial bundle $\ggot_{\C}\times\Zcal$ has two irreducible complex 
spinor bundles $\wedge^{\bullet}_{\C}\ggot_{\C}\times\Zcal$ and 
$\wedge^{\bullet}_{J_{\Dcal}}\ggot_{\C}\times\Zcal$. Thus 
\begin{equation}\label{eq.L.D}
\wedge^{\bullet}_{J_{\Dcal}}\ggot_{\C}\times\Zcal = 
\wedge^{\bullet}_{\C}\ggot_{\C}\times\Zcal\otimes \pi^*L_{\Dcal}
\end{equation}
where $\pi^*L_{\Dcal}\to\Zcal$ is the line bundle equal to
$\Hom_{Cl_{\C}}(\wedge^{\bullet}_{\C}\ggot_{\C}\times\Zcal,
\wedge^{\bullet}_{J_{\Dcal}}\ggot_{\C}\times\Zcal)$ : at $z\in\Zcal$, 
$\pi^*L_{\Dcal}\vert_{z}$ is the complex vector space of linear 
maps $\wedge^{\bullet}_{\C}\ggot_{\C}\to\wedge^{\bullet}_{J_{\Dcal}(z)}\ggot_{\C}$ 
commuting with the Clifford actions (see \cite{Schroder}). Note that 
$\wedge^{\pm}_{J_{\Dcal}}\ggot_{\C}\times\Zcal = 
\wedge^{\pm}_{\C}\ggot_{\C}\times\Zcal\otimes \pi^*L_{\Dcal}$ if the 
orientation of $J_{\Dcal}$ coincide with those defined by $\imath$ (i.e. 
$\det\Dcal>0$). If $\det\Dcal<0$, we have 
$\wedge^{\pm}_{J_{\Dcal}}\ggot_{\C}\times\Zcal = 
\wedge^{\mp}_{\C}\ggot_{\C}\times\Zcal\otimes \pi^*L_{\Dcal}$.

\begin{prop}\label{prop.RR.T.0.bis}
    Suppose that the decomposition (\ref{eq.J.induit}) holds, and  let 
    \break 
    $RR^{J_{red}}(\Mcal_{red},-)$ be the quantization map given by 
    $J_{red}$. For every $G$-equivariant vector bundle $E\to M$, we have  
    \begin{equation}\label{eq.RR.G.0.bis}
    \left[RR_{0}^{^{G}}(M,E)\right]^G =\pm\,
    RR^{J_{red}}(\Mcal_{red},E_{red}\otimes L_{\Dcal})\ ,
    \end{equation}
    where $\pm$ is the sign of $\det\Dcal$.
\end{prop}

{\em Proof of Proposition \ref{prop.RR.T.0}} : Following Definition 
\ref{def.RR.beta}, the map $RR_{0}^{^{G}}(M,-)$ is defined by 
$\Thom_{G,[0]}^f(M)\in K_{G}(\T_{G}\Ucal^{^{G,0}})$, where 
$\Ucal^{^{G,0}}$ is a (small) neighbourhood of $\Zcal$ in $M$.  Since 
$0$ is a regular value of $f_{_{G}}$, $\Ucal^{^{G,0}}$ is 
diffeomorphic to $\Zcal\times \ggot^{*}$, and the moment map is equal 
to the projection $f:\Zcal\times \ggot^{*}\to \ggot^{*}$ in a 
neighbourhood of $\Zcal$ in $\Zcal\times \ggot^{*}$.  We denote by 
$\sigma_{\Zcal}\in K_{G}(\T_{G}(\Zcal\times \ggot^{*}))$ the symbol 
corresponding to $\Thom_{G,[0]}^f(M)$ through the diffeomorphism 
$\Ucal^{^{G,0}}\cong \Zcal\times \ggot^{*}$.  Let 
$\indice_{\Zcal\times \ggot^{*}}^{G}: K_{G}(\T_{G}(\Zcal\times 
\ggot^{*}))\to R^{-\infty}(G)$ be the index map on $\Zcal\times 
\ggot^{*}$.  The map $RR_{0}^{^{G}}(M,-)$ is defined by 
$RR_{0}^{^{G}}(M,E)=\indice_{\Zcal\times \ggot^{*}}^{G} 
(\sigma_{\Zcal}\otimes f^{*}(E_{\vert\Zcal}))$.

Following Atiyah \cite{Atiyah.74}[Theorem 4.3], the inclusion map
$j:\Zcal\croc \Zcal\times \ggot^{*}$ induces an $R(G)$-module morphism
$j_{!}:K_{G}(\T_{G}\Zcal)\to K_{G}(\T_{G}(\Zcal\times \ggot^{*}))$, with 
the commutative diagram
\begin{equation}\label{j.point.bis}
\xymatrix@C=2cm{
 K_{G}(\T_{G}\Zcal)\ar[r]^{j_{!}} \ar[dr]_{\indice^{G}_{\Zcal}} & 
 K_{G}(\T_{G}( \Zcal\times \ggot^{*})) 
 \ar[d]^{\indice_{\Zcal\times \ggot^{*}}^{G}}\\
     & R^{-\infty}(G)}.
\end{equation} 

More generally, the map $i_{!}:K_{G}(\T_{G}\Zcal)\to K_{G}(\T_{G}\Ycal)$ is 
defined by Atiyah for any embedding $i:\Zcal\croc\Ycal$ of $G$-manifolds 
with $\Zcal$ compact. 

Consider now the case where $i$ is the zero-section of a 
$G$-vector bundle $\Ecal \to \Zcal$. 
In general the map $i_{!}$ is {\em not} an isomorphism. 
If furthermore the $G$-action is {\em locally free} over $\Zcal$, then
$\T_{G}\Zcal$, $\T_{G}\Ecal$ are respectively  
subbundles of $\T\Zcal\to \Zcal$, $\T\Ecal\to \Ecal$, and 
the projection $\T_{G}\Ecal\to\T_{G}\Zcal$ is a vector bundle
isomorphic to $s^{*}(\T\Ecal)$ (where $s:\T_{G}\Zcal\croc \T\Zcal$ is 
the inclusion). Hence the vector bundle $\T_{G}\Ecal\to\T_{G}\Zcal$
inherits a complex structure over the fibers (coming from the complex
vector bundle $\T\Ecal\to\T\Zcal$). In this situation, the map
$i_{!}:K_{G}(\T_{G}\Zcal)\to K_{G}(\T_{G}\Ecal)$ is the
Thom isomorphism. 

In the case of the (trivial) vector bundle $\Zcal\times 
\ggot^{*}\to\Zcal$, the map $j_{!}:K_{G}(\T_{G}\Zcal)\to$ 
$K_{G}(\T_{G}(\Zcal\times \ggot^{*}))$ is then an {\em isomorphism}.  
Take $\tilde{\sigma}_{\Zcal}=(j_{!})^{-1}(\sigma_{\Zcal})$, and from 
the commutative diagram (\ref{j.point.bis}) we have 
$RR_{0}^{^{G}}(M,E)= \indice_{\Zcal}^{G} 
\left(\tilde{\sigma}_{\Zcal}\otimes E\vert_{\Zcal}\right)$.  From 
Theorem \ref{thm.atiyah.1} we get $$ 
\indice_{\Zcal}^{G}(\tilde{\sigma_{\Zcal}}\otimes E\vert_{\Zcal})= 
\sum_{\mu\in\Lambda^{*}_{+}}\indice_{\Mcal_{red}}(\sigma^{red}\otimes 
E_{red} \otimes\underline{V_{\mu}}^{*}) .  V_{\mu}\ , $$ where 
$\sigma^{red}\in K_{orb}(\T\Mcal_{red})$ corresponds to 
$\tilde{\sigma}_{\Zcal}=(j_{!})^{-1}(\sigma_{\Zcal})$ through the 
isomorphism $\pi^{*}:K_{orb}(\T\Mcal_{red})\to K_{G}(\T_{G}\Zcal)$.  
Proposition \ref{prop.RR.T.0.bis} follows immediately from the
 
\begin{lem}\label{lem.RR.O.canonique}
  We have 
  $$
  j_{!}\circ(\pi)^{*}\Big(\sthom(\Mcal_{red})\Big)=
  \sigma_{\Zcal}
  $$
  in $K_{G}(\T_{G}(\Zcal\times\ggot^{*}))$.
\end{lem} 

{\em Proof :} Let $\Scal(M)$ the irreducible spinor bundle defined by 
the almost complex structure $J$.  Let $\widetilde{J}$ be the almost 
complex structure on $\Zcal\times\ggot^*$, equal to $J$ on $\Zcal$, 
and which is constant on the fibers of the projection $\Zcal\times 
\ggot\to \Zcal$.  Since the almost complex structures $J$ and 
$\widetilde{J}$ are {\em homotopic} near $\Zcal$, the complex 
$\sigma_{\Zcal}$ can be defined on $\Zcal\times \ggot$ with 
$\widetilde{J}$ : we take $\Scal(M)\vert_{\Zcal}\times \ggot^*$ for 
bundle of spinors over $\Zcal\times\ggot^*$.  Following 
(\ref{eq.TM.decompose.bis}) and (\ref{eq.TM.decompose}), for 
$(z,\xi)\in\Zcal\times\ggot^*$ a vector $v\in 
\T_{(z,\xi)}(\Zcal\times\ggot^*)$ decomposes into $v=v_1 +X+\imath Y$, 
where $v_1\in\pi^*(\T M_{\xi})$, and $X+\imath Y\in \ggot_{\C}$.  The 
map $\sigma_{\Zcal}(z,\xi;v)$ acts on $\Scal(M)_z$ by the Clifford 
action pushed by the vector field\footnote{The tangent vector 
$\Hcal^{^G}(z,\xi)\in\ggot_{\Zcal}\vert_z$ is equal to 
$\imath\xi\in\ggot_{\C}\times\Zcal$.} $\Hcal^{^G}(z,\xi)=\imath\,\xi$ :
$$
\sigma_{\Zcal}(z,\xi;v)=\Clif_z(v_1 +X+\imath \,(Y-\xi))\ .
$$
Using now Lemma \ref{lem.spinc.induit}, we see that $\sigma_{\Zcal}$ is homotopic to 
the symbol $\sigma'_{\Zcal}$ which acts on the product 
$(\pi^{*}\Scal(\Mcal_{red})\otimes 
\wedge^{\bullet}_{\C}\ggot_{\C}\times\Zcal)\times\ggot^*$ by 
$$
\sigma'_{\Zcal}(z,\xi;v)=\Clif_z(v_1)\odot \Clif(X+\imath (Y-\xi))\ .
$$
Now we see that the map
$\Clif_{z}(v_{1})\odot \Clif(X+\imath(Y-\xi))$ is homotopic, as a 
$G$-transversally elliptic symbol, to 
$\Clif_{z}(v_{1})\odot \Clif(\xi+\imath X)$. The $K$-theory class of 
this former symbol is equal to $(\pi)^{*}(\sthom(\Mcal_{red}))\odot 
k_!(\C)$ (where $k:\{0\}\croc\ggot^*$) which is the symbol
map of $j_{!}\circ(\pi)^{*}\left(\sthom(\Mcal_{red})\right)$ 
(see the construction of the map $j_{!}$ in \cite{Atiyah.74}[Lecture 4]).
We have shown that $j_{!}\circ(\pi)^{*}\left(\sthom(\Mcal_{red})\right)
=\sigma_{\Zcal}$ in $K_{G}(\T_{G}(\Zcal\times \ggot^*))$. 
$\Box$
 
\bigskip

{\em Proof of Proposition \ref{prop.RR.T.0.bis}} : Here the proof is 
similar to the former proof but we use Lemma \ref{lem.J.induit} 
instead of Lemma \ref{lem.spinc.induit}. One as to show that 
$$
j_{!}\circ(\pi)^{*}\Big(\sthom(\Mcal_{red})\otimes L_{\Dcal}\Big)=\pm
  \sigma_{\Zcal}
$$
in $K_{G}(\T_{G}(\Zcal\times\ggot^{*}))$, where $\pm$ is the sign of 
$\det\Dcal$. By Lemma \ref{lem.J.induit}, we see as before that 
$\sigma_{\Zcal}$ is homotopic to the product
\begin{equation}\label{eq.clif.J.red}
\Clif_{z}(v_{1})\odot \Clif_{J_{\Dcal}}(\xi+\imath X)
\end{equation}
acting on $(\wedge^{\bullet}_{J_{red}}\pi^*(\T\Mcal_{red})
\otimes\wedge^{\bullet}_{J_{\Dcal}}\ggot_{\C}\times\Zcal)\times\ggot^*$.
Now we use the isomorphism of irreducible complex spinor bundles 
(\ref{eq.L.D}) where we have two different orientations 
$o(J_{\Dcal})$ and $o(\imath)$ on $\ggot_{\C}\times\Zcal$:
$o(J_{\Dcal})=\pm o(\imath)$ where $\pm$ is the sign of $\det\Dcal$. 
Hence the transversally elliptic symbol (\ref{eq.clif.J.red}) is 
equal to 
$$
\pm\ 
\Clif_{z}(v_{1})\odot \Clif(\xi+\imath X)\odot {\rm Id}_{L_{\Dcal}}
$$
acting on $(\wedge^{\bullet}_{J_{red}}\pi^*(\T\Mcal_{red})
\otimes\wedge^{\bullet}_{\C}\ggot_{\C}\times\Zcal\otimes 
L_{\Dcal})\times\ggot^*$. $\Box$


\subsection{The map \protect $RR^{^{G}}_{\beta}$ when $G_{\beta}=G$}
\label{subsec.RR.G.beta}

When $\beta\in \Bcal_{_{G}}-\{0\}$ is in the center of $\ggot$,
the map $RR^{^G}_{\beta}(M,-)$ is the Riemann-Roch character localized near
$M^{\beta}\cap f_{_{G}}^{-1}(\beta)$. In this subsection we prove 
that $[RR^{^G}_{\beta}(M,E)]^G=0$ if $E$ is a $f_{_{G}}$-strictly positive 
complex vector bundle.

The almost complex structure $J$ and the abstract
moment map $f_{_{G}} : M\to \ggot$ restrict on $M^{\beta}$ to an almost complex 
structure $J_{\beta}$ and a abstract moment map $f_{_{G}}\vert_{M^{\beta}}$. 
The set $M^{\beta}\cap f_{_{G}}^{-1}(\beta)=
(f_{_{G}}\vert_{M^{\beta}})^{-1}(\beta)$ is a component
of the critical set of $C^{f_{_{G}}\vert_{M^{\beta}}}$, and
we denote by $RR_{\beta}^{^{G}}(M^{\beta},-):K_{G}(M^{\beta})\to 
R^{-\infty}(G)$ the Riemann-Roch character on $M^{\beta}$ 
localized near the component $(f_{_{G}}\vert_{M^{\beta}})^{-1}(\beta)$ (see
Definition \ref{def.RR.beta}).

Here we proceed as in section \ref{sec.loc.M.beta}.
Let $p:\Ncal\to M^{\beta}$ be the normal bundle of $M^{\beta}$ in $M$.
The torus $\tore_{\beta}\croc G$ acts linearly on the fibers of the 
complex vector bundle $\Ncal$, thus we associate, as in Theorem
\ref{th.localisation.pt.fixe}, the polarized complex 
$G$-vector bundles $\Ncal^{+,\beta}$ and $(\Ncal\otimes\C)^{+,\beta}$.

\begin{prop}\label{prop.loc.RR.G.beta}
    For every $E\in K_{G}(M)$, we have the following equality in 
    \break $\widehat{R}(G)$ :
    $$
    RR_{\beta}^{^{G}}(M,E)=(-1)^{r_{\Ncal}}\sum_{k\in\N}
    RR_{\beta}^{^{G}}(M^{\beta},E\vert_{M^{\beta}}\otimes\det
    \Ncal^{+,\beta}\otimes S^k((\Ncal\otimes\C)^{+,\beta}) \ ,
    $$
    where $r_{\Ncal}$ is the locally constant function on $M^{\beta}$
    equal to the complex rank of $\Ncal^{+,\beta}$.
\end{prop}

Consider the $G\times\tore_{\beta}$-Riemann-Roch character
$RR_{\beta}^{^{G\times{\rm T}_{\beta}}}(M^{\beta},-)$ localized near 
$M^{\beta}\cap f_{_{G}}^{-1}(\beta)$. It
can be extended trivially to a map, still denoted by 
$RR_{\beta}^{^{G\times{\rm T}_{\beta}}}(M^{\beta},-)$, from 
$K_{G}(M^{\beta})\,\widehat{\otimes}\, R(\tore_{\beta})$ to
$R^{-\infty}(G)\,\widehat{\otimes}\, R(\tore_{\beta})$.
Following Definition \ref{wedge.V.inverse}
the element $\wedge_{\C}^{\bullet}\overline{\Ncal}\in 
K_{G\times \tore_{\beta}}(M^{\beta})\simeq
K_{G}(M^{\beta})\otimes R(\tore_{\beta})$ admits a polarized inverse
$\left[\wedge_{\C}^{\bullet}\overline{\Ncal}\,\right]^{-1}_{\beta}\, \in\, 
K_{G}(M^{\beta})\,\widehat{\otimes}\,R(\tore_{\beta})$.
Finally the result of Proposition \ref{prop.loc.RR.G.beta} can be 
written as the following equality in 
$R^{-\infty}(G)\,\widehat{\otimes}\, R(\tore_{\beta})$ :

\begin{equation}\label{eq.loc.RR.G.beta.simplifie}
RR^{^G}_{\beta}(M,E)=RR^{^{G\times{\rm T}_{\beta}}}_{\beta}
\left(M^{\beta},E\vert_{M^{\beta}}\otimes
\left[\wedge_{\C}^{\bullet}\overline{\Ncal}\,\right]^{-1}_{\beta}\right)\ .
\end{equation}

\medskip

Consider the decomposition of $RR_{\beta}^{^G}(M,E) 
=\sum_{\lambda}m_{\beta,\lambda}(E)\,\chi_{_{\lambda}}^{_G}$ in 
irreducible characters $\chi_{_{\lambda}}^{_G},\ 
\lambda\in\Lambda_{+}^{*}$.  Let $E$ be a $f_{_{G}}$-strictly positive 
complex vector bundle over $M$, and let $\eta_{_{E,\beta}}> 0$ be the 
constant defined in Definition \ref{eq.mu.positif}.  If $\Zcal$ is a 
connected component of $M^{\beta}$ which intersects 
$f_{_{G}}^{-1}(\beta)$, every weight $a$ of the $\tore_{\beta}$-action 
on the fibers of the complex vector bundle 
$E^{\stackrel{k}{\otimes}}\vert_{\Zcal}\otimes\det 
\Ncal^{+,\beta}\otimes S^k((\Ncal\otimes\C)^{+,\beta})$ satisfy 
$\langle a,\beta\rangle\geq k.\eta_{_{E,\beta}}$.  Lemma 
\ref{lem.multiplicites.tore} and Corollary 
\ref{coro.multiplicites.tore}, applied to this situation, show that 
\begin{equation}\label{eq.m.E.k} 
m_{\beta,\lambda}(E^{\stackrel{k}{\otimes}})\neq 0 \ \Longrightarrow \ 
\langle \lambda,\beta\rangle \geq\, k.\eta_{_{E,\beta}}\ .  
\end{equation}

In particular $[RR_{\beta}^{^G}(M,E)]^G=m_{\beta,0}(E)=0$, so we 
have proved the 

\begin{coro}\label{coro.G.beta.egal.G}
    Let $E$ be a $f_{_{G}}$-{\em strictly positive} complex vector bundle 
    over $M$ (see Def. \ref{eq.mu.positif}). For any 
    $\beta\in\Bcal_{_{G}}-\{0\}$, with $G_{\beta}=G$,
    the $G$-invariant part of $RR_{\beta}^{^{G}}(M,E)$ is equal to $0$.
\end{coro}

\medskip

{\em Proof of Proposition \ref{prop.loc.RR.G.beta}} : 

\medskip

Here we proceed 
as in the proof of Theorem \ref{th.localisation.pt.fixe}. The almost 
complex structure $J$ induces an almost complex structure $J_{\beta}$
on $M^{\beta}$ and a complex structure $J_{\Ncal}$ on the fibers of 
$\Ncal$.  The 
$G\times\tore_{\beta}$-vector bundle 
$p:\Ncal\to M^{\beta}$ is isomorphic
to $R\times_{U}N\to M^{\beta}=R/U$, where $R$ is the 
$\tore_{\beta}$-equivariant unitary frame of $(\Ncal,J_{\Ncal})$ 
framed on $N$.

Let $\Ucal^{^{G,\beta}}$ be a neighbourhood of $C^{^G}_{\beta}$
in $M$, and consider the $G$-transversally elliptic symbol 
$\Thom^f_{G,[\beta]}(M)\in K_{G}(\T_{G}\Ucal^{^{G,\beta}})$ introduced in  
Definition \ref{def.thom.beta.f}. Here we choose $\Ucal^{^{G,\beta}}$ 
diffeomorphic to an open subset of $\Ncal$ of the form
$\Vcal:=\{ n=(x,v)\in\Ncal,\ x\in \Ucal\ {\rm and}\ \vert 
v\vert<\esp\}$, where $\Ucal$ is a neighbourhood of
$(f_{_{G}}\vert_{M^{\beta}})^{-1}(\beta)$ in $M^{\beta}$.
The moment map $f_{_{G}}$, the vector field $\Hcal^{^G}$, and
$\Thom^f_{G,[\beta]}(M)$ are transported by this
diffeomorphism to $\Vcal$ (we keep the same symbol for these 
elements). 

We define now the homogeneous vector field $\widetilde{\Hcal}^{^G}$ on 
$\Ncal$ by
\begin{equation}\label{def.H.tilde}
\widetilde{\Hcal}^{^G}_{n}:=\Big(f_{_{G}}(p(n))\Big)_{\Ncal}(n),
\ n\ \in\ \Ncal\ .
\end{equation}

Using the isomorphism $\T\Ncal\tilde{\to} p^{*} 
(\T M^{\beta}\oplus \Ncal)$ (see (\ref{eq:trivialisation.T.N})) 
the manifold $\Ncal$ is endowed with the almost complex structure 
$\widetilde{J}:=p^{*}(J_{\beta}\oplus J_{\Ncal})$. With the data 
$(\widetilde{J},\ \widetilde{\Hcal}^{^G})$, we construct
the following $G$-transversally elliptic symbol over $\Ncal$ :
\begin{equation}\label{Thom.G.f.modifie}
\Thom^f_{G,[\beta]}(\Ncal)(n,w):=
\Thom_{G}(\Ncal,\widetilde{J})(n,w-\widetilde{\Hcal}^{^{G}}_{n}),\quad
{\rm for}\quad (n,w)\ \in\ \T\Ncal\ .
\end{equation}

Let us now verify that 
$$
\Thom^f_{G,[\beta]}(M)=\Thom^f_{G,[\beta]}(\Ncal)
\quad {\rm in}\  K_{G}(\T_{G}\Vcal)\ .
$$ 

The invariance of the Thom class after the modification of the almost 
complex structure is carried out in Lemma \ref{lem.J.modifie} : the class of 
$\Thom^f_{G,[\beta]}(M)$ is equal in $K_{G}(\T_{G}\Vcal)$ to the class of 
the symbol
$$
\sigma_{1}(n,w):=\Thom_{G}(\Ncal,\widetilde{J})(n,w-\Hcal^{^G}_{n}),\quad 
(n,w)\ \in\ \T\Vcal\ .
$$

Using now the family of vectors field  
$\Hcal_{t}^{^G}(n):=\Big(f_{_{G}}(x,t.v)\Big)_{\Vcal}(n)$, $\, t\in[0,1]$,
$\, n=(x,v)\in \Vcal$, we construct the homotopy
$$
\sigma_{t}(n,w):=\Thom_{H}(\Ncal,\widetilde{J})(n,w-\Hcal_{t}^{^G}(n)),\quad 
(n,w)\ \in\ \T\Vcal
$$ of $G$-transversally elliptic symbol between
$\sigma_{1}$ and $\Thom^f_{G,[\beta]}(\Ncal)$ (one easily verifies that
$\Char(\sigma_{t})\cap\T_{G}\Vcal = C^{^G}_{\beta}$ for every
$t\in[0,1]$). Finally, we have shown that 
$\Thom^f_{G,[\beta]}(\Ncal)=\Thom^f_{G,[\beta]}(M)$ in 
$K_{G}(\T_{G}\Vcal)$, thus
$$
RR_{\beta}^{^G}(E)=\indice_{\Ncal}^{G}\left(\Thom^f_{G,[\beta]}(\Ncal)\otimes
p^{*}(E_{\vert M^{\beta}})\right) 
$$ 
for every $E\in K_{G}(M)$.

\medskip

Now we proceed as follows. For every $(n,w)\in\T\Vcal$, the Clifford 
action \break
$\Thom^f_{G,[\beta]}(\Ncal)(n,w)=Cl_{n}(w-\widetilde{\Hcal}^{^G}_{n})$ 
on $\wedge_{\C}^{\bullet}\T_{n}\Vcal$ is equal to
the exterior product 
\begin{equation}\label{eq.produit.clifford}
    Cl_{x}(w_{1}-[\widetilde{\Hcal}^{^G}_{n}]_{1})
\odot Cl_{x}(w_{2}-[\widetilde{\Hcal}^{^G}_{n}]_{2})
\end{equation}
acting on $\wedge_{\C}^{\bullet}\T_{x}M^{\beta}\otimes
\wedge_{\C}^{\bullet}\Ncal\vert_{x}$, where
$x=p(n)$. Here $w\to w_{1},\ \T_{n}\Vcal\to\T_{x}M^{\beta}$
is the tangent map $\T p\vert_{n}$, and $w\to w_{2}=[w]^{V},
\ T_{n}\Vcal\to\Ncal\vert_{x}$ is the `vertical' map. We see that
$[\widetilde{\Hcal}^{^G}_{n}]_{1}=\Hcal^{^G}_{x}$ is the vector field on 
$M^{\beta}$ generated by the moment map $f_{_{G}}\vert_{M^{\beta}}$
(see Definition \ref{def.H.vector}).

Suppose that the exterior product (\ref{eq.produit.clifford}) can be 
modified in 
\begin{equation}\label{eq.produit.clifford.bis}
Cl_{x}(w_{1}-\Hcal^{^G}_{x})
\odot Cl_{x}(w_{2}-\beta_{\Ncal}\vert_{n}),
\end{equation}
without changing the K-theoretic class. This will prove a modified version of 
(\ref{eq:Thom.egalite}) in
$K_{G\times\tore_{\beta}\times U}
(\T_{G\times\tore_{\beta}\times U}(R\times N))$ :

\begin{equation}
    \pi^{*}_{N}\Thom^f_{G,[\beta]}(\Ncal)=
    \pi^{*}\Thom^f_{G,[\beta]}(M^{\beta})\odot 
    \Thom^{\beta}_{\tore_{\beta}\times U}(N)\ ,
    \label{eq:Thom.egalite.bis}
\end{equation}
where $\pi_{N}:R\times N\to R \times_{U}N=\Ncal$, 
$\pi:R\to R/U= M^{\beta}$ are the quotient maps relative 
to the free $U$-action, and $\odot$ is the product

\begin{equation}\label{eq:G.T.beta.U.bis}
   K_{G\times U}(\T_{G\times U}R)
   \times
   K_{\tore_{\beta}\times U}(\T_{\tore_{\beta}}N)
   \longrightarrow
   K_{G\times\tore_{\beta}\times U}
   (\T_{G\times\tore_{\beta}\times U}(R\times N)) .
\end{equation}

The symbols $\Thom^f_{G,[\beta]}(\Ncal)$,
$\Thom^f_{G,[\beta]}(M^{\beta})$ and $\Thom^{\beta}_{\tore_{\beta}\times U}(N)$
belong respectively to 
$K_{G\times\tore_{\beta}}
(\T_{G\times\tore_{\beta}}(R\times_{U} N))$, 
$K_{G}(\T_{G}(R/U))$, and 
$K_{\tore_{\beta}\times U}(\T_{\tore_{\beta}\times U}N)$. The  Proposition 
\ref{prop.loc.RR.G.beta} follows after taking the index, 
and the $U$-invariants, in (\ref{eq:Thom.egalite.bis}).

\medskip

Finally we explain why the change of $[\widetilde{\Hcal}^{^G}_{n}]_{2}$ 
in $\beta_{\Ncal}\vert_{n}$ can be done in 
(\ref{eq.produit.clifford}) without changing the class of 
$\Thom^f_{G,[\beta]}(\Ncal)$.

\medskip

Let $\mu^{\Ncal}:\ggot\to \Gamma(M^{\beta},\End(\Ncal))$ be the `moment'  
relative to the choice of a connection on $\Ncal\to M^{\beta}$
(see Definition 7.5 in \cite{B-G-V}). Then, for every $X\in\ggot$ we have 
$$
[X_{\Ncal}(x,v)]^{V}=-\mu^{\Ncal}(X)\vert_{x}.v,\quad
(x,v)\ \in\ \Ncal 
$$
(see Proposition 7.6 in \cite{B-G-V}). When $X=\beta$, the vector 
field $\beta_{\Ncal}$ is vertical, hence we have $\mu^{\Ncal}(\beta)\vert_{x}.v
=\Lcal^{\Ncal}(\beta)\vert_{x}.v=-\beta_{\Ncal}(x,v)$,
where $\Lcal^{\Ncal}(\beta)$ is the infinitesimal action of $\beta$
on the fiber of $\Ncal\to M^{\beta}$. We have also 
$[\widetilde{\Hcal}^{^G}_{n}]_{2}=-\mu^{\Ncal}(f_{_{G}}(x))\vert_{x}.v$, 
for every $n=(x,v)\in\Ncal$.

Note that the quadratic form
$v\in\Ncal_{x}\to\vert\Lcal^{\Ncal}(\beta)\vert_{x}.v\vert^{2}$ 
is positive definite for $x\in M^{\beta}$. Hence, for every
$X\in\ggot$ close enough to $\beta$, the quadratic form
$v\in\Ncal_{x}\to(\mu^{\Ncal}(\beta)\vert_{x}.v,
\mu^{\Ncal}(X)\vert_{x}.v)$ 
is positive definite for $x\in M^{\beta}$.

Consider now the homotopy
$$
\sigma^{t}(n,w):= Cl_{x}(w_{1}-\Hcal^{^G}_{x})
\odot Cl_{x}(w_{2}-t.[\widetilde{\Hcal}^{^G}_{n}]_{2}
-(1-t).\beta_{\Ncal}\vert_{n}),\quad (n,v)\in\Vcal\quad t\in[0,1].
$$
We see that $(n,w)\in\Char(\sigma^{t})\cap\T_{G}\Vcal$ if and only if

\noindent i) $w_{1}=\Hcal^{^G}_{x}$, $w_{2}=t[\widetilde{\Hcal}^{^G}_{n}]_{2}
+(1-t)\beta_{\Ncal}(n)$, and 

\noindent ii) $(w_{1},X_{M^{\beta}}(x)) 
+(w_{2},[X_{\Ncal}(x,v)]^{V})=0$ for all $X\in\ggot$.

Take now $X=f_{_{G}}(x)$ in ii). Using i), we get
\begin{equation}\label{eq.somme.positive}
\left| \Hcal^{^G}_{x}\right|^{2}+
t.\left| \mu^{\Ncal}(f_{_{G}}(x))\vert_{x}.v\right|^{2}+
(1-t).\Sigma(x,v)= 0\ ,
\end{equation}
with $\Sigma(x,v):=
(\mu^{\Ncal}(\beta)\vert_{x}.v,\mu^{\Ncal}(f_{_{G}}(x))\vert_{x}.v)$.

If $x\in M^{\beta}$ is sufficiently close to 
$(f_{_{G}}\vert_{M^{\beta}})^{-1}(\beta)$ , the term $\Sigma(x,v)$ is 
positive for all $v\in\Ncal_{x}$. 
In this case, (\ref{eq.somme.positive}) gives 
$\Hcal^{^G}_{x}=0$ and $\Sigma(x,v)=0$, which insures
that $x\in C^{^G}_{\beta}$ and $v=0$. 

We have proved that $\Char(\sigma^{t})\cap\T_{G}\Vcal=C^{^G}_{\beta}$ 
for every $t\in [0,1]$ if $\Vcal$ is `small' enough.
Hence $\sigma^{t}$ is an homotopy of $G$-transversally
elliptic symbols over $\T\Vcal$ between the exterior products 
(\ref{eq.produit.clifford}) and (\ref{eq.produit.clifford.bis}). 
$\Box$


\subsection{Induction formula}\label{subsec.induction.G.H}

This section is concerned by an induction formula which compare
the map $RR_{\beta}^{^{G}}(M,-)$ with the similar
localized Riemann-Roch characters defined for the maximal torus, and 
the stabilizer $G_{\beta}$. 
The idea of this induction comes from a previous paper of the author
\cite{pep2} where a similar induction formula in the context of 
equivariant cohomology was proved. 

\medskip

Consider the restriction $f_{_{H}}:M\to \hgot$ of the moment map
$f_{_{G}}$ to the maximal torus $H$. In this situation we use the 
vector field $\Hcal^{^{H}}\vert_{m}=
f_{_{H}}(m)_{M}\vert_{m}, m\in M$ to decompose the map
$RR^{^{H}}(M,-):K_{H}(M)\to R(H)$ near the set 
$C^{f_{_{H}}}=\{\Hcal^{^{H}}=0\}$. From Lemma \ref{lem.C.f.G}
there exists a finite subset $\Bcal_{_{H}}\subset \hgot$, such that 
$C^{f_{_{H}}}=\bigcup_{\beta\in\Bcal_{_{H}}}C^{^{H}}_{\beta}$, with 
$C^{^{H}}_{\beta}= M^{\beta}\cap f_{_{H}}^{-1}(\beta)$.
As in Definition \ref{def.RR.beta}, we define for
every $\beta\in \Bcal_{_{H}}$, the map $RR_{\beta}^{^{H}}(M,-):
K_{H}(M)\to R^{-\infty}(H)$ which is 
the Riemann-Roch character localized near $C^{^{H}}_{\beta}$.

Let $W$ be the Weyl group of $(G,H)$. Note that $\Bcal_{_{H}}$ is a
$W$-stable subset of $\hgot$, and that $\Bcal_{_{G}}\subset\Bcal_{_{H}}
\cap\hgot_{+}$.

\begin{theo}\label{th.induction.G.H}
We have, for every $\beta\in\Bcal_{_{G}}$, the following induction 
formula between $RR_{\beta}^{^{G}}(M,-)$ and $RR_{\beta}^{^{H}}(M,-)$.
For every $E\in K_{G}(M)$, we have\footnote{See Equations 
(\ref{eq.holomorphe.G.H}) and (\ref{eq:Hol-G-beta}) in Appendix B
 for the definition of the holomorphic induction maps 
$\HolH$ and $\HolB$.}
\begin{eqnarray*}
RR_{\beta}^{^{G}}(M,E)&=&\frac{1}{\vert W_{\beta}\vert}\HolH\left(
RR_{\beta}^{^{H}}(M,E)\,\wedge_{\C}^{\bullet}\overline{\ggot/\hgot}\right)\\
&=&\frac{1}{\vert W_{\beta}\vert}\sum_{w\in W}\HolH\left(
w.RR_{\beta}^{^{H}}(M,E)\right)\\
&=& \sum_{\beta'\in W.\beta}\HolH\left(RR_{\beta'}^{^{H}}(M,E)\right)
\end{eqnarray*}
where $W_{\beta}$ is the stabilizer of $\beta$ in $W$.
\end{theo}

We can use the previous induction formula between $G$ and $H$ index 
maps to produce an induction formula between $G$ and $G_{\beta}$ index 
maps. Consider the restriction $f_{_{G_{\beta}}}:M\to\ggot_{\beta}$ 
of the moment map to the stabiliser $G_{\beta}$ of $\beta$ in $G$. 
Let $RR_{\beta}^{^{G_{\beta}}}(M,-)$ be 
the Riemann-Roch character localized near $C^{^{G_{\beta}}}_{\beta}=
M^{\beta}\cap f_{_{G}}^{-1}(\beta)$\footnote{Note that 
$M^{\beta}\cap f_{_{G_{\beta}}}^{-1}(\beta)=
M^{\beta}\cap f_{_{G}}^{-1}(\beta)$ because $f_{_{G_{\beta}}}=
f_{_{G}}$ on $M^{\beta}$.}.

\begin{coro}\label{coro.induction.G.G.beta}
For every $\beta\in\Bcal_{_{G}}$ and every $E\in K_{G}(M)$, we have
$$
RR_{\beta}^{^{G}}(M,E)=\HolB\left(
RR_{\beta}^{^{G_{\beta}}}(M,E)\,\wedge_{\C}^{\bullet}
\overline{\ggot/\ggot_{\beta}}\right)
\quad{\rm in}\quad R^{-\infty}(G)\ .
$$
\end{coro}

{\em Proof of the Corollary }: It comes immediately by applying the 
induction formula of Theorem \ref{th.induction.G.H} to the couples
$(G,H)$ and $(G_{\beta},H)$.

\begin{coro}\label{coro.multiplicites.beta}
Let $E$ be a $f_{_{G}}$-{\em strictly positive} complex vector bundle over 
$M$ (see Def. \ref{eq.mu.positif}). We have
$[RR_{\beta}^{^{G}}(M,E^{\stackrel{k}{\otimes}})]^{G}=0$,
if $k.\eta_{_{E,\beta}}>\langle \theta,\beta\rangle$. Here 
$\theta=\sum_{\alpha>0}\alpha$ is the sum of the 
positive roots of $G$, and $\eta_{_{E,\beta}}$ is the strictly 
positive constant defined in Definition \ref{eq.mu.positif}.
\end{coro}

\medskip
 
{\em Proof of Corollary \ref{coro.multiplicites.beta}} : 

Let us first write the decomposition\footnote{We choose a set 
$\Lambda_{+,\beta}^{*}$ of dominant weight for $G_{\beta}$ that 
contains the set $\Lambda_{+}^{*}$ of
dominant weight for $G$.} 
$RR_{\beta}^{^{G_{\beta}}}(M,E^{\stackrel{k}{\otimes}})=
\sum_{\lambda\in\Lambda_{\beta}^{+}}
m_{\lambda,\beta}(E^{\stackrel{k}{\otimes}})
\chi_{_{\lambda}}^{^{G_{\beta}}}$, 
in irreducible character of $G_{\beta}$.
We know from (\ref{eq.m.E.k}) that
$m_{\lambda,\beta}(E^{\stackrel{k}{\otimes}})\neq 0 \ 
\Longrightarrow \ 
\langle \lambda,\beta\rangle \geq k.\eta_{_{E,\beta}}$. 
Each irreducible character 
$\chi_{_{\lambda}}^{^{G_{\beta}}}$ is equal to $\HolHB(h^{\lambda})$, 
so from Corollary \ref{coro.induction.G.G.beta} we have 
$RR_{\beta}^{^{G}}(M,E^{\stackrel{k}{\otimes}})=
\HolH\Big((\sum_{\lambda}m_{\lambda,\beta}(E^{\stackrel{k}{\otimes}})\, 
h^{\lambda})\Pi_{\alpha\in\Delta(\ggot/\ggot_{\beta})}(1-h^{-\alpha})\Big)$
where $\Delta(\ggot/\ggot_{\beta})$ is the set of $H$-weight on
$\ggot/\ggot_{\beta}$\footnote{The complex structure on 
$\ggot/\ggot_{\beta}$ is defined by $\beta$, so that 
$\langle\alpha,\beta\rangle>0$ for all $\alpha\in 
\Delta(\ggot/\ggot_{\beta})$.}. Finally , we see that 
$RR_{\beta}^{^{G}}(M,E^{\stackrel{k}{\otimes}})$
is a sum of terms of the form 
$m_{\lambda,\beta}(E^{\stackrel{k}{\otimes}})\,
\HolH(h^{\lambda-\alpha_{I}})$ where $\alpha_{I}=\sum_{\alpha\in 
I}\alpha$ and $I$ is a subset of $\Delta(\ggot/\ggot_{\beta})$.

We know from Appendix B that 
$\HolH(h^{\lambda'})$ is either $0$ or the character of an 
irreducible representation; in particular $\HolH(h^{\lambda'})$ is 
equal to $\pm 1$  only if 
$\langle\lambda',X\rangle\leq 0$ for every $X\in\hgot_{+}$ 
(see Remark \ref{hol.egale.1}). So 
$[RR_{\beta}^{^{G}}(M,E^{\stackrel{k}{\otimes}})]^G\neq 0$ only if there 
exists a weight $\lambda$ such that 
$m_{\lambda,\beta}(E^{\stackrel{k}{\otimes}})\neq 0$ and 
$\HolH(h^{\lambda-\alpha_{I}})=\pm 1$. The first condition imposes 
$\langle \lambda,\beta\rangle \geq k.\eta_{_{E,\beta}}$ and the second gives 
$\langle \lambda,\beta\rangle \leq \langle \alpha_{I},\beta\rangle$, 
and combining the two we end with $k.\eta_{_{E,\beta}}\leq \langle 
\alpha_{I},\beta\rangle\leq 
\sum_{\alpha\in \Delta(\ggot/\ggot_{\beta})}\langle\alpha,\beta\rangle
=\langle \theta,\beta\rangle$. We have proved that 
$[RR_{\beta}^{^{G}}(M,E^{\stackrel{k}{\otimes}})]^G = 0$ 
if $k.\eta_{_{E,\beta}} >\langle \theta,\beta\rangle$. $\Box$

\medskip

{\em Proof of Theorem \ref{th.induction.G.H}} : 

\medskip

The first two equalities of the Theorem can be deduced from the third 
one, that is $RR_{\beta}^{^{G}}(M,E)= \sum_{\beta'\in 
W.\beta}\HolH\left(RR_{\beta'}^{^{H}}(M,E)\right)$.  First, it is easy 
to see that \break $RR_{w.\beta}^{^{H}}(M,E)=w.RR_{\beta}^{^{H}}(M,E)$ 
for every $w\in W$ and $\beta\in \Bcal_{_{H}}$.  After, the relation 
$\HolH(\phi\,\wedge_{\C}^{\bullet}\overline{\ggot/\hgot})=\sum_{w\in 
W} \HolH(w.\phi)$, which is true for every $\phi\in R^{-\infty}(H)$ 
(see Remark \ref{wedge.C-wedge.R}), gives the first equality of the 
Theorem.

The map $RR_{\beta}^{^{G}}(M,-)$ is defined 
through the symbol
$\Thom_{G,[\beta]}^f(M)\in K_{G}(\T_{G}\Ucal^{^{G,\beta}})$ where
$i^{^{G,\beta}}:\Ucal^{^{G,\beta}}\to M$ is any $G$-invariant neighbourhood 
of $C^{^{G}}_{\beta}$ such that $\overline{\Ucal^{^{G,\beta}}}\cap 
C^{f_{_{G}}}=C^{^{G}}_{\beta}$ (see Definition \ref{def.thom.beta.f}). 
We define in the same way the localized Thom complex 
$\Thom_{H,[\beta]}^f(M)\in K_{H}(\T_{H}\Ucal^{^{H,\beta}})$.

For notational convenience, we will note in the same way the direct image of 
$\Thom_{G,[\beta]}^f(M)$ (resp. $\Thom_{H,[\beta]}^f(M)$) in
$K_{G}(\T_{G}M)$ (resp. $K_{H}(\T_{H}M)$) via 
$i^{^{G,\beta}}_{*}: K_{G}(\T_{G}\Ucal^{^{G,\beta}})\to 
K_{G}(\T_{G}M)$ (resp. $i^{^{H,\beta}}_{*}: K_{H}(\T_{H}\Ucal^{^{H,\beta}})\to 
K_{H}(\T_{H}M)$). 

Then we have $RR_{\beta}^{^{G}}(M,E)=\indice^{G}_{M}(\Thom_{G,[\beta]}^f(M)
\otimes E)$ for $E\in K_{G}(M)$. The Weyl group acts on $K_{H}(\T_{H}M)$
and we remark that $w.\Thom_{H,[\beta]}^f(M)=\Thom_{H,[w.\beta]}^f(M)$
for every $\beta\in \Bcal_{_{H}}$, and $w\in W$. After taking the 
index we see that $RR_{w.\beta}^{^{H}}(M,E)=w.RR_{\beta}^{^{H}}(M,E)$ 
for every $G$-vector bundle $E$.

Consider the map $r^{\gamma}_{_{G,H}}:
K_{G}(\T_{G}M)\to K_{H}(\T_{H}M)$ defined with $\gamma\in\hgot$ in the 
interior of the Weyl chamber, so that $G_{\gamma}=H$ 
(see subsection \ref{subsec.reduction}). The third equality of the Theorem 
is an immediate consequence of the next Lemma.
 
\begin{lem}\label{lem.rest.thom.beta.f}
We have 
$$
r^{\gamma}_{_{G,H}}\left(\Thom_{G,[\beta]}^f(M)\right)=
\sum_{\beta'\in W.\beta}\Thom_{H,[\beta']}^f(M)\otimes 
\wedge_{\C}^{\bullet}\ggot/\hgot
\quad {\rm in}\quad K_{H}(\T_{H}M)\ .
$$
\end{lem}

\medskip

{\em Proof of Lemma \ref{lem.rest.thom.beta.f} :} 

Consider a
$G$-invariant open neighbourhood $\Ucal^{^{G,\beta}}$ of $C^{^{G}}_{\beta}$ 
such that  $\overline{\Ucal^{^{G,\beta}}}\cap C^{f_{_{G}}}=
C^{^{G}}_{\beta}$. We know from Proposition \ref{prop.restriction.bis} 
that the class $r^{\gamma}_{_{G,H}}(\Thom_{G,[\beta]}^f(M))$ is 
represented by the restriction to $\T\Ucal^{^{G,\beta}}$ of the symbol 
$$
\sigma_{I}(m,v)=Cl_{m}(v-\Hcal^{^{G}}_{m})\odot Cl(\mu_{_{G/H}}(v)),\quad
\quad (m,v)\in \T M\ .
$$
Here $\mu_{_{G/H}}:\T M\to \ggot/\hgot$ is the $\ggot/\hgot$ part of
the Hamiltonian moment map $\mu_{_{G}}:\T M\to \ggot$. Let $f_{_{G/H}}:
M\to\ggot/\hgot$ (resp. $f_{_{H}}:M\to\hgot$) be
the $\ggot/\hgot$-part (resp. the $\hgot$-part) of the moment 
map $f_{_G}$. We will use in our proof the relation 
\begin{equation}\label{relation.p.s}
  (\mu_{_{G/H}}(\Hcal^{^{G}}), f_{_{G/H}})_{_{\ggot}}=
||\Hcal^{^{G}}||^{2}_{_{M}}-(\Hcal^{^{G}},\Hcal^{^{H}})_{_{M}}\ .
\end{equation}
Consider the family of $H$-equivariant symbols 
$\sigma_{\theta},\ \theta\in [0,1]$ 
defined on $\T M$ by 
$$
\sigma_{\theta}(m,v)=Cl_{m}(v-\Hcal^{^{G}}_{m})\odot 
Cl\left(\theta\mu_{_{G/H}}(v)+(1-\theta)f_{_{G/H}}(m)\right),\quad
\quad (m,v)\in \T M\ .
$$
We see that $(m,v)\in\Char(\sigma_{\theta}) \Longleftrightarrow 
v=\Hcal^{^{G}}_{m}\ {\rm and}\ 
\theta\mu_{_{G/H}}(\Hcal^{^{G}}_{m})+(1-\theta)f_{_{G/H}}(m)=0$.  
Combining (\ref{relation.p.s}) with the fact that the vector field 
$\Hcal^{^{H}}$ belongs to the $H$-orbits, we see that 
$\Char(\sigma_{\theta})\cap\T_{H}M\subset \{\Hcal^{^G}=0\}$, for every 
$\theta\in [0,1]$.  By this way we have proved that 
$\sigma_{I}\vert_{\Ucal^{^{G,\beta}}}$ is homotopic to the 
$H$-transversally elliptic symbol 
$\sigma_{II}\vert_{\Ucal^{^{G,\beta}}}$ where 
$$ 
\sigma_{II}(m,v)=Cl_{m}(v-\Hcal^{^{G}}_{m})\odot 
Cl(f_{_{G/H}}(m)),\quad \quad (m,v)\in \T M\ .  
$$
We transform now $\sigma_{II}$ via the following homotopy of
$H$-transversally elliptic symbols 
$$
\sigma^u(m,v):=Cl_{m}(v-\Hcal^{^{H}}_{m}- u.\Hcal^{^{G/H}}_{m} )
\odot Cl(f_{_{G/H}}(m)),\quad (m,v)\in \T M\ ,
$$
for  $u\in [0,1]$. Here $\Char(\sigma^u)\cap\T_{H}M=
\{\Hcal^{^G}=0\}\cap\{ f_{_{G/H}}=0\}$ for all $u\in [0,1]$, hence 
$\sigma_{II}\vert_{\Ucal^{^{G,\beta}}}$ is homotopic to the 
$H$-transversally elliptic symbol 
$\sigma_{III}\vert_{\Ucal^{^{G,\beta}}}$ 
where
$$
\sigma_{III}(m,v)=Cl_{m}(v-\Hcal^{^H}_{m})\odot 
Cl(f_{_{G/H}}(m)),\quad
 \quad (m,v)\in \T M\ .
$$
At this stage we have proved that 
$\sigma_{I}\vert_{\Ucal^{^{G,\beta}}}= 
\sigma_{III}\vert_{\Ucal^{^{G,\beta}}}$ in 
$K_{H}(\T_{H}\Ucal^{^{G,\beta}})$.  Note that

\begin{eqnarray*}
\Char(\sigma_{III}\vert_{\Ucal^{^{G,\beta}}})\cap\T_{H}\Ucal^{^{G,\beta}}
&=& G.(M^{\beta}\cap f_{_{G}}^{-1}(\beta))\bigcap \{ f_{_{G/H}}=0\}\\
&=& W.(M^{\beta}\cap f_{_{G}}^{-1}(\beta))\ ,
\end{eqnarray*}
because $G.\beta\cap \hgot= W.\beta$. Let 
$i: \Ucal^{^{G,\beta}}\croc\Ucal$ be a $H$-invariant neighbourhood
of $W.(M^{\beta}\cap f_{_{H}}^{-1}(\beta))$ 
such that $\overline{\Ucal}\cap \{\Hcal^{^H}=0\}=
W.(M^{\beta}\cap f_{_{H}}^{-1}(\beta))$. The 
symbol $\sigma_{III}\vert_{\Ucal}$ is $H$-transversally elliptic and 
\begin{equation}\label{eq.U.G.beta.U}
i_{*}(\sigma_{III}\vert_{\Ucal})=
\sigma_{III}\vert_{\Ucal^{^{G,\beta}}}=\sigma_{I}\vert_{\Ucal^{^{G,\beta}}} 
\quad {\rm in}
\quad K_{H}(\T_{H}\Ucal^{^{G,\beta}})\ .
\end{equation}
As in the proof of Proposition
\ref{prop.localisation}, (\ref{eq.U.G.beta.U}) is an 
immediate consequence of the excision property. 

The symbol
$(m,v) \to Cl_{m}(v-\Hcal^{^{H}}_{m})$ is $H$-transversally
elliptic on $\T\Ucal$, and equal (by definition) to 
$\sum_{\beta'\in W.\beta}\Thom_{H,[\beta']}^f(M)$. 
Hence $\sigma_{III}\vert_{\Ucal}$ is 
homotopic, in $K_{H}(\T_{H}\Ucal)$, to 
$(m,v) \to Cl_{x}(v-\Hcal^{^{H}}_{m})
\odot 0_{\ggot/\hgot}$, where
$0_{\ggot/\hgot}$ is the zero map from $\wedge_{\C}^{even}\ggot/\hgot$ 
to $\wedge_{\C}^{odd}\ggot/\hgot$. Finally we have shown that
$\sigma_{III}\vert_{\Ucal}=\sum_{\beta'\in W.\beta}
\Thom_{H,[\beta']}^f(M)\otimes \wedge_{\C}^{\bullet}\ggot/\hgot$ in 
$K_{H}(\T_{H}\Ucal)$, and then (\ref{eq.U.G.beta.U}) 
finishes the proof.
$\Box$

\bigskip


\section{The Hamiltonian case}\label{sec.Hamiltonien}

\medskip

In this section, we assume that $(M,\omega)$ is a compact 
symplectic manifold with a Hamiltonian action of a compact 
connected Lie group $G$. The corresponding moment map 
$\mu_{_{G}} :M\to \ggot^{*}$ 
is defined by 
\begin{equation}\label{eq:def.application.moment}
d\langle\mu_{_G},X\rangle=-\,\omega(X_{M},-),\quad \forall\ X\in\ggot.
\end{equation}

The symplectic $2$-form $\omega$ insures the existence of a 
$G$-invariant almost complex structure $J$ {\em compatible} with 
$\omega$, i.e, such that :
$$
(v,w)\to\omega_{x}(v,J_{x}w),\quad v,w \in \T_{x}M\ 
$$
is symmetric and positive definite for all $x\in M$. 
We fix once and for all a $G$-invariant {\em compatible} almost complex 
structure $J$, and we denote by $(-,-)_{_{M}}:=\omega(-,J-)$ the 
corresponding Riemannian metric. Let 
$RR^{^{G}}(M,-)$ be the quantization map defined with 
the {\em compatible} almost complex structures $J$. Since two 
compatible almost complex structure are homotopic 
\cite{Salamon-McDuff}, the map $RR^{^{G}}(M,-)$ does not depend of 
this choice (see Lemma \ref{lem.inv.homotopy}). 

Here the vector field $\Hcal^{^{G}}$ is the Hamiltonian vector field 
of the function\footnote{Equality \ref{eq:def.application.moment} 
gives $\frac{-1}{2}d||\mu_{_{G}}||^{2}=\omega(\Hcal^{^{G}},-)$}
$\frac{-1}{2}||\mu_{_{G}}||^{2}:M\to\R$,
and $\{\Hcal^{^{G}}=0\}$ is the set of critical points of 
$||\mu_{_{G}}||^2$. We know from 
the beginning of section \ref{sec.Localisation.f} that we 
have the decomposition $RR^{^{G}}(M,-)=$ \break 
$\sum_{\beta\in\Bcal_{G}}RR^{^{G}}_{\beta}(M,-)$, 
where $RR^{^{G}}_{\beta}(M,-): K_{G}(M)\to R^{-\infty}(G)$ is 
the Riemann-Roch character localized near the critical set
$C_{\beta}^{^{G}}=G(M^{\beta}\cap\mu_{_G}^{-1}(\beta))$.
In this section we prove  the following
Theorem for the $\mu_{_{G}}$-{\em positive} vector bundles 
(see Def. \ref{eq.mu.positif}).

\begin{theo}\label{QR=RQ-regulier}
Let $E\to M$ be a $G$-equivariant vector bundle over $M$.
For all $\beta\in\Bcal_{_G}-\{0\}$, the $G$-invariant part of 
$RR^{^{G}}_{\beta}(M,E)$ is equal to $0$ if  $E$ is 
$\mu_{_{G}}${\em -positive} and $\mu_{_{G}}^{-1}(0)\neq\emptyset$, or
if $E$ is $\mu_{_{G}}$-{\em strictly positive}. If $0$ is a regular 
value of $\mu_{_{G}}$, the $G$-invariant part of $RR^{^{G}}_{0}(M,E)$ 
is equal to $RR(\Mcal_{red},E_{red})$.
\end{theo}

In subsection \ref{non-regulier}, we consider the general case where 
$0$ is not necessarily a regular value of $\mu_{_{G}}$, and  
$E=L$ a moment bundle for $\mu_{_{G}}$ (see Def. \ref{moment.bundle}). 
With our $K$-theoritic approach we recover the following 

\begin{theo} [Meinrenken-Sjamaar]\label{QR=RQ-singulier}
Let $L\to M$ be a $\mu_{_{G}}$-moment bundle, and let $\tau$ be the 
principal face of $M$.  The 
$G$-invariant part of $RR^{^{G}}(M,L)$ is equal to 
$RR(\Mcal_{a},L_{a})$ for every generic value of 
$\tau\cap\mu_{_{G}}(M)$ sufficiently close to $0$ (see 
subsection \ref{non-regulier} for the notations). 
\end{theo}

\medskip 

\subsection{The map $RR^{^{G}}_{0}$.}\label{subsec.RR.O.hamilton}

We assume that $0$ is a regular value of $\mu_{_{G}}$. The orbifold 
space $\Mcal_{red}:=\mu_{_{G}}^{-1}(0)/G$ inherits a symplectic 
structure $\omega_{red}$. Let $\Dcal(X)=-d\,\mu_{_{G}}(J(X_{M}))$ be 
the endomorphism of the trivial bundle $\mu_{_{G}}^{-1}(0)\times 
\ggot$ defined in (\ref{eq.Dcal}). The compatibility of $J$ with 
$\omega$ gives 
$$
(\Dcal(X),X)= \omega(X_{M},J(X_{M}))_{_{M}}
=\parallel X_{M}\parallel^{2}\ ,
$$
thus decomposition (\ref{eq.J.induit}) holds. A small check shows that 
the induced almost complex structure $J_{red}$ on $\Mcal_{red}$ is 
compatible with $\omega_{red}$. Moreover 
$t\mapsto t\Dcal +(1-t)Id$ is an homotopy of invertible maps between 
$\Dcal$ and the identity, hence the line bundle $L_{\Dcal}\to\Mcal_{red}$ 
defined in (\ref{eq.L.D}) is trivial. The map $RR^{^G}_{0}$ is 
determined by the Proposition \ref{prop.RR.T.0.bis} ; in particular 
$$
    \left[RR_{0}^{^{G}}(M,E)\right]^G =
    RR^{J_{red}}(\Mcal_{red},E_{red})\ ,
$$ 
for any $E\in K_{G}(M)$.

\subsection{The map $RR^{^{G}}_{\beta}$ when $G_{\beta}=G$.}

When $\beta\in\Bcal_{_{G}}-\{0\}$ is in the center of $\ggot$, we 
proved in Corollary \ref{coro.G.beta.egal.G}, that the $G$-invariant 
part of $RR_{\beta}^{^{G}}(M,E)$ is equal to $0$ when $E$ is 
$\mu_{_{G}}$-{\em strictly} positive.  In the Hamiltonian case we 
extend this result for the $\mu_{_{G}}$-positive bundles.

\begin{lem}\label{lem.N.beta.plus}
Let $(\Xcal,\omega)$ be a connected symplectic manifold with a 
$G$-action, and a proper moment map 
$\mu :\Xcal\to \ggot$. Let $J$ be a $G$-invariant almost complex 
structure on $\Xcal$ compatible with $\omega$. Let $\beta$ be 
a $G$-invariant element in a Weyl chamber $\hgot_{+}$ of the 
Lie group $G$, such that 
$\Xcal^{\beta}\cap \mu^{-1}(\beta)\neq\emptyset$. 
Let $\Ncal^{+,\beta}$ be the polarized normal bundle of 
$\Xcal^{\beta}$ in $\Xcal$ (see Def. \ref{wedge.V.inverse} and 
Theorem \ref{th.localisation.pt.fixe}).

If $\Ncal^{+,\beta}=0$, we have 
$$
\mu(\Xcal)\cap\hgot_{+}\subset\{ X\in\hgot_{+},\ (X,\beta)\geq 
\parallel \beta\parallel^2\, \}\ ,
$$
implying in particular that $\parallel \beta\parallel^2$ is the minimal 
value of $\parallel \mu\parallel^2$ on $\Xcal$.
\end{lem}

{\em Proof of the Lemma :}  Let $\Zcal$ be a
connected component of $\Xcal^{\beta}$ which intersects
$\mu^{-1}(\beta)$, and consider the 
set of weights $\{\alpha_{i},\ i\in I\}$ for the action of 
$\tore_{\beta}$ 
on the fibers of the vector bundle $\Ncal\to \Zcal$. We have then the 
following description of the function $(\mu,\beta)$ in the 
neighbourhood of $\Zcal$. For $v\in\Ncal_{x}$, with the 
decomposition $v=\oplus_{i}v_{i}$,  we have for $\vert v\vert$ small
enough
$(\mu,\beta)_{(x,v)}=\vert\beta\vert^{2}-\frac{1}{2}
\sum_{i\in I}\langle\alpha_{i},\beta\rangle \vert v_{i}\vert^{2}$.
If $\langle\alpha_{i},\beta\rangle<0$ for every $i\in I$, we have
\begin{equation}\label{eq.mu.beta}
    (\mu,\beta)\geq\parallel \beta\parallel^2\quad {\rm in\ a\ 
    neighbourhood}\ \Vcal\ {\rm of}\ \Zcal.
\end{equation}
As $\mu^{-1}(\beta)$ is connected and intersect $\Zcal$, the last 
inequality imposes $\mu^{-1}(\beta)\subset \Zcal$. Take 
$X\in \mu(\Xcal)\cap\hgot_{+}$, and consider 
$\Kcal:=\mu^{-1}([X,\beta])$. From the convexity Theorem
\cite{Atiyah.82,Guillemin-Sternberg82.bis,Kirwan.84.bis,L-M-T-W}, the 
set $\Kcal$ is connected. Then $\Vcal\cap\Kcal$ contains, but is not 
equal to $\mu^{-1}(\beta)$ : there exists $m\in \Vcal\cap\Kcal$ with
$\mu(m)\in [X,\beta)$. So $\mu(m)=\beta + t(X-\beta)$ with $t>0$, and 
$(\mu(m),\beta)\geq\parallel \beta\parallel^2$. This two conditions 
imply that $(X,\beta)\geq\parallel \beta\parallel^2$. $\Box$


\begin{lem}\label{lem.RR.beta.G.invariant}
    Let $\beta\in\Bcal_{_{G}}-\{0\}$ be a G-invariant element such 
    that $\parallel \beta\parallel^2$ is not the minimal value of
    $\parallel \mu_{_{G}}\parallel^2$ on $M$. Then for every 
    $\mu_{_{G}}$-positive vector bundle $E$ over $M$ we have 
    the decomposition
    $RR_{\beta}^{^{G}}(M,E)=\sum_{\lambda}m_{\beta,\lambda}(E)\, 
    \chi_{_{\lambda}}^{_G}$ in irreducible characters with
    $$
    m_{\beta,\lambda}(E)\neq 0\Longrightarrow 
    \langle\lambda,\beta\rangle >0\ .
    $$
    In particular, if $\mu_{_{G}}^{-1}(0)$ is not empty, the
    $G$-invariant part of $RR_{\beta}^{^{G}}(M,E)$ is equal to $0$ for
    every $G$-invariant $\beta\in\Bcal_{_{G}}-\{0\}$.
    The result remains when $M$ is non-compact, and the moment map
    $\mu_{_{G}}$ is proper.
\end{lem}

{\em Proof :} Recall the localization formula 
on $M^{\beta}$ obtained in Proposition
\ref{prop.loc.RR.G.beta}. For every complex $G$-vector bundle $E$ 
over $M$, we have the 
following equality in $\widehat{R}(G)$
\begin{equation}\label{eq.non.compacte}
    RR_{\beta}^{^{G}}(M,E)=(-1)^{r_{\Ncal}}\sum_{k\in\N}
    RR_{\beta}^{^{G}}(M^{\beta},E\vert_{M^{\beta}}\otimes\det
    \Ncal^{+,\beta}\otimes S^k((\Ncal\otimes\C)^{+,\beta}) \ .
\end{equation}
Suppose that $M$ is non-compact and that the moment map 
$\mu_{_{G}}$ is proper as a map from a $G$-invariant open neighborhood
of $\mu_{_{G}}^{-1}(\beta)$ in $M$ to a $G$-invariant open neighborhood
of $\beta$ in $\ggot$. Each terms of (\ref{eq.non.compacte}) 
are well defined and the equality remains valid in this case (It is 
not difficult to extend the proof given in subsection \ref{subsec.RR.G.beta}
to this situation).

If $\parallel \beta\parallel^2$ is not the minimal value of
$\parallel \mu_{_{G}}\parallel^2$, we know from Lemma \ref{lem.N.beta.plus}, 
that the vector bundle $\Ncal^{+,\beta}$ is not trivial over 
each connected component $\Zcal$ of $M^{\beta}$ that intersects 
$\mu^{-1}(\beta)$. Then every $\tore_{\beta}$-weight $a$ 
on the fibers of the complex vector bundle 
$E\vert_{\Zcal}\otimes\det\Ncal^{+,\beta}\otimes S^k((\Ncal\otimes\C)^{+,\beta}$
satisfies $\langle a,\beta\rangle>0$.
Lemma \ref{lem.multiplicites.tore} and Corollary 
\ref{coro.multiplicites.tore}, applied to this situation, show that 
$RR_{\beta}^{^{G}}(M,E)=\sum_{\lambda}m_{\beta,\lambda}(E)\,
    \chi_{_{\lambda}}^{_G}$ with  $m_{\beta,\lambda}(E)\neq 0$
    only if $\langle\lambda,\beta\rangle >0$. $\Box$

\medskip

\subsection{The map $RR^{^{G}}_{\beta}$ when $G_{\beta}\neq G$.}

Let $\sigma$ be the unique open face of $\hgot_{+}$ which contains 
$\beta$. The stabilizer subgroup $G_{\xi}$ does not depend on the 
choice of $\xi\in\sigma$, and is denoted by $G_{\sigma}$. Let 
$\ggot_{\sigma}$ be the Lie algebra of $G_{\sigma}$, and let 
$U_{\sigma}$ the $G_{\sigma}$-invariant open subset of 
$\ggot_{\sigma}$ defined by $U_{\sigma}=G_{\sigma}\cdot
\{y\in \hgot_{+}\vert G_{y}\subset G_{\sigma}\}$.

The symplectic cross-section Theorem 
\cite{Guillemin-Sternberg84,L-M-T-W}
asserts that the pre-image $\Ycal_{\sigma}=\mu_{_{G}}^{-1}(U_{\sigma})$
is a symplectic submanifold of $M$ provided with a Hamiltonian action 
of $G_{\sigma}$. We denote by $\omega_{\sigma}$ the symplectic $2$-form on
$\Ycal_{\sigma}$, and $\mu_{\sigma}:\Ycal_{\sigma}\to
\ggot_{\sigma}$  the moment map. Let $J_{\sigma}$ be a 
$G_{\sigma}$-invariant almost complex structure on $\Ycal_{\sigma}$, 
which is compatible with $\omega_{\sigma}$. 
The vector field $\Hcal^{\sigma}$ 
on $\Ycal_{\sigma}$ generated by $\mu_{\sigma}$ vanishes on 
$C^{\sigma}_{\beta}:=
\mu^{-1}_{\sigma}(\beta)\cap (\Ycal_{\sigma})^{\beta}
=\mu^{-1}_{_G}(\beta)\cap M^{\beta}$ 
(see Definition \ref{def.H.vector}). We denote by\footnote{
For a non-compact $G$-manifold $\Xcal$, we denote by 
$\tilde{K}_{G}(\Xcal)$ the equivariant $K$-theory of $\Xcal$ with 
non-compact support.}
$$ 
RR^{^{G_{\sigma}}}_{\beta}(\Ycal_{\sigma},-):
\tilde{K}_{G_{\sigma}}(\Ycal_{\sigma})\to R^{-\infty}(G_{\sigma})
$$
the Riemann-Roch character on $\Ycal_{\sigma}$ localized near 
the {\em compact} subset $C^{\sigma}_{\beta}$ by the vector filed 
$\Hcal^{\sigma}$. It   
is well defined even since $\mu_{\sigma}$ is a proper map 
(see Definition \ref{def.RR.beta}).

\medskip

\begin{theo}\label{th.RR.G.beta.hamilton}
For every $E\in K_{G}(M)$, we have 
$$
RR^{^{G}}_{\beta}(M,E)=\Hols\left(
RR^{^{G_{\sigma}}}_{\beta}(\Ycal_{\sigma},E\vert_{\Ycal_{\sigma}})
\right)\quad {\rm in}\quad R^{-\infty}(G)\ , 
$$
\end{theo}

\medskip

\begin{coro} Let $\beta\in\Bcal_{_{G}}$ with $G_{\beta}\neq G$.
   If $\mu_{_{G}}^{-1}(0)\neq \emptyset$, we have 
$[ RR^{^{G}}_{\beta}(M,E)]^{G}=0$, for every 
$\mu_{_{G}}$-positive vector bundle $E\to M$. In general,
$[ RR^{^{G}}_{\beta}(M,E)]^{G}=0$, for every 
$\mu_{_{G}}$-{\em strictly} positive vector bundle $E$.

\end{coro}

\medskip

{\em Proof of the Corollary :} 
The moment map $\mu_{\sigma}$ is proper
as a map from a $G_{\sigma}$-invariant open neighborhood
of $\mu_{\sigma}^{-1}(\beta)$ in $\Ycal_{\sigma}$ to a 
$G_{\sigma}$-invariant open neighborhood of $\beta$ in $\ggot_{\sigma}$.
If $0\in \mu_{_{G}}(M)$ we see that $t\beta\in 
\mu_{\sigma}(\Ycal_{\sigma})$ for any $0<t<1$, hence 
$\parallel \beta\parallel^2$ is not the minimal value of
$\parallel \mu_{\sigma}\parallel^2$.

Proposition \ref{lem.RR.beta.G.invariant} can be 
used for the map $RR^{^{G_{\sigma}}}_{\beta}(\Ycal_{\sigma},-)$. 
For any $\mu_{_{G}}$-positive vector bundle $E$, we have
$RR^{^{G_{\sigma}}}_{\beta}(\Ycal_{\sigma},E\vert_{\Ycal_{\sigma}})=
\sum_{\lambda}m_{\beta,\lambda}(E)\,\chi_{_{\lambda}}^{_{G_{\sigma}}}$ 
with  $m_{\beta,\lambda}(E)\neq 0$ only if 
$\langle\lambda,\beta\rangle >0$ (the same holds when 
$0\notin \mu_{_{G}}(M)$ and $E$ is $\mu_{_{G}}$-{\em strictly} positive).
With the induction formula of Theorem \ref{th.RR.G.beta.hamilton} we 
get\footnote{$\Hols(\chi_{_{\lambda}}^{_{G_{\sigma}}})=\HolH(h^{\lambda})$ 
since $\chi_{_{\lambda}}^{_{G_{\sigma}}}=\HolHs(h^{\lambda})$.}
$RR^{^{G}}_{\beta}(M,E) =\sum_{\lambda}m_{\beta,\lambda}(E)\,
\HolH(h^{\lambda})$. But $\HolH(h^{\lambda})=\pm 1$ only if 
$\langle \lambda,X\rangle\leq 0$ for every $X$ in the Weyl chamber
(see Remark \ref{hol.egale.1}). 
This shows 
$$
\HolH(h^{\lambda})=\pm 1\Longrightarrow \langle 
\lambda,\beta\rangle\leq 0
\Longrightarrow m_{\beta,\lambda}(E)=0\ .
$$
We have then proved that $[RR^{^{G}}_{\beta}(M,E)]^{G}=0$. $\Box$

\bigskip


\underline{\bf  Proofs of Theorem \ref{th.RR.G.beta.hamilton} :}

We propose here two different proofs for this induction formula.  Both 
of them use the same technical remark. 

The set $G\cdot \Ycal_{\sigma}\cong 
G\times_{G_{\sigma}}\Ycal_{\sigma}$ is a $G$-invariant open 
neighborhood  of the critical set $C_{\beta}^{^G}$ in $M$.
The symplectic form $\omega$, when restricted to 
$G\times_{G_{\sigma}}\Ycal_{\sigma}$, 
can be written in terms of the moment map
$\mu_{\sigma}$ and the symplectic form $\omega_{\sigma}$:
\begin{equation}\label{eq.omega.slice.bis}
\omega_{[g,y]}(X+v,Y+w)=-(\mu_{\sigma}(y),[X,Y])+ 
\omega_{\sigma}\vert_{y}(v,w)\ , 
\end{equation}
where $X,Y\in \ggot/\ggot_{\beta}$, and $v,w\in\T_{y}
\Ycal_{\sigma}$\footnote{We use here the identification 
$\T(G\times_{G_{\sigma}}\Ycal_{\sigma})\cong G\times_{G_{\sigma}}
(\ggot/\ggot_{\sigma}\oplus \T\Ycal_{\sigma})$ 
(see (\ref{eq.espace.tangent})).}.
With the complex structure 
$J_{G/G_{\sigma}}$ on $G/G_{\sigma}$ determined by $\beta$, 
we form the almost complex structure $\widetilde{J}:= 
J_{G/G_{\sigma}}\times J_{\sigma}$ on 
$G\times_{G_{\sigma}}\Ycal_{\sigma}$. Equation 
(\ref{eq.omega.slice.bis}) shows that  $\widetilde{J}$ is compatible 
with $\omega$ in a neighborhood of $C_{\beta}^{^G}$, 
hence $\widetilde{J}$ is homotopic
to $J$ in a neighborhood of $C_{\beta}^{^G}$ in 
$G\times_{G_{\sigma}}\Ycal_{\sigma}$.

\begin{rem}\label{rem.J.tilde}
  The almost complex structures $J$ and $\widetilde{J}$ are 
  homotopic in a neighborhood of $C_{\beta}^{^G}$, so as in Lemma 
  \ref{lem.inv.homotopy} we see that the computation of the localized 
  Riemann-Roch character $RR^{^{G}}_{\beta}(M,E)$ can be done with 
  $\widetilde{J}$ instead of $J$.
\end{rem}

\medskip


{\bf First proof of Theorem \ref{th.RR.G.beta.hamilton} :} 
We will show here that 
Theorem \ref{th.RR.G.beta.hamilton} is a 
consequence of the induction formula proved in Theorem 
\ref{th.induction.G.H}
and of the localization formula obtained in Proposition \ref{prop.loc.RR.G.beta}.
The induction of Corollary \ref{coro.induction.G.G.beta} shows that
$RR^{^{G}}_{\beta}(M,E)=\Hols(RR^{^{G_{\sigma}}}_{\beta}(M,E)
\wedge^{\bullet}\overline{\ggot/\ggot_{\sigma}})$.
So we have to prove the following equality 
\begin{equation}\label{equation1.th.7.4}
RR^{^{G_{\sigma}}}_{\beta}(\Ycal_{\sigma},E\vert_{\Ycal_{\sigma}})=
RR^{^{G_{\sigma}}}_{\beta}(M,E)
\wedge^{\bullet}\overline{\ggot/\ggot_{\sigma}}\ .
\end{equation}
First we use the localization formula on both 
sides of the equality. For the map 
$RR^{^{G_{\sigma}}}_{\beta}(M,-)$ 
this gives
\begin{equation}\label{equation2.th.7.4}
RR^{^{G_{\sigma}}}_{\beta}(M,E)=
RR^{^{G_{\sigma}\times{\rm T}_{\beta}}}_{\beta}
\left(M^{\beta},E\vert_{M^{\beta}}\otimes
\left[\wedge_{\C}^{\bullet}\overline{\Ncal}\,
\right]^{-1}_{\beta}\right)\ ,
\end{equation}
and for $RR^{^{G_{\sigma}}}_{\beta}(\Ycal_{\sigma},-)$ we have
\begin{equation}\label{equation3.th.7.4}
RR^{^{G_{\sigma}}}_{\beta}(\Ycal_{\sigma},E\vert_{\Ycal_{\sigma}})=
RR^{^{G_{\sigma}\times{\rm T}_{\beta}}}_{\beta}
\left((\Ycal_{\sigma})^{\beta},E\vert_{(\Ycal_{\sigma})^{\beta}}
\otimes\left[\wedge_{\C}^{\bullet}\overline{\Ncal'}\,
\right]^{-1}_{\beta}
\right)\ .
\end{equation}

Here $\Ncal$ and $\Ncal'$ are respectively  the normal bundle of 
$M^{\beta}$ in $M$, and the normal bundle of 
$(\Ycal_{\sigma})^{\beta}$ in $\Ycal_{\sigma}$. The complex 
structures on the fibers of $\Ncal$ and $\Ncal'$ are  induced 
respectively by the almost complex 
structure $\widetilde{J}$, and by the almost complex structure 
$J_{\sigma}$ (see Remark \ref{rem.J.tilde}).

Now we remark that $(\Ycal_{\sigma})^{\beta}$ is an 
open neighborhood of $M^{\beta}\cap\mu_{_{G}}^{-1}(\beta)$ 
in $M^{\beta}$, thus we have 
$RR^{^{G_{\sigma}}}_{\beta}(M^{\beta},F)=
RR^{^{G_{\sigma}}}_{\beta}((\Ycal_{\sigma})^{\beta},
F\vert_{(\Ycal_{\sigma})^{\beta}})$ for any equivariant vector bundle
$F$. So 
(\ref{equation2.th.7.4}) and (\ref{equation3.th.7.4}) 
shows us that (\ref{equation1.th.7.4}) is equivalent
to the following
\begin{eqnarray}\label{equation4.th.7.4}
\lefteqn{RR^{^{G_{\sigma}\times{\rm T}_{\beta}}}_{\beta}
\left((\Ycal_{\sigma})^{\beta},E\vert_{(\Ycal_{\sigma})^{\beta}}
\otimes\left[\wedge_{\C}^{\bullet}\overline{\Ncal}\,
\right]^{-1}_{\beta}\otimes \left[\wedge^{\bullet}
\overline{\ggot/\ggot_{\sigma}}\,\right]\right)
=}\\
& & RR^{^{G_{\sigma}\times{\rm T}_{\beta}}}_{\beta}
\left((\Ycal_{\sigma})^{\beta},E\vert_{(\Ycal_{\sigma})^{\beta}}
\otimes\left[\wedge_{\C}^{\bullet}\overline{\Ncal'}\,
\right]^{-1}_{\beta}\right)\ ,\nonumber
\end{eqnarray}
where $\left[\wedge^{\bullet}\overline{\ggot/\ggot_{\sigma}}\,\right]$
is the trivial bundle 
$\wedge^{\bullet}\overline{\ggot/\ggot_{\sigma}}\times 
(\Ycal_{\sigma})^{\beta}\to (\Ycal_{\sigma})^{\beta}$. 

To finish the proof, we notice that the normal bundle $\Ncal\to 
M^{\beta}$, when restricted to $(\Ycal_{\sigma})^{\beta}$, can be 
decomposed as
$\Ncal\vert_{(\Ycal_{\sigma})^{\beta}}=\Ncal'
\oplus[\ggot/\ggot_{\sigma}]$. 
Here $[\ggot/\ggot_{\sigma}]\to (\Ycal_{\sigma})^{\beta}$ is 
the trivial complex vector bundle defined  by 
$[\ggot/\ggot_{\sigma}]_{m}=\{ X_{(\Ycal_{\sigma})^{\beta}}\vert_{m},\ 
X\in \ggot/\ggot_{\sigma}\}$ for any $m\in (\Ycal_{\sigma})^{\beta}$.
This decomposition gives first  
the equality $\wedge_{\C}^{\bullet}\overline{\Ncal}=
\wedge_{\C}^{\bullet}\overline{\Ncal'}\otimes
[\wedge_{\C}^{\bullet}\overline{\ggot/\ggot_{\sigma}}\,]$ 
and after\footnote{The product of 
$\left[\wedge_{\C}^{\bullet}\overline{\Ncal'}\,\right]^{-1}_{\beta}$ 
and $\left[\wedge_{\C}^{\bullet}
\overline{\ggot/\ggot_{\sigma}}\,\right]^{-1}_{\beta}$
is well defined in $\tilde{K}_{G_{\sigma}}((\Ycal_{\sigma})^{\beta})
\widehat{\otimes}\,R(\tore_{\beta})$ since these elements are
polarized by $\beta$: each of them is a sum over the set of 
weights of $\tore_{\beta}$ of the form 
$\sum_{\alpha}E_{\alpha}h^{\alpha}$ such that 
$E_{\alpha}\neq 0$ only if $\langle\alpha,\beta\rangle\geq 0$, 
and for any $\delta'>\delta\geq 0$ the sum 
$\sum_{\delta\leq\langle\alpha,\beta\rangle\leq\delta'}
E_{\alpha}h^{\alpha}$ is 
finite (see definition \ref{wedge.V.inverse}).}
$\left[\wedge_{\C}^{\bullet}\overline{\Ncal}\,\right]^{-1}_{\beta}=
\left[\wedge_{\C}^{\bullet}\overline{\Ncal'}\,\right]^{-1}_{\beta}
\otimes\left[\wedge_{\C}^{\bullet}\overline{\ggot/\ggot_{\sigma}}
\,\right]^{-1}_{\beta}\ ,
$
which implies
$
\left[\wedge_{\C}^{\bullet}\overline{\Ncal}\,\right]^{-1}_{\beta}
\otimes
\left[\wedge_{\C}^{\bullet}\overline{\ggot/\ggot_{\sigma}}\,\right]=
\left[\wedge_{\C}^{\bullet}\overline{\Ncal'}\,\right]^{-1}_{\beta}\ .
$
(\ref{equation4.th.7.4}) is then proved. $\Box$

\medskip


{\bf  Second proof of Theorem \ref{th.RR.G.beta.hamilton} :}
A $G$-invariant neighborhood $\Ucal^{^{G,\beta}}$ of the critical set 
$C_{\beta}^{^G}$ in $M$ can be taken of the form $\Ucal^{^{G,\beta}}=
G\times_{G_{\sigma}}\Ucal^{^{\sigma,\beta}}$ where 
$\Ucal^{^{\sigma,\beta}}$ a relatively compact $G_{\sigma}$-invariant 
neighborhood of $\mu^{-1}_{_G}(\beta)\cap M^{\beta}$ in 
$\Ycal_{\sigma}$ such that 
$\overline{\Ucal^{^{\sigma,\beta}}}\cap\{\Hcal^{\sigma}=0\}= 
\mu^{-1}_{_G}(\beta)\cap M^{\beta}$.

The maps $RR^{^{G}}_{\beta}(M, -)$ and 
$RR^{^{G_{\sigma}}}_{\beta}(\Ycal_{\sigma},-)$ are respectively 
defined by the localized Thom symbols 
$\Thom_{G,[\beta]}^{\mu}(M)\in K_{G}(\T_{G}\Ucal^{^{G,\beta}})$ and  
$\Thom_{G_{\sigma},[\beta]}^{\mu}(\Ycal_{\sigma})
\in K_{G_{\sigma}}(\T_{G_{\sigma}}\Ucal^{^{\sigma,\beta}})$ 
(see Definition \ref{def.thom.beta.f}). The inclusion 
$i:G_{\sigma}\croc G$ induces an isomorphism
$i_{*}:K_{G_{\sigma}}(\T_{G_{\sigma}}\Ucal^{^{\sigma,\beta}})\to
K_{G}(\T_{G}(G\times_{G_{\sigma}}\Ucal^{^{\sigma,\beta}}))$ 
(see subsection \ref{subsec.induction.def}). 

\begin{lem}\label{lem.induction.slice}
We have the following equality
$$
i_{*}\left(\Thom_{G_{\sigma},[\beta]}^{\mu}(\Ycal_{\sigma})
\, \wedge^{\bullet}_{\C}\ggot/\ggot_{\sigma}\right)
=\Thom_{G,[\beta]}^{\mu}(M)\ .
$$
\end{lem}

\medskip

This Lemma, combined with Theorem \ref{thm.atiyah.2}, shows that 
$RR^{^{G}}_{\beta}(M, E)=$\break
$\inds\left(RR^{^{G_{\sigma}}}_{\beta}(\Ycal_{\sigma},E\vert_{\Ycal_{\sigma}})
\,\wedge^{\bullet}_{\C}\ggot/\ggot_{\sigma}\right)=
\Hols\left(RR^{^{G_{\sigma}}}_{\beta}(\Ycal_{\sigma},
E\vert_{\Ycal_{\sigma}})\right)$ 
for any $G$-complex vector bundle $E\to M$. The proof of Theorem 
\ref{th.RR.G.beta.hamilton} is then completed. $\Box$

\medskip

{\em Proof of Lemma \ref{lem.induction.slice} :}

\medskip

Through the identification 
$\T(G\times_{G_{\sigma}}\Ucal^{^{\sigma,\beta}}) \cong 
G\times_{G_{\sigma}}(\ggot/\ggot_{\sigma}\oplus 
\T\Ucal^{^{\sigma,\beta}})$, the vector fiels $\Hcal^{\sigma}$ and 
$\Hcal^{^{G}}$ satisfy the relation $\Hcal^{^{G}}_{[g,y]}\cong 
\Hcal^{\sigma}_{y},\ [g,y]\in \Ucal^{^{G,\beta}}$.  The symbol 
$\sigma_{[g,y;X+v]}$ of $\Thom_{G,[\beta]}^{\mu}(M)$ at $[g,y;X+v]\in 
G\times_{G_{\sigma}}(\ggot/\ggot_{\sigma}\oplus 
\T\Ucal^{^{\sigma,\beta}})$ acts on 
$\wedge^{\bullet}_{\widetilde{J}}\T_{[g,y]}\Ucal^{^{G,\beta}} \cong 
\wedge^{\bullet}\ggot/\ggot_{\sigma}\otimes 
\wedge^{\bullet}_{J_{\sigma}}\T_{y}\Ucal^{^{\sigma,\beta}}$ as the 
product
$$
\sigma_{[g,y;X+v]}=Cl(X)\odot Cl_{y}(v-\Hcal^{\sigma}_{y})\ .
$$
Now we see that $[g,y;X+v]\to Cl(X)\odot Cl_{y}(v-\Hcal^{\sigma}_{y})$ 
is homotopic, as $G$-transversally elliptic symbol, to 
$\widetilde{\sigma}: [g,y;X+v]\to Cl(0)\odot 
Cl_{y}(v-\Hcal^{\sigma}_{y})$, and $\widetilde{\sigma}$ is, by 
definition, the image of 
$\Thom_{G_{\sigma},[\beta]}^{\mu}(\Ycal_{\sigma}) 
\wedge^{\bullet}_{\C}\ggot/\ggot_{\sigma}$ by $i_{*}$.  The proof of 
Lemma \ref{lem.induction.slice} is then completed.  $\Box$

\medskip 

\subsection{The singular case.}\label{non-regulier}

In this section, we do not assume that $0$ is a regular value of
$\mu_{_G}$, and we use the `shifting trick' to
compute $[RR^{^G}(M,L)]^G$ in term of reduced manifolds 
of the type $\mu_{_G}^{-1}(a)/G_{a}$, for every $\mu_{_G}$-moment 
bundle $L$. We know from Theorem \ref{QR=RQ-regulier} that $[RR^{^G}(M,L)]^G=0$ 
if $0\notin \mu_{_G}(M)$ since every moment bundle is
strictly positive (see Lemma \ref{L-a-positif}). 
So, we assume for the rest of this section that $0\in \mu_{_G}(M)$.

Let $\Ocal_a$ be the coadjoint orbit through $a\in\ggot^*
$. It has a canonical symplectic 2-form and the moment map $\Ocal_a\to\ggot^*$
for the $G$-action is the inclusion. We denote by $\overline{\Ocal_a}$
the coadjoint orbit $\Ocal_a$ with the opposite symplectic form. The
product $M\times\overline{\Ocal_a}$ is a symplectic manifold with a
Hamiltonian moment map
\begin{eqnarray*}
    \mu_a :M\times\overline{\Ocal_a}&\longrightarrow&\ggot^*\\ 
(m,\xi)&\longmapsto& \mu_{_G}(m)-\xi\ .
\end{eqnarray*}

On the symplectic manifold $M\times\overline{\Ocal_a}$ we have a
quantization map $RR^{^G}(M\times\overline{\Ocal_a},-)$ with the
following property: for any $G$-vector bundles $E$ and $F$ over
$M$ and $\Ocal_a$ respectively, we have
$RR^{^G}(M\times\overline{\Ocal_a},\pi_a^*(E)\otimes(\pi_a')^*(F))=
RR^{^G}(M,E)\cdot RR^{^G}(\overline{\Ocal_a},F)$ in $R(G)$. Here
we denote by $\pi_a:M\times\overline{\Ocal_a}\to M$ the projection to 
the first factor and $\pi_a'$ the projection to the second factor.
Since $RR^{^G}(\overline{\Ocal_a},\C)=1$ we have
\begin{equation}\label{shift-a}
RR^{^G}(M\times\overline{\Ocal_a},\pi_a^*(L))= RR^{^G}(M,L)\ .
\end{equation}

We can now compute $[RR^{^G}(M,L)]^G$ by  localizing the character
$RR^{^G}(M\times\overline{\Ocal_a},\pi_a^*(L))$ with  the
moment map $\mu_a$. We need the following Lemma which was proved 
by Tian-Zhang \cite{Tian-Zhang97} for the prequantum line bundles.
\begin{lem}\label{L-a-positif}
Let $L$ be a $\mu_{_G}$-moment bundle over $M$.  There exists 
$\epsilon>0$ such that for any $\vert a\vert<\epsilon$, the vector 
bundle $\pi_a^*(L)$ is $\mu_a$-positive.  For $a=0$, the bundle 
$L=\pi_0^*(L)$ is $\mu_{_{G}}$-strictly positive.
\end{lem}

Let $RR^{^G}_0(M\times\overline{\Ocal_a},-)$ be the Riemann-Roch
character localized near $\mu_a^{-1}(0)\simeq \mu_{_G}^{-1}(\Ocal_a)$.
Theorem \ref{QR=RQ-regulier}, Equality \ref{shift-a}, and Lemma \ref{L-a-positif}
show that 
\begin{equation}\label{shift-a-0}
[RR^{^G}(M,L)]^G=[RR^{^G}_0(M\times\overline{\Ocal_a},\pi_a^*(L))]^G\ ,
\end{equation}
for any moment bundle $L$ if $a\in\mu_{_G}(M)$ is close enough to $0$.

There exists a unique open face $\tau$ of the Weyl chamber $\hgot_+$ such that 
$\mu_{_G}(M)\cap\tau$ is dense in $\mu_{_G}(M)\cap\hgot_+$.
The face $\tau$ is called the principal face of $(M,\mu_{_G})$ \cite{L-M-T-W}.
All points in the open face $\tau$ have the same connected 
centralizer $G_{\tau}$. Let $A_{\tau}$ be the identity component of 
the center of $G_{\tau}$ and $[G_{\tau},G_{\tau}]$ its semisimple part. Note 
that we have an identification between 
the Lie algebra $\agot_{\tau}$  of $A_{\tau}$ and the 
linear span of the face
$\tau$. The Principal-cross-section Theorem \cite{L-M-T-W} tells us that
$Y_{\tau}:=\mu_{_G}^{-1}(\tau)$ is a symplectic $G_{\tau}$-manifold, 
with a trivial action of $[G_{\tau},G_{\tau}]$. So, the restriction
of $\mu_{_G}$ on $\Ycal_{\tau}$ is a moment map 
$\mu_{\tau}: \Ycal_{\tau}\to\agot_{\tau}$ for the Hamiltonian action of 
the torus $A_{\tau}$. We decompose the torus $A_{\tau}$ 
in a product of two subtorus $A_{\tau}=A_{\tau}^1\times A_{\tau}^2$ 
where $A_{\tau}^{1}$ is the identity component of the principal 
stabilizer for the action of $A_{\tau}$ on $\Ycal_{\tau}$. 

We take now $a$ with value in $\tau\cap\mu_{_G}(M)$. 
For generic values $a\in \tau\cap\mu_{_G}(M)$, 
$\mu_{_G}^{-1}(a)=\mu_{\tau}^{-1}(a)$ is
a smooth manifold of $M$  with a locally free action of $A_{\tau}^{2}$, hence
the quotient $\Mcal_a:=\mu_{_G}^{-1}(a)/G_a=
\mu_{\tau}^{-1}(a)/(A_{\tau}^2)$ 
is a symplectic orbifold. We denote by $RR(\Mcal_{a},-)$ the quantization
map defined by the choice of a compatible almost complex structure. 
If $L$ is a $\mu_{_{G}}$-moment bundle on $M$, $L\vert_{\Ycal_{\tau}}$
is a $\mu_{\tau}$-moment bundle: the action of 
$A_{\tau}^{1}[G_{\tau},G_{\tau}]$ on $L\vert_{\Ycal_{\tau}}$ is trivial.
Then the quotient $L\vert_{\mu_{\tau}^{-1}(a)}/G_{a}=
L\vert_{\mu_{\tau}^{-1}(a)}/(A_{\tau}^2)$ is an orbifold line bundle
over $\Mcal_a$ for generic $a$.
 
We compare now the  Riemann-Roch character $RR^{^{G_{\tau}}}_0(\Ycal_{\tau},-)$ 
localized near $\mu_{\tau}^{-1}(a)$ by the moment map $\mu_{\tau}-a$ 
and the  Riemann-Roch character $RR^{^G}_0(M\times\overline{\Ocal_a},-)$
localized near $\mu_{a}^{-1}(0)=G\cdot(\mu_{\tau}^{-1}(a)\times\{a\})$.
All we need is contained in the following 
\begin{prop}\label{RR-a}
Let $E$ be a G-vector bundle over $M$, and take $a\in\tau$. We have
$RR^{^G}_0(M\times\overline{\Ocal_a},\pi_{a}^{*}E)=
\indT\left(RR^{^{G_{\tau}}}_0(\Ycal_{\tau},E\vert_{\Ycal_{\tau}})\right)$,
in particular 
$[RR^{^G}_0(M\times\overline{\Ocal_a},\pi_{a}^{*}E)]^G=
[RR^{^{G_{\tau}}}_0(\Ycal_{\tau},E\vert_{\Ycal_{\tau}})]^{G_{\tau}}$.
\end{prop}

If $L$ is a $\mu_{_{G}}$-moment bundle, the action of 
$A_{\tau}^{1}[G_{\tau},G_{\tau}]$ on $L\vert_{\Ycal_{\tau}}$ is trivial, then 
$[RR^{^{G_{\tau}}}_0(\Ycal_{\tau},L\vert_{\Ycal_{\tau}})]^{G_{\tau}}=
[RR^{^{A^2_{\tau}}}_0(\Ycal_{\tau},L\vert_{\Ycal_{\tau}})]^{A^{2}_{\tau}}$. 
Finally, for every generic 
value $a\in\tau\cap\mu_{_G}(M)$, the quotient  
$L_{a}:=L\vert_{\mu_{\tau}^{-1}(a)}/A^{2}_{\tau}$ is an orbifold line 
bundle over $\Mcal_{a}$, so from subsection \ref{subsec.RR.O.hamilton} 
we get $[RR^{^{A^{2}_{\tau}}}_0(\Ycal_{\tau},L\vert_{\Ycal_{\tau}})]^{A^{2}_{\tau}}=
RR(\Mcal_{a},L_{a})$. 

With this last equality, Proposition \ref{L-a-positif}, and
equality (\ref{shift-a-0}) we have proved the central result of 
this section
\begin{prop}\label{RR-0-singulier}
  Suppose that $0\in\mu_{_{G}}(M)$. If $L$ is a $\mu_{_{G}}$-moment bundle, 
  there exist $\epsilon>0$, such that
  $$
  [RR^{^{G}}(M,L)]^{G}= RR(\Mcal_{a},L_{a})\ ,
  $$
  for every generic value $a\in \tau\cap\mu_{_{G}}(M)$ with $\vert 
  a\vert<\epsilon$.
\end{prop}

\subsubsection{Proof of Lemma \ref{L-a-positif}}

Let $L$ be a $\mu_{_{G}}$-moment bundle over $M$, where
$\mu_{_{G}}:M\to\ggot^{*}$ is a Hamiltonian moment map. 
Recall that the Lie 
algebra $\ggot$ is identified to $\ggot^{*}$ trough
an invariant scalar product $(-,-)$. Let $H$ be a maximal 
torus of $G$ with Lie algebra $\hgot$.

\begin{lem}\label{critere-L-positif} For $\beta\in \hgot$ and 
    $m\in M^{\beta}\cap\mu_{_{G}}^{-1}(\gamma)$, 
    the weight $\alpha$ for the action of $\tore_{\beta}$ on $L_{m}$
    satifies $(\alpha,\beta)=(\gamma,\beta)$.
\end{lem}    
    
{\em Proof }: Let $N$ be the connected component of $M^{\beta}$ 
containing $m$, and let $m'$ be a point of $N^{H}$.  Since $N$ is 
connected, $\alpha$ is also the weight for the action of 
$\tore_{\beta}$ on $L_{m'}$, and $\mu_{_{G}}(m')$ is the 
weight for the action of $H$ on $L_{m'}$: then 
$(\alpha,X)=(\mu_{_{G}}(m'),X)$ for every $X\in 
Lie(\tore_{\beta})$.  But the map $x\to(\mu_{_{G}}(x),\beta)$ 
is constant on $N$, then 
$(\gamma,\beta)=(\mu_{_{G}}(m),\beta)= 
(\mu_{_{G}}(m'),\beta)=(\alpha,\beta)$.  $\Box$

\medskip
    
The element $a$ is taken in $\hgot$.  The critical set of the function
$||\mu_{a}||^2:M\times\Ocal_{a}\to\R$ admits the following
decomposition $\Cr(||\mu_{_a}||^2)= G\cdot(\Cr(||\mu_{_a}||^2)\cap
(M\times \{a\}))= G\cdot
\Big(\left(\Cr(||\mu_{_{G_{a}}}-a||^2)\cap\mu_{_{G}}^{-1}(\ggot_{a})\right)
\times \{a\}\Big)$,
where $\mu_{_{G_{a}}}:M\to\ggot_{a}$ is the moment map for the action
of $G_{a}$.  Let $\Bcal_{a}$ the finite subset of $\hgot$ defined
by $\Bcal_{a}=\{\beta\in\hgot,\ M^{\beta}\cap\mu_{_{G}}^{-1}(\beta
+a)\neq \emptyset\}$.  Finally we have the decomposition

$$
\Cr(||\mu_{_a}||^2)=\bigcup_{\beta\in\Bcal_{a}}
G\cdot\left(M^{\beta}\cap\mu_{_{G}}^{-1}(\beta 
+a)\times\{a\}\right)\ .
$$

Using Lemma \ref{critere-L-positif}, we see that $\pi_{a}^{*}L$ is 
$\mu_{a}$-positive if and only if 
\begin{equation}\label{L-a-bis}
    (\beta + a,\beta)\geq 0\quad {\rm for\  every}\quad
\beta\in\Bcal_{a}\ . 
\end{equation}
We first see that it is trivially true
if $a=0$: in this case $L$ is strictly positive.

Let $\mu_{_{H}}:M\to\hgot$ be the moment map for the maximal torus $H$. 
Consider the finite set $\Bcal_{H,a}$ which parametrizes the critical 
set of $||\mu_{_H}-a||^2$: $\Bcal_{H,a}=
\{\beta\in\hgot,\ M^{\beta}\cap\mu_{_{H}}^{-1}(\beta +a)\neq
\emptyset\}$. We have obviously the inclusion $\Bcal_{a}\subset
\Bcal_{H,a}$, so it suffices to show (\ref{L-a-bis}) for 
$\Bcal_{H,a}$. 

To finish our proof we use now a characterisation of the set
$\Bcal_{H,a}$ we gave in \cite{pep1}. There exists a finite collection 
$\Bcal$ of affine subspaces of $\hgot$ such that 
$$
\Bcal_{H,a}\subset\{P_{\Delta}(a)-a,\Delta\in \Bcal\}
$$
for every $a\in\hgot$. Here $P_{\Delta}:\hgot\to\hgot$ is the
orthogonal projection on $\Delta$. It is now easy to 
compute the sign of $(P_{\Delta}(a),P_{\Delta}(a)-a)$ for
all $\Delta\in\Bcal$. A simple computation gives
$(P_{\Delta}(a),P_{\Delta}(a)-a)=\vert P_{\Delta}(0)\vert^{2}
-(a,P_{\Delta}(0))$. Hence, either $0\in \Delta$ and then 
$(P_{\Delta}(a),P_{\Delta}(a)-a)$ is equal to $0$ for all $a\in\hgot$
 , or $0\notin \Delta$ and then $(P_{\Delta}(a),P_{\Delta}(a)-a)>0$ if
$\vert a\vert<\vert P_{\Delta}(0)\vert$. We can take $\epsilon=
\inf_{0\notin \Delta}\vert P_{\Delta}(0)\vert$ in 
Lemma \ref{L-a-positif}. $\Box$

\subsubsection{Proof of Proposition \ref{RR-a}}

Since the point $a$ takes value in $\tau$ we identify the 
coadjoint orbit $\Ocal_{a}$ with $G/G_{\tau}$. Let $\Hcal^{a}$ be the 
Hamiltonian vector field of the function 
$\frac{1}{2}\parallel\mu_{a}\parallel^{2}:M\times G/G_{\tau}\to\R$. 
To simplify the notations, 
$\Ycal_{\tau}$ will denote a small neighborhood 
of $\mu^{-1}_{_{G}}(a)$ in the symplectic slice 
$\mu^{-1}_{_{G}}(\tau)$ such that 
the open subset $\Ucal:=(G\times_{G_{\tau}}\Ycal_{\tau})\times 
G/G_{\tau}$ is then a neighborhood of 
$\mu_{a}^{-1}(0)=G\cdot(\mu_{\tau}^{-1}(a)\times\{\bar{e}\})$ which 
satisfies $\overline{\Ucal}\cap\{ \Hcal^{a}=0\}=\mu_{a}^{-1}(0)$.
Following Definition \ref{def.thom.beta.f}, the localized Riemann-Roch 
character $RR^{^G}_{0}(M\times G/G_{\tau},-)$ is computed by means of 
the Thom class $\Thom_{G,[0]}^{\mu_{a}}(M\times G/G_{\tau})\in 
K_{G}(\T_{G}\Ucal)$. On the other hand, the localized Riemann-Roch 
character $RR^{^{G_{\tau}}}_{0}(\Ycal_{\tau},-)$ is computed by means of 
the Thom class $\Thom_{G_{\tau},[0]}^{\mu_{\tau}-a}(\Ycal_{\tau})\in 
K_{G_{\tau}}(\T_{G_{\tau}}\Ycal_{\tau})$. 

Proposition \ref{RR-a} will follow from a simple relation 
between $\Thom_{G,[0]}^{\mu_{a}}(M\times G/G_{\tau})$ and 
$\Thom_{G_{\tau},[0]}^{\mu_{\tau}-a}(\Ycal_{\tau})$.

First, one considers  the isomorphism  
\begin{eqnarray}\label{def-U-prime}
    \phi:\Ucal&\to&\Ucal'\\
    ([g;y],[h])&\to&[g;[g^{-1}h],y]\ ,\nonumber
\end{eqnarray}    
with $\Ucal':=G\times_{G_{\tau}}(G/G_{\tau}\times\Ycal_{\tau})$, and
let $\phi^{*}:K_{G}(\T_{G}\Ucal')\to K_{G}(\T_{G}\Ucal)$ be
the induced isomorphism. After one consider the inclusion 
$i:G_{\tau}\croc G$ which induces an isomorphism 
$i_{*}:K_{G_{\tau}}(\T_{G_{\tau}}(G/G_{\tau}\times\Ycal_{\tau}))\to
K_{G}(\T_{G}\Ucal')$ (see subsection \ref{subsec.induction.def}). Let 
$j:\Ycal_{\tau}\croc G/G_{\tau}\times\Ycal_{\tau}$ be the 
$G_{\tau}$-invariant inclusion map defined by $j(y):=(\bar{e},y)$. We 
have then a pushforward map $j_{!}:
K_{G_{\tau}}(\T_{G_{\tau}}\Ycal_{\tau})\to 
K_{G_{\tau}}(\T_{G_{\tau}}(G/G_{\tau}\times\Ycal_{\tau}))$. 
Finally we have produced a map $\Theta:=\phi^{*}\circ i_{*}\circ 
j_{!}$ from $K_{G_{\tau}}(\T_{G_{\tau}}\Ycal_{\tau})$ to 
$K_{G}(\T_{G}\Ucal)$, such that $\indice_{\Ucal}^{G}(\Theta(\sigma))=
\indT(\indice_{\Ycal_{\tau}}^{G_{\tau}}(\sigma))$ for every 
$\sigma\in K_{G_{\tau}}(\T_{G_{\tau}}\Ycal_{\tau})$.

Proposition \ref{RR-a} is an immediate consequence of the following 
\begin{lem}
    We have the equality 
    $$
 \Theta\left(\Thom_{G_{\tau},[0]}^{\mu_{\tau}-a}(\Ycal_{\tau})\right)=
 \Thom_{G,[0]}^{\mu_{a}}(M\times G/G_{\tau})\ .
    $$
 \end{lem} 
 
{\em Proof }: Let $\mu_{a}':=\mu_{a}\circ\phi^{-1}$ be the moment map 
on $\Ucal'$, and let $\Hcal^{',a}$ be the Hamiltonian vector field 
of $\parallel\mu_{a}'\parallel$. For the tangent manifold $\T\Ucal'$ 
we have the decomposition 
$$
\T\Ucal'\simeq G\times_{G_{\tau}}\left(\ggot/\ggot_{\tau}\oplus
G\times_{G_{\tau}}(\overline{\ggot/\ggot_{\tau}})\oplus 
\T\Ycal_{\tau}\right)\ .
$$ 
A small computation gives $\Hcal^{',a}(m)=pr_{\ggot/\ggot_{\tau}}(ha)+ 
R(m) + \Hcal^{\tau}_{a}(y)+ S(m)$
for $m=[g;y,[h]]\in \Ucal'$, where $R(m)\in 
\overline{\ggot/\ggot_{\tau}}$ and  $S(m)\in \T_{y}\Ycal_{\tau}$ vanishes
when $m\in G\times_{G_{\tau}}(\{\bar{e}\}\times\Ycal_{\tau})$, i.e.
$[h]=\bar{e}$. Here $\Hcal^{\tau}_{a}$ is the Hamiltonian vector field 
of the function $\frac{1}{2}\parallel\mu_{\tau}-a\parallel^{2}:
\Ycal_{\tau}\to\R$. 

The transversally elliptic symbol $\sigma_{1}:=
(\phi^{-1})^{*}(\Thom_{G,[0]}^{\mu_{a}}(M\times G/G_{\tau}))$ 
is equal to the exterior product
$$
\sigma_{1}(m, \xi_{1}+\xi_{2}+ v)=
Cl(\xi_{1}-pr_{\ggot/\ggot_{\tau}}(ha))\odot
Cl(\xi_{2}-R(m))\odot
Cl(v-\Hcal^{\tau}_{a}- S(m))\ ,
$$
with $\xi_{1}\in\ggot/\ggot_{\tau}$,   
$\xi_{2}\in\overline{\ggot/\ggot_{\tau}}$,  $v\in\T\Ycal_{\tau}$.

Now we simplify the symbol $\sigma_{1}$ whithout changing its 
$K$-theoretic class. Since $\Char(\sigma_{1})\cap\T_{G}\Ucal'
= G\times_{G_{\tau}}(\{\bar{e}\}\times\Ycal_{\tau})$, we can transform
$\sigma_{1}$ through the 
$G_{\tau}$-invariant diffeomorphism $h=e^{X}$ from a neighborhood 
of  $0$ in $\ggot/\ggot_{\tau}$ to a neighborhood of $\bar{e}$ in 
$G/G_{\tau}$. This gives $\sigma_{2}\in 
K_{G}(\T_{G}(G\times_{G_{\tau}}(\ggot/\ggot_{\tau}\times 
\Ycal_{\tau})))$ defined by 
\begin{eqnarray*}
\lefteqn{\sigma_{2}([g,X,y], \xi_{1}+\xi_{2}+ v)
=}\\
& & 
Cl(\xi_{1}-pr_{\ggot/\ggot_{\tau}}(e^{X}a))\odot
Cl(\xi_{2}-R(m))\odot
Cl(v-\Hcal^{\tau}_{a}- S(m))\ .
\end{eqnarray*}    
Now trivial homotopies link $\sigma_{2}$ 
with the symbol $\sigma_{3}$, where we have removed the terms
$R(m)$ and $S(m)$, and where we have replaced 
$pr_{\ggot/\ggot_{\tau}}(e^{X}a)=[X,a]+o([X,a])$ by the term
$[X,a]$:
$$
\sigma_{3}([g,X,y], \xi_{1}+\xi_{2}+ v)= 
Cl(\xi_{1}-[X,a])\odot
Cl(\xi_{2})\odot
Cl(v-\Hcal^{\tau}_{a})\ .
$$
Now, we get $\sigma_{3}=i_{*}( \sigma_{4})$ where the symbol 
$\sigma_{4}\in K_{G_{\tau}}(\T_{G_{\tau}}(\ggot/\ggot_{\tau}\times 
\Ycal_{\tau}))$ is defined by 
$$
\sigma_{4}(X,y;\xi_{2}+ v)= 
Cl(-[X,a])\odot Cl(\xi_{2})\odot Cl(v-\Hcal^{\tau}_{a})\ .
$$    
So $\sigma_{4}$ is equal to the exterior product of $(y,v)\to
Cl(v-\Hcal^{\tau}_{a})$, which is 
$\Thom_{G_{\tau},[0]}^{\mu_{\tau}-a}(\Ycal_{\tau})$, with the 
transversally elliptic symbol on $\ggot/\ggot_{\tau}$:
$(X,\xi_{2})\to Cl(-[X,a])\odot Cl(\xi_{2})$. As in 
Lemma \ref{lem.indice.wedge.V.inverse}, we see that the 
$K$-theoretic class of this former symbol is  equal to $k_{!}(\C)$ 
where $k:\{0\}\croc\ggot/\ggot_{\tau}$. This shows that  
$$
\sigma_{4}=k_{!}(\C)\odot
\Thom_{G_{\tau},[0]}^{\mu_{\tau}-a}(\Ycal_{\tau})=
j_{!}(\Thom_{G_{\tau},[0]}^{\mu_{\tau}-a}(\Ycal_{\tau}))\ .
$$
$\Box$

\section{Appendix A: G=SU(2)}

We restrict our attention to an action of $G=SU(2)$ on a compact
manifold $M$.  We suppose that $M$ is endowed with a $G$-invariant
almost complex structure $J$ and an abstract moment map $f:M\to\ggot$. 
In this situation, the decomposition $RR^{^{G}}(M,-)=
\sum_{\beta\in\Bcal_{_{G}}}RR_{\beta}^{^{G}}(M,-)$ becomes simple.

Let $S^{1}$ be the maximal torus of $SU(2)$, and $f_{_{S^{1}}}:M\to
\R$ the induced moment map for the $S^{1}$-action.  The critical set
$\{\Hcal^{^G}=0\}$ has a particularly simple expression: $
\{\Hcal^{^G}=0\}=f^{-1}(0)\cup G.M^{S^{1}}_+$, 
where $M^{S^1}_{+}$ is the union of the connected components 
$F\subset M^{S^1}$ with $f_{S^1}(F)>0$. Note that the critical set 
$\{\Hcal^{^{S^1}}=0\}$ is equal to 
$f_{_{S^{1}}}^{-1}(0)\cup M^{S^{1}}$,

\medskip

\underline{The non-symplectic case}

\medskip

Here the induction formula of Theorem \ref{th.induction.G.H}, and 
Proposition \ref{prop.loc.RR.G.beta} gives
\begin{equation}\label{exemple1}
RR^{^{G}}(M,E)=RR^{^{G}}_{0}(M,E)+ 
    \HolS\Big(\Theta(E)(t).(1-t^{-2})\Big)
\end{equation}
where $\Theta(E)\in R^{-\infty}(S^1)$ is determined by 
\begin{equation}\label{exemple1.1}
\Theta(E)=  (-1)^{r_{\Ncal}}\sum_{k\in\N}
    RR^{^{S^1}}(M^{S^1}_{+},E\vert_{M^{S^1}_{+}}\otimes\det
    \Ncal^{+}\otimes S^k((\Ncal\otimes\C)^{+}))\ . 
\end{equation}
Here $\Ncal\to M^{S^1}_{+}$ is the normal bundle of 
$M^{S^1}_{+}$ in $M$.

\medskip

\underline{The Hamiltonian case}

\medskip

Here we suppose that $(M,\omega)$ is a symplectic manifold, with
moment map $\mu$ and a $\omega$-compatible almost complex structure $J$.
Let $\Ycal=\mu^{-1}(\R_{>0})$ be the symplectic slice associated 
to the interior of the Weyl chamber $\R_{>0}\subset Lie(S^1)$.

The induction formula of Theorem \ref{th.RR.G.beta.hamilton} 
gives
\begin{equation}\label{exemple2}
RR^{^{G}}(M,E)=RR^{^{G}}_{0}(M,E)+ 
    \HolS\Big(\widetilde{\Theta}(E)\Big)
\end{equation}
where $\widetilde{\Theta}(E)\in R^{-\infty}(S^1)$ is determined by   
\begin{equation}\label{exemple2.2}
\widetilde{\Theta}(E)=  (-1)^{r_{\widetilde{\Ncal}}}\sum_{k\in\N}
    RR^{^{S^1}}(M^{S^1}_{+},E\vert_{M^{S^1}_{+}}\otimes\det
    \widetilde{\Ncal}^{+}\otimes S^k((\widetilde{\Ncal}\otimes\C)^{+}))\ .
\end{equation}
Here $\widetilde{\Ncal}\to M^{S^1}_{+}$ is the normal bundle  of
$M^{S^1}_{+}$ in $\Ycal$.

\medskip

Recall that the irreducible characters $\phi_{n}$ of $G=SU(2)$ are 
labelled by $\Z_{\geq 0}$, and are completely determined by the 
relation $\phi_{n}=\HolS(t^n)$ in $R(G)$
(See Lemma \ref{lem.Weyl}). Hence the component 
$\HolS\Big(\Theta(E)(t).(1-t^{-2})\Big)$
of (\ref{exemple1}) does not contain the trivial character $\phi_{0}$
if $\Theta(E)=\sum_{n\in \Z}a_{n}t^n$ with
\begin{equation}\label{plus.grand.3}
a_{n}\neq 0\Longrightarrow n\geq 3\ .
\end{equation}
(\ref{exemple1.1}) tells us that 
(\ref{plus.grand.3}) is satisfied if the weights for the
action of $S^1$ in the fibers of the complex vector bundle 
$E\vert_{M^{S^1}_{+}}\otimes\det\Ncal^{+}$ are all 
bigger than $3$. 

\medskip

The conditions are weaker in the `Hamiltonian' situation. The term 
\break $\HolS(\widetilde{\Theta}(E))$
of (\ref{exemple2}) does not contain the trivial character $\phi_{0}$
if $\widetilde{\Theta}(E)=\sum_{n\in \Z}a_{n}t^n$ with
\begin{equation}\label{plus.grand.1}
a_{n}\neq 0\Longrightarrow n\geq 1\ ,
\end{equation}
and this condition is fulfilled if the weights for the
action of $S^1$ in the fibers of the complex vector bundle 
$E\vert_{M^{S^1}_{+}}\otimes\det\widetilde{\Ncal}^{+}$ 
are all bigger than $1$. Here we have another important difference: 
the vector bundle $\widetilde{\Ncal}^{+}\to M^{S^1}_{+}$ is not equal to 
the zero bundle if $0\in\mu(M)$ (see Lemma \ref{lem.N.beta.plus}).

We see finally that, in the Hamiltonian case, the condition 
`{\em $E$ is $\mu$-positive}' implies 
$$
0\in\mu(M)\Longrightarrow \left[ RR^{^{G}}(M,E)\right]^{G}=
\left[ RR^{^{G}}_{0}(M,E)\right]^G\ .
$$


\section{Appendix B: Induction map and multiplicities}

Let $G$ be a compact connected Lie group, with maximal torus $H$, and 
$\hgot_{+}^{*}\subset \hgot^{*}=(\ggot^{*})^{H}$ some choice of  
positive Weyl chamber. We denote by $\Rgot_{+}$ the associated system of 
positive roots, and we label the irreducible representations of $G$  
by the set $\Lambda^{*}_{+}=\Lambda^{*}\cap \hgot^{*}_{+}$ of dominant 
weights. For any weights $\alpha\in \Lambda^{*}$ we denote by 
$H\to \C^{*},\ h\mapsto h^{\alpha}$ the corresponding 
character : $(\exp(X))^{\alpha}=e^{\imath\langle \alpha,X\rangle}$ 
for $X\in \hgot$.

Let $W$ be the Weyl group of $(G,H)$, and $\carre(H)$ be the vector space of 
square integrable complex functions on $H$. For $f\in\carre(H)$, we consider
$J(f)$ \break $=\sum_{w\in W}(-1)^{w}\, w.f$, 
where $W\to \{1,-1\},\ w\to (-1)^{w}$, is the signature operator and 
$w.f\in \carre(H)$ is defined by $w.f(h)=f(w^{-1}.h),\ h\in H$ 
(see Section 7.4 of \cite{Bourbaki.Lie.9}). The map
$\frac{1}{\vert W\vert}J$ is the orthogonal projection from $\carre(H)$
to the space of $W$-anti-invariant elements of $\carre(H)$.

Let $\rho\in\hgot^{*}$ be the half sum of the positive roots. The 
function $H\to \C^{*},\ h\mapsto h^{\rho}$ is well defined as an 
element of $\carre(H)$ (even if $\rho$ is not a weight). The Weyl's character 
formula can be written in the following way. For any dominant weight 
$\lambda\in \Lambda^{*}_{+}$, the restriction 
$\chi_{_{\lambda}}^{_{G}}\vert_{H}$ 
of the irreducible $G$-character $\chi_{_{\lambda}}^{_G}$ satisfies 
\begin{equation}\label{eq.Weyl.1}
J(h^{\rho}).\chi_{_{\lambda}}^{_G}\vert_{H}= J(h^{\lambda+\rho})\quad {\rm 
in}\quad \carre(H)\ .
\end{equation}

For our purpose we give an expression of the character $\chi_{_{\lambda}}^{_G}$ 
through the induction map $\indH:\fgene(H)\to\fgene(G)^{G}$ (see 
(\ref{eq:fonction.induction})). 
Consider the affine action of the Weyl group on 
the set of weights : $w\circ \lambda =w.(\lambda+\rho)-\rho$ for $w\in 
W$ and $\lambda\in \Lambda^{*}$.

\begin{lem}\label{lem.Weyl}
   1) For any dominant weight $\lambda\in \Lambda^{*}_{+}$, the 
    character $\chi_{_{\lambda}}^{_G}$ is determined by the relation
    $\chi_{_{\lambda}}^{_G}= \indH\Big(h^{\lambda}\prod_{\alpha\in\Rgot_{+}}
    (1-h^{\alpha})\Big)$ in $\fgene(G)^{G}$.

    \noindent 2) For $\lambda\in\Lambda^{*}$ and $w\in W$, we have 
    $\indH(h^{w\circ\lambda}\Pi_{\alpha\in\Rgot_{+}}(1-h^{\alpha}))=$ \break
    $(-1)^{w}\indH(h^{\lambda}\Pi_{\alpha\in\Rgot_{+}}(1-h^{\alpha}))$.
    
    \noindent 3) For any weight $\lambda$, the following statements are 
    equivalent :
    
    a) $\ \indH(h^{\lambda}\Pi_{\alpha\in\Rgot_{+}}(1-h^{\alpha}))=0$,
    
    b) $\ W\circ \lambda\cap \Lambda^{*}_{+}=\emptyset$,
    
    c) The element $\lambda+\rho$ is not a regular element of $\hgot^{*}$.
\end{lem}

{\em Proof of 1) :}To prove it, we need the 
following relations proved in \cite{Bourbaki.Lie.9}[section 7.4] :

i) $\ \overline{J(h^{\rho})}=h^{-\rho}\prod_{\alpha\in\Rgot_{+}}(1-h^{\alpha})$,
\quad 
ii) $\ J(h^{\rho}).\overline{J(h^{\rho})}=\prod_{\alpha\in\Rgot}(1-h^{\alpha})$.

Let $dg$ and $dt$ be respectively the normalized Haar measures on $G$ and $H$. 
For any $f\in\f(G)^{G}$ we have
\begin{eqnarray*}
    \int_{G}\chi^{_G}_{_{\lambda}}(g)\, f(g)\, dg &=&
    \frac{1}{\vert W\vert}\int_{H}\chi_{_{\lambda}}^{_G}\vert_{H}(h)\, 
    \Pi_{\alpha\in\Rgot}(1-h^{\alpha})\, f\vert_{H}(h)\, dh
     \hspace{2cm} [1]       \\
    &=&
    \frac{1}{\vert W\vert}\int_{H}J(h^{\lambda+\rho})\, 
    \overline{J(h^{\rho})}\, f\vert_{H}(h)\, dh
    \hspace{3cm} [2] \\
    &=&
    \int_{H}h^{\lambda+\rho}\,\overline{J(h^{\rho})}\, f\vert_{H}(h)\,dh
    \hspace{4cm} [3] \\
    &=&
    \int_{H} h^{\lambda}\,\Pi_{\alpha\in\Rgot_{+}}(1-h^{\alpha})\,
    f\vert_{H}(h)\, dh \ .
    \hspace{3cm} [4] 
\end{eqnarray*}
The first equality is the Weyl integration formula.  The equality 
$[2]$ comes from ii) and (\ref{eq.Weyl.1}).  Since $\frac{1}{\vert 
W\vert} J$ is the orthogonal projection on 
$\carre(H)^{W-anti-invariant}$ and $h\mapsto \overline{J(h^{\rho})}\, 
f\vert_{H}(h)$ is $W$-anti-invariant we obtain the third equality.  
The equality $[4]$ comes from i).

{\em Proof of 2) : }  From $i)$, wee see that 
$h^{w\circ\lambda}\Pi_{\alpha\in\Rgot_{+}}(1-h^{\alpha})=
h^{w(\lambda+\rho)}\overline{J(h^{\rho})}=$ \break 
$(-1)^{w}\, w^{-1}.( h^{\lambda+\rho}\overline{J(h^{\rho})})= (-1)^{w}\, 
w^{-1}.(h^{\lambda}\Pi_{\alpha\in\Rgot_{+}}(1-h^{\alpha}))$, hence 
the relation $2)$ is proved since $\indH$ is $W$-invariant.

{\em Proof of 3) :} The implication $a)\Longrightarrow b)$ is an 
immediate consequence of $1)$ and $2)$. Proposition 3 in section 7.4 of 
\cite{Bourbaki.Lie.9} tells us that $\{J(h^{\lambda'+\rho}),\ 
\lambda'\in \Lambda^{*}_{+}\}$ is an orthogonal basis of the Hilbert 
space $\carre(H)^{W-anti-invariant}$. For $\lambda\in \Lambda^{*}$ and 
$\lambda'\in\Lambda^{*}_{+}$ we have 
$<J(h^{\lambda+\rho}),J(h^{\lambda'+\rho})>_{\carre}=$ 
$\vert W\vert<J(h^{\lambda+\rho}),h^{\lambda'+\rho}>_{\carre}=
\vert W\vert\sum_{w\in W}(-1)^w\int_{T}t^{w\circ\lambda-\lambda'}dt$. 
Thus,  the condition 
$W\circ \lambda\cap \Lambda^{*}_{+}=\emptyset$ is equivalent to  
$J(h^{\lambda+\rho})=0$. But the equality $[2]$ gives
$\indH(h^{\lambda}\Pi_{\alpha\in\Rgot_{+}}(1-h^{\alpha}))=$ \break
$\frac{1}{\vert W\vert}
\indH(J(h^{\lambda+\rho})h^{-\rho}\Pi_{\alpha\in\Rgot_{+}}(1-h^{\alpha}))$, 
hence $J(h^{\lambda+\rho})=0$ implies the point $a)$. We have proved 
that $b)\Longrightarrow a)$. Finally we see that 
$J(h^{\lambda+\rho})=0 \Longleftrightarrow \exists w\in W, 
w.(\lambda+\rho)=\lambda+\rho \Longleftrightarrow \lambda+\rho\ {\rm 
is\ not\ a\ regular\ value\ of\ } \hgot^{*}$. We have proved 
that $b)\Longleftrightarrow c)$. $\Box$

\bigskip

From the previous Lemma, we see that $v\mapsto
\indH(v(h)\Pi_{\alpha\in\Rgot_{+}}(1-h^{\alpha}))$ is the holomorphic 
induction map 
\begin{equation}\label{eq.holomorphe.G.H}
    \HolH:R(H)\to R(G)\ .
\end{equation}    
We keep the same notation for the extended map 
$\HolH:R^{-\infty}(H)\to R^{-\infty}(G)$. Note that the choice of a 
positive Weyl chamber $\hgot^{*}_{+}$ determines a complex structure on
$\ggot/\hgot$, and $\Pi_{\alpha\in\Rgot_{+}}(1-h^{\alpha})$ is the 
trace of the virtual $H$-representation 
$\wedge_{\C}^{\bullet}\ggot/\hgot\, \in R(H)$. Then the map $\HolH$ 
will be defined simply by the relation 
$\HolH(v)=\indH(v\,\wedge_{\C}^{\bullet}\ggot/\hgot)$.

\begin{rem}\label{wedge.C-wedge.R}
The relations i) and ii)  used in the proof of the past lemma show  
that $\sum_{w\in W}w. \prod_{\alpha>0}(1-h^{\alpha})= 
\sum_{w\in W}w.(\overline{J(h^{\rho})}h^{\rho})=
\overline{J(h^{\rho})}.J(h^{\rho})=
\prod_{\alpha}(1-h^{\alpha})$. In other words 
$\sum_{w\in W}w.\wedge_{\C}^{\bullet}\ggot/\hgot=
(\wedge_{\R}^{\bullet}\ggot/\hgot)\otimes\C=
\wedge_{\C}^{\bullet}\ggot/\hgot\,
\wedge_{\C}^{\bullet}\overline{\ggot/\hgot}$ 
in $R(H)$. These equalities give 
\begin{equation}\label{eq.Weyl.Hol}
  \indH\Big((\sum_w w.\phi)\, \wedge_{\C}^{\bullet}\ggot/\hgot\Big)=
  \indH(\phi\,\wedge_{\R}^{\bullet}\ggot/\hgot)
\end{equation}
since $\indH$ is $W$-invariant. The Weyl integration formula is usually
  state as the relation $f=\frac{1}{\vert W\vert}
\indH(f\vert_{H}\,\wedge_{\R}^{\bullet}\ggot/\hgot)$ for any
$f\in \f(G)^G$. But $f\vert_{H}$ is $W$-invariant, so 
(\ref{eq.Weyl.Hol}) gives   
$\frac{1}{ \vert W\vert}
\indH(f\vert_{H}\,\wedge_{\R}^{\bullet}\ggot/\hgot)=
\indH(f\vert_{H}\,\wedge_{\C}^{\bullet}\ggot/\hgot)$.
Finally, for any $\phi\in R(G)$, the Weyl integration formula is
equivalent to the following equality in $R(G)$:
$$
\phi=\HolH(\phi\vert_{H}) .
$$
\end{rem}

\begin{rem}\label{hol.egale.1}
A weight $\lambda$ satisfies $\HolH(h^{\lambda})=\pm 1$ if and only if
$0\in\,W\circ \lambda\cap\Lambda^*_{+}$, that is $\lambda=
-(\rho -w.\rho)$ for some $w\in W$. But a small computation shows that 
$\rho -w.\rho=\sum_{\alpha>0,w^{-1}.\alpha<0}\alpha$, hence 
$\langle \rho -w.\rho,X\rangle\geq 0$ for any $X\in\hgot_{+}$. 
Finally the equality $\HolH(h^{\lambda})=\pm 1$ implies that
$\langle \lambda,X\rangle\leq 0$ for any $X\in\hgot_{+}$.
\end{rem}

\medskip

Consider now the stabiliser $G_{\beta}$ of the non-zero element 
$\beta\in\hgot_{+}$. The subgroup $H$ is also a maximal torus of $G_{\beta}$. 
The Weyl group $W_{\beta}$ of $(G_{\beta},H)$ is identified with 
$\{w\in W,\ w.\beta=\beta\}$. We consider a Weyl chamber 
$\hgot^{*}_{+,\beta}\subset\hgot^{*}$ for $G_{\beta}$ that contains 
the Weyl chamber $\hgot^{*}_{+}$ of $G$. The irreducible representations
$\chi_{_{\lambda}}^{_{G_{\beta}}},\ \lambda\in\Lambda_{+,\beta}^{*}$ of
$G_{\beta}$ are labelled by the set $\Lambda_{+,\beta}^{*}=
\Lambda^{*}\cap\hgot_{+,\beta}^{*}$ of dominant weights.

We have a unique `holomorphic' induction map $\HolB:R(G_{\beta})\to 
R(G)$ such that $\HolH=\HolB\circ\HolHB$. This map is defined 
precisely by the equation\footnote{We take on $\ggot/\ggot_{\beta}$ 
the complex structure defined by $\beta$.} 
\begin{equation}\label{eq:Hol-G-beta}
 \HolB(v)=\indB\left(v\,\wedge_{\C}^{\bullet}\ggot/\ggot_{\beta}\right)\ ,
\end{equation}
 for every $v\in R(G_{\beta})$.

\bigskip

We finish this appendix with some general remarks about 
$P$-transversally elliptic symbols on a compact manifold $M$, when a 
subgroup $\tore$ in the center of $P$ acts trivially on $M$.

More precisely, let $H$ be a compact maximal torus in $P$, $\hgot_{+}$ 
be a choice of a positive Weyl chamber in the Lie algebra $\hgot$ of 
$H$, and let $\beta\in \hgot_{+}$ be a non-zero element in the center 
of the Lie algebra $\pgot$ of $P$\footnote{The Lie group $P$ is 
supposed connected then $\beta\in(\pgot)^{P}$.}.  We suppose here that 
the subtorus $\tore\subset H$, which is equal to the closure of 
$\{\exp(t.\beta),\ t\in\R\}\ $, acts trivially on $M$.

Every $P$-equivariant complex vector bundle $E\to M$ can be decomposed 
relatively to the $\tore$-action: 
$E=\oplus_{a\in\hat{\tore}}E^a\otimes\C_{a}$, where 
$E^{a}:=\hom_{\tore}(E,\C^*_{a})$\footnote{The torus $\tore$ acts on 
the complex line $\C_{a}$ with the representation $t\to t^a$.} is a 
$P$-complex vector bundle with a trivial action of $\tore$.  Then, 
each $P$-equivariant symbol $\sigma:p^{*}(E_{1})\to p^{*}(E_{2})$ 
where $E_{1},E_{2}$ are $P$-equivariant complex vector bundles over 
$M$, and where $p:\T M\to M$ is the canonical projection, admits a 
finite $P\times\tore$-equivariant decomposition
\begin{equation}\label{eq:sigma-a}
    \sigma=\sum_{a\in\hat{\tore}}\sigma^a\otimes\C_{a}.
\end{equation}
Here $\sigma^a:p^{*}(E^a_{1})\to p^{*}(E^a_{2})$ is a $P$-equivariant 
symbol, trivial for the $\tore$-action.

Let us consider the inclusion map $i:\tore\croc H$, 
with the induced maps $i:Lie(\tore) \to \hgot$ at the level of Lie algebra and
$i^{*}:\hgot^{*}\to Lie(\tore) ^{*}$. Note that $i^{*}(\lambda)$ is a 
weight for $\tore$ if $\lambda$ is a weight for $H$.

\begin{lem}\label{lem.multiplicites.tore}
   Let $M$ be a $P$-manifold with the same properties as above.
   Let $\sigma:p^{*}(E_{1})\to p^{*}(E_{2})$ be a $P$-equivariant
   transversally elliptic symbol over $M$ and denote by 
   $m_{\lambda}(\sigma),\, \lambda\in \Lambda^{*}_{P,+}$, the
   multiplicities of its index : $\indice_{M}^P(\sigma)=
   \sum_{\lambda\in \Lambda^{*}_{P,+}}
   m_{\lambda}(\sigma)\chi^{_P}_{_{\lambda}}$. Then, if 
   $m_{\lambda}(\sigma)\neq 0$, the weight $a=i^{*}(\lambda)$
   occurs in the decomposition (\ref{eq:sigma-a}).
\end{lem}

\begin{coro}\label{coro.multiplicites.tore}   
   Suppose that the weights $a\in \hat{\tore}$ which occur in the 
   decomposition (\ref{eq:sigma-a}) satisfy 
   $\langle a,\beta\rangle\geq \eta$ for some fixed $\eta\in\R$. 
   Then, for the multiplicities, we get 
   $$
   m_{\lambda}(\sigma)\neq 0\Longrightarrow 
   \langle \lambda,\beta\rangle\geq \eta\ .
   $$
   In particular, $\indice_{M}^P(\sigma)$ does not contain the 
   trivial representation when $\eta>0$. 
\end{coro}   
   
\begin{rem}The previous Lemma and Corollary remain true if $M$ is a 
    $P$-invariant open subset of a compact $P$-manifold.
\end{rem}

 For the Corollary, we have just to notice 
that\footnote{We use the same notations for 
$\beta\in Lie(\tore) $ and $i(\beta)\in\hgot$.}
$\langle\lambda,\beta\rangle
=\langle a,\beta\rangle$ for 
$a=i^{*}(\lambda)$. Then, if we have 
$\langle a,\beta\rangle\geq \eta$ for all $\tore$-weights  
occurring in $\sigma$, we get $\langle\lambda,\beta\rangle\geq \eta $
for every $\lambda$ such that $m_{\lambda}(\sigma)\neq 0$.

\medskip

{\em Proof of Lemma \ref{lem.multiplicites.tore}}: 
Let $P'$ be a Lie subgroup of $P$ such that 
$r:\tore\times P'\to P,\ r(t,g)=t.g $, is a finite covering of $P$.
The map $r$ induces $r^{*}:K_{P}(\T_{P}M)\to K_{\tore\times 
P'}(\T_{P'}M)$\footnote{Note that $\T_{P'}M=\T_{P}M$ because $\tore$ 
acts trivially on $M$.} and an injective map $r^{*}:R^{-\infty}(P)\to
R^{-\infty}(\tore\times P')$, such that 
$\indice_{M}^{\tore\times P'}(r^{*}\sigma)=r^{*}(\indice_{M}^P(\sigma))$.

The decomposition (\ref{eq:sigma-a}) can be read through the 
identification $K_{\tore\times P'}(\T_{P'}M)=$\break
$K_{P'}(\T_{P'}M)\otimes R(\tore)$: we have 
$r^*\sigma=\sum_{a\in\hat{\tore}}\sigma^a\otimes\C_{a}$
with $\sigma^a\in  K_{P'}(\T_{P'}M)$. Hence 
\begin{equation}\label{decomposition-1}
\indice_{M}^{\tore\times P'}(r^{*}\sigma)(t,g)=\sum_{a\in\hat{\tore}}
\indice_{M}^{P'}(\sigma^a)(g).\, t^a\ ,\quad (t,g)\in \tore\times P'\ .
\end{equation}
The irreducible characters $\chi^{_P}_{_{\lambda}}$ satisfy
$r^{*}\chi^{_P}_{_{\lambda}}(t,g)=\chi^{_P}_{_{\lambda}}\vert_{P'}(g).\, 
t^{i^{*}(\lambda)}$. If we start from the 
decomposition  $\indice_{M}^P(\sigma)=\sum_{\lambda\in \Lambda^{*}_{P,+}}
m_{\lambda}(\sigma)\chi^{_P}_{_{\lambda}}$ relative to the 
irreducible characters of $P$, we get
\begin{equation}\label{decomposition-2}
r^{*}\left(\indice_{M}^{\tore\times P'}(\sigma)\right)(t,g)=
\sum_{a\in\hat{\tore}}\left(\sum_{i^{*}(\lambda)=a}
m_{\lambda}(\sigma)\chi^{_P}_{_{\lambda}}\vert_{P'}(g)\right) .\, 
t^a\ , 
\end{equation}
for any $(t,g)\in \tore\times P'$.
If we compare (\ref{decomposition-1}) and 
(\ref{decomposition-2}), we get $\indice_{M}^{P'}(\sigma^a)=
\sum_{i^{*}(\lambda)=a}m_{\lambda}(\sigma)
\chi^{_P}_{_{\lambda}}\vert_{P'}$. The map $r^{*}:R^{-\infty}(P)\to
R^{-\infty}(\tore\times P')$ is injective, so $\sum_{i^{*}(\lambda)=a}
m_{s\lambda}(\sigma)\chi^{_P}_{_{\lambda}}\vert_{P'}=0$ if and only 
if $m_{\lambda}(\sigma)=0$ for every $\lambda$ satisfying 
$i^{*}(\lambda)=a$. Hence if the multiplicity 
$m_{\lambda}(\sigma)$ is non zero, the element $a=i^{*}(\lambda)$
is a weight for the action of $\tore$ on  
$\sigma:p^{*}(E_{1})\to p^{*}(E_{2})$. $\Box$

\bigskip



\end{document}